\documentclass[12pt,leqno]{amsart}
\usepackage{euscript}
\newtheorem{theorem}{Theorem}
\newtheorem{definition}[theorem]{Definition}

\newtheorem{example}[theorem]{Example}

\newtheorem{proposition}[theorem]{Proposition}

\def\gen#1#2{
\if #11
    \if #22 \jednadva \else \fi
\else
\fi
\if #12
    \if #22 \dvadva \else \fi
\else
\fi
\if #12
    \if #21 \dvajedna \else \fi
\else
\fi
\if #13
    \if #22 \tridva \else \fi
\else
\fi
\if #13
    \if #21 \trijedna \else \fi
\else
\fi
\if #12
    \if #23 \dvatri \else \fi
\else
\fi
\if #11
    \if #23 \jednatri \else \fi
\else
\fi
\if #11
    \if #24 \jednactyri \else \fi
\fi
\if #14
    \if #21 \ctyrijedna \else \fi
\fi
}
\def\jednadva{{
\unitlength=.4pt
\begin{picture}(24.00,20.00)(-2.00,0.00)
\bezier{20}(10.00,10.00)(15.00,5.00)(20.00,0.00)
\bezier{20}(10.00,10.00)(5.00,5.00)(0.00,0.00)
\put(10.00,20.00){\line(0,-1){10.00}}
\end{picture}}
}

\def\dvajedna{{
\unitlength=.4pt
\begin{picture}(24.00,20.00)(-2.00,0.00)
\put(10.00,10.00){\line(0,-1){10.00}}
\bezier{20}(10.00,10.00)(15.00,15.00)(20.00,20.00)
\bezier{20}(0.00,20.00)(5.00,15.00)(10.00,10.00)
\end{picture}}
}

\def\bZbbZ{
{
\unitlength=.27pt
\begin{picture}(48.00,30.00)(-4,0.00)
\bezier{34}(20.00,20.00)(30.00,10.00)(40.00,0.00)
\bezier{34}(20.00,20.00)(10.00,10.00)(0.00,0.00)
\bezier{20}(30.00,10.00)(25.00,5.00)(20.00,0.00)
\put(20.00,30.00){\line(0,-1){10.00}}
\end{picture}}
}

\def\ZvvZv{{
\unitlength=.27pt
\begin{picture}(48.00,30.00)(-4,-30.00)
\bezier{30}(20.00,-20.00)(30.00,-10.00)(40.00,0.00)
\bezier{30}(20.00,-20.00)(10.00,-10.00)(0.00,0.00)
\bezier{20}(10.00,-10.00)(15.00,-5.00)(20.00,0.00)
\put(20.00,-30.00){\line(0,1){10.00}}
\end{picture}}}

\def\vZvvZ{{
\unitlength=.27pt
\begin{picture}(48.00,30.00)(-4,-30.00)
\bezier{34}(20.00,-20.00)(30.00,-10.00)(40.00,0.00)
\bezier{34}(20.00,-20.00)(10.00,-10.00)(0.00,0.00)
\bezier{20}(30.00,-10.00)(25.00,-5.00)(20.00,0.00)
\put(20.00,-30.00){\line(0,1){10.00}}
\end{picture}}
}

\def\dvojiteypsilon{{
\unitlength=.3pt
\begin{picture}(24.00,30.00)(0.00,3.00)
\put(10.00,20.00){\line(0,-1){10.00}}
\bezier{20}(10.00,10.00)(15,5)(20.00,0.00)
\bezier{20}(10.00,10.00)(5,5)(0.00,0.00)
\bezier{20}(10.00,20.00)(15,25)(20.00,30.00)
\bezier{20}(0.00,30.00)(5,25)(10.00,20.00)
\end{picture}}}

\def\Jac#1#2#3{{
\unitlength=.4pt
\begin{picture}(48.00,30.00)(-4,0.00)
\bezier{34}(20.00,20.00)(30.00,10.00)(40.00,0.00)
\bezier{34}(20.00,20.00)(10.00,10.00)(0.00,0.00)
\bezier{20}(10.00,10.00)(15.00,5.00)(20.00,0.00)
\put(20.00,30.00){\line(0,-1){10.00}}
\put(0,-5){\makebox(0,0)[t]{\scriptsize $#1$}}
\put(20,-5){\makebox(0,0)[t]{\scriptsize $#2$}}
\put(40,-5){\makebox(0,0)[t]{\scriptsize $#3$}}
\end{picture}}}

\def\coJac#1#2#3{{
\unitlength=.4pt
\begin{picture}(48.00,30.00)(-4,-30.00)
\bezier{30}(20.00,-20.00)(30.00,-10.00)(40.00,0.00)
\bezier{30}(20.00,-20.00)(10.00,-10.00)(0.00,0.00)
\bezier{20}(10.00,-10.00)(15.00,-5.00)(20.00,0.00)
\put(20.00,-30.00){\line(0,1){10.00}}
\put(0,5){\makebox(0,0)[b]{\scriptsize $#1$}}
\put(20,5){\makebox(0,0)[b]{\scriptsize $#2$}}
\put(40,5){\makebox(0,0)[b]{\scriptsize $#3$}}
\end{picture}}}

\def\dvojiteypsilonvetsi#1#2#3#4{{
\unitlength=.46pt
\begin{picture}(24.00,30.00)(0.00,3.00)
\put(10.00,20.00){\line(0,-1){10.00}}
\bezier{20}(10.00,10.00)(15,5)(20.00,0.00)
\bezier{20}(10.00,10.00)(5,5)(0.00,0.00)
\bezier{20}(10.00,20.00)(15,25)(20.00,30.00)
\bezier{20}(0.00,30.00)(5,25)(10.00,20.00)
\put(0,-5){\makebox(0,0)[t]{\scriptsize $#1$}}
\put(20,-5){\makebox(0,0)[t]{\scriptsize $#2$}}
\put(0,35){\makebox(0,0)[b]{\scriptsize $#3$}}
\put(20,35){\makebox(0,0)[b]{\scriptsize $#4$}}
\end{picture}}}

\def\pravaplastevvetsi#1#2#3#4{
\unitlength=.45pt
\begin{picture}(38,30)(-4,4)
\put(0.00,0.00){\line(0,1){10}}
\put(20.00,0.00){\line(0,1){10}}
\put(10.00,20.00){\line(0,1){10}}
\put(30.00,20.00){\line(0,1){10}}
\put(0,10){\bezier{20}(0.00,0.00)(5.00,5.00)(10.00,10.00)}
\put(20,10){\bezier{20}(0.00,0.00)(5.00,5.00)(10.00,10.00)}
\put(20,10){\bezier{20}(0.00,0.00)(-5.00,5.00)(-10.00,10.00)}
\put(0,-5){\makebox(0,0)[t]{\scriptsize $#1$}}
\put(20,-5){\makebox(0,0)[t]{\scriptsize $#2$}}
\put(10,35){\makebox(0,0)[b]{\scriptsize $#3$}}
\put(30,35){\makebox(0,0)[b]{\scriptsize $#4$}}
\end{picture}
}

\def\levaplastevvetsi#1#2#3#4{
\unitlength=.45pt
\begin{picture}(38,30)(-34,4)
\put(0.00,0.00){\line(0,1){10}}
\put(-20.00,0.00){\line(0,1){10}}
\put(-10.00,20.00){\line(0,1){10}}
\put(-30.00,20.00){\line(0,1){10}}
\put(0,10){\bezier{20}(0.00,0.00)(-5.00,5.00)(-10.00,10.00)}
\put(-20,10){\bezier{20}(0.00,0.00)(-5.00,5.00)(-10.00,10.00)}
\put(-20,10){\bezier{20}(0.00,0.00)(5.00,5.00)(10.00,10.00)}
\put(0,-5){\makebox(0,0)[t]{\scriptsize $#2$}}
\put(-20,-5){\makebox(0,0)[t]{\scriptsize $#1$}}
\put(-10,35){\makebox(0,0)[b]{\scriptsize $#4$}}
\put(-30,35){\makebox(0,0)[b]{\scriptsize $#3$}}
\end{picture}
}

\def\hPROP{$\frac12${\small\sc PROP}}

\def\sfLieB{{\sf LieB}}

\def\motylek{{
\unitlength=.3pt
\begin{picture}(66.00,60.00)(-3.00,20.00)
\put(10.00,50.00){\line(0,1){0.00}}
\put(50.00,10.00){\line(0,-1){10.00}}
\put(10.00,10.00){\line(0,-1){10.00}}
\put(50.00,60.00){\line(0,-1){10.00}}
\put(10.00,60.00){\line(0,-1){10.00}}
\put(60.00,40.00){\line(0,-1){20.00}}
\put(0.00,40.00){\line(0,-1){20.00}}
\bezier{20}(50.00,10.00)(55.00,15.00)(60.00,20.00)
\bezier{20}(0.00,20.00)(5.00,15.00)(10.00,10.00)
\bezier{20}(50.00,50.00)(55.00,45.00)(60.00,40.00)
\bezier{20}(10.00,10.00)(15,15)(22,22)
\put(28,28){\bezier{20}(10.00,10.00)(15,15)(22,22)}
\put(10.00,50.00){\line(1,-1){40.00}}
\bezier{20}(10.00,50.00)(5.00,45.00)(0.00,40.00)
\end{picture}}}

\def\bZbbZ{
{
\unitlength=.27pt
\begin{picture}(48.00,30.00)(-4,0.00)
\bezier{34}(20.00,20.00)(30.00,10.00)(40.00,0.00)
\bezier{34}(20.00,20.00)(10.00,10.00)(0.00,0.00)
\bezier{20}(30.00,10.00)(25.00,5.00)(20.00,0.00)
\put(20.00,30.00){\line(0,-1){10.00}}
\end{picture}}
}

\def\ZbbZb{{
\unitlength=.27pt
\begin{picture}(48.00,30.00)(-4,0.00)
\bezier{34}(20.00,20.00)(30.00,10.00)(40.00,0.00)
\bezier{34}(20.00,20.00)(10.00,10.00)(0.00,0.00)
\bezier{20}(10.00,10.00)(15.00,5.00)(20.00,0.00)
\put(20.00,30.00){\line(0,-1){10.00}}
\end{picture}}}

\def\ef{
\unitlength 1mm 
\linethickness{0.4pt}
\begin{picture}(2,2)(16,16.5)
\put(17,20){\line(0,-1){4}}
\put(17,18){\makebox(0,0)[cc]{\small $\bullet$}}
\end{picture}
}

\def\glab#1#2{
\unitlength 1mm
\linethickness{0.4pt}
\begin{picture}(8,4)(15,15)
\put(19,20){\line(0,-1){2}}
\multiput(19,18)(-.0333333,-.0333333){60}{\line(0,-1){.0333333}}
\multiput(19,18)(.0333333,-.0333333){60}{\line(0,-1){.0333333}}
\put(17,14){\makebox(0,0)[cc]{\scriptsize $#1$}}
\put(21,14){\makebox(0,0)[cc]{\scriptsize $#2$}}
\end{picture}
}

\def\levaplastev{
\unitlength=.35pt
\begin{picture}(38,30)(-34,4)
\put(0.00,0.00){\line(0,1){10}}
\put(-20.00,0.00){\line(0,1){10}}
\put(-10.00,20.00){\line(0,1){10}}
\put(-30.00,20.00){\line(0,1){10}}
\put(0,10){\bezier{20}(0.00,0.00)(-5.00,5.00)(-10.00,10.00)}
\put(-20,10){\bezier{20}(0.00,0.00)(-5.00,5.00)(-10.00,10.00)}
\put(-20,10){\bezier{20}(0.00,0.00)(5.00,5.00)(10.00,10.00)}
\end{picture}
}

\def\pravaplastev{
\unitlength=.35pt
\begin{picture}(38,30)(-4,4)
\put(0.00,0.00){\line(0,1){10}}
\put(20.00,0.00){\line(0,1){10}}
\put(10.00,20.00){\line(0,1){10}}
\put(30.00,20.00){\line(0,1){10}}
\put(0,10){\bezier{20}(0.00,0.00)(5.00,5.00)(10.00,10.00)}
\put(20,10){\bezier{20}(0.00,0.00)(5.00,5.00)(10.00,10.00)}
\put(20,10){\bezier{20}(0.00,0.00)(-5.00,5.00)(-10.00,10.00)}
\end{picture}
}

\def\longhookleftarrow{
\begin{picture}(6,2)
\unitlength.7em
\put(0.00,0.00){\vector(-1,0){5}}
\put(0.00,0.00){\bezier{50}(0,0)(.5,0)(.5,.5)}
\put(.5,.5){\bezier{50}(0,0)(0,.5)(-.5,.5)}
\end{picture}
}

\def\cathPROP{{\textstyle\frac12}{\tt PROP}}
\def\catOper{{\tt Oper}} \def\catProper{\tt Proper}
\def\catDiop{{\tt diOp}}
\def\catProper{\tt Proper}
\def\catPROP{{\tt PROP}}

  \def\sfR{{\sf R}}
\def\sfr{{\sf r}} \def\sfb{{\sf b}}
\def\freehPROP{{\Gamma_{\frac12}}} \def\subc{{\tt SGr}}

\def\modcoll#1{\{#1(g,n)\}_{g\geq 0,n\geq -1}}
\def\n#1{\{0,\ldots,#1\}} \def\sfP{{\sf P}}

\def\Graphcat#1#2{{\tt MGr}(#1,#2)} \def\MGr{\Graphcat}
\def\GR{{\tt Gr}}\def\UGR{\tt UGr}
\def\Gr#1#2{{\GR}(#1,#2)}\def\UGr#1#2{{\UGR}(#1,#2)}
\def\GRH{{\tt Gr}_{\frac12}}
\def\Grh#1#2{{\GRH(#1,#2)}}
\def\UGRC{{{\tt UGr}_c}}      \def\UGRD{{\tt UGr}_{\tt D}}
\def\UGrc#1#2{{\UGRC(#1,#2)}} \def\UGrd#1#2{\UGRD(#1,#2)}
 
\def\Modtriple{{\mathbb M}\hskip .1em} 

\def\collmod{\modcoll}\def\MCol{{\tt MMod}}\def\sfI{{\sf I}}
\def\MMod{\MCol}\def\sfQ{{\sf Q}}
\newcommand{\Bij}{{\it Bij}}
\def\lzavorkaobyc{(}\def\rzavorkaobyc{)} 
\def\hatcalM{{\widehat {\mathfrak M}}} \def\calE{{\mathcal E}}
\def\stablecoll#1{\{#1(g,n)\}_{(g,n) \in \frakS}} \def\frakS{{\mathfrak S}}

\def\Sigmap{\Sigma^+} \def\Hom#1#2{{\it Lin\/}(#1,#2)}
\def\calT{{\EuScript T}}\def\edge{{\it edge}} \def\leg{{\it Leg}}
 \def\Flag{{\it Flag}}
\def\pre-Lie{{\mbox {\it pre-\Lie}}} \def\Leib{\hbox{${\mathcal L}${\it eib}}}
\def\Span{{\it Span}}\def\eucalC{{\EuScript C}} \def\sgn{{\rm sgn}}
\def\lrubercurly#1pt{{\left\{\rule{0pt}{#1pt}\right.}}
\def\rrubercurly#1pt{{\left.\rule{0pt}{#1pt}\right\}}}

\def\bbjedna{{\sf 1}}\def\freePROP{{\Gamma_{\tt P}}}
\def\freeproperad{{\Gamma_{\tt c}}}\def\freedioperad{{\Gamma_{\tt D}}}

\def\enil(#1,#2){\line(#1,-#2)} \def\Lie{\hbox{{$\mathcal L$}{\it ie\/}}}
 \def\sfB{{\sf B}}
\def\Poiss{\mathcal{P}{\hskip -.15em \it oiss}}
\def\sLie{\hbox{{\scriptsize $\mathcal L$}{\scriptsize \it ie\/}}}
 \def\tup(#1,#2){\put(#1,-#2)}
\def\Ass{\hbox{{$\mathcal A$}{\it ss\/}}}\def\Aut{{\it Aut\/}}
\def\sAss{\hbox{{\scriptsize $\mathcal A$}{\scriptsize \it ss\/}}}
\def\sCom{\hbox{{\scriptsize $\mathcal C$}\hskip -.2mm{\scriptsize\it om\/}}}
\def\calC{{\mathcal C}}
\def\c#1{\circ_#1} \def\circij#1#2{{\mbox{${}_{#1}\circ_{#2}$}}}
\def\)){{)\hskip -.8mm)}}\def\(({{(\hskip -.8mm(}}
\def\Modk{{\tt Mod}_\bfk}
\def\tree{{\mathcal T}\hskip -.7mm{\it ree}}
\def\en{[n]}\def\TTT{{\tt Tree}}\def\RR{{\mathbb R}}
\def\UTTT{{\tt UTree}} \def\MTTT{{\tt MTree}}
\def\colim#1{\mathop{{\rm colim}}%
             \limits_{\rule{0em}{1em}\mbox{\scriptsize $#1$}}}

\def\bicol#1#2{\{#1(m,n)\}_{m,n \geq #2}}

\def\vlra{{\hbox{$-\hskip-1mm-\hskip-2mm\longrightarrow$}}}

\def\Set{{\tt Set}_f}
\def\In{{\it in}} \def\Out{{\it out}}
\def\Vert{{\it Vert}}\def\vert{{\it vert}}
\def\val{{\it val}}\def\ext{\it ext}
\def\DD{{\mathbb D}}\def\SS{{\mathbb S}}\def\wcalS{{\widetilde{\calS}}}
\def\intU{{\hskip 3.235pt\raisebox{10pt}{\scriptsize\rm o} \hskip -7pt U}}
\def\PROP{{\small PROP}}\def\Leg{{\it Leg}}
\def\Leaf{{\it Leaf}}\def\Edg{{\it Edg}}
\def\forget{
\setlength{\unitlength}{.5em}
\begin{picture}(1.2,1)
\put(0,0){\line(0,1){1}}
\put(0,0){\line(1,0){1}}
\put(1,1){\line(0,-1){1}}
\put(1,1){\line(-1,0){1}}
\end{picture}
}
\def\sigmamod{{\mbox{$\Sigma$-{\tt mod}}}}
\def\sigmabimod{{\mbox{$\Sigma$-{\tt bimod}}}}
\def\sigmamodp{{\mbox{$\Sigma^+$-{\tt mod}}}}
\def\ssigmamod{{\mbox{\scriptsize $\Sigma$-{\tt mod}}}}

\def\coll#1#2{{\{#1(n)\}_{n \geq #2}}}\def\May{{\mathrm {May}\/}}
\def\Mar{{\mathrm {Mar}\/}}
\def\calK{{\EuScript K}}\def\Mor{{\it Mor\/}}
\def\calV{{\EuScript V}}\def\calA{{\EuScript A}}
\def\ns{non-$\Sigma$}\def\Rada#1#2#3{#1_{#2},\dots,#1_{#3}}
\def\ucalP{{\underline{\EuScript P}}}\def\ugamma{\underline{\gamma}}
\def\ueta{{\underline{\eta}}}\def\barcalM{{\overline{\EuScript M}}}
\def\End{{\mathcal E} \hskip -.1em{\it nd\/}}\def\calS{{\EuScript S}}
\def\calP{{\EuScript P}} \def\ot{\otimes} \def\calQ{{\EuScript Q}}
\def\bfk{{\mathbf k}} \def\ot{\otimes} \def\calD{{\EuScript D}}
\def\id{{\it id}} \def\Lin{{\it Lin\/}} \def\calM{{\EuScript M}}
\def\Com{{\EuScript C}{\it om}}  \def\Associative{{\EuScript A}{\it ss}}
\def\uAss{{\underline{\Associative}}}\def\CPj{{\mathbb {CP}}^1}
\def\rada#1#2{{#1,\ldots,#2}}\def\calR{{\EuScript R}}\def\calI{{\EuScript I}}
\def\ccalP{{\widehat \calP}}
\def\barcalP{{\overline \calP}} \def\Ker{{\it Ker\/}}
\def\otexp#1#2{#1^{\otimes #2}} \def\rr{\rule{0em}{.9em}}
\def\cases#1#2#3#4{
                  \left\{
                         \begin{array}{ll}
                           #1,\ &\mbox{#2}
                           \\
                           #3,\ &\mbox{#4}
                          \end{array}
                   \right.
}
\def\bijections#1#2#3#4{{
{
\unitlength=.5pt
\thinlines
\begin{picture}(110.00,40.0)(0.00,13.00)
\put(49.50,0.00){\makebox(0.00,0.00)[t]{\scriptsize $#4$}}
\put(49.50,43.50){\makebox(0.00,0.00)[b]{\scriptsize $#3$}}
\put(97.00,20.00){\makebox(0.00,0.00)[lc]{$#2$}}
\put(5.00,20.00){\makebox(0.00,0.00)[rc]{$#1$}}
\put(20.00,8.5){\vector(-2,1){10.00}}
\put(76.00,32.00){\vector(2,-1){10.00}}
\bezier{100}(79.50,10.50)(49.50,0.50)(16.0,10.50)
\bezier{100}(16.0,30.50)(49.50,40.50)(79.50,30.50)
\end{picture}}
}}

\def\ovI{
\unitlength .37mm
\linethickness{0.4pt}
\ifx\plotpoint\undefined\newsavebox{\plotpoint}\fi 
\begin{picture}(13,40)(0,0)
\bezier{183}(5,36)(5,40)(9,40)
\bezier{183}(9,40)(13,40)(13,36)
\put(5,36){\line(0,-1){15}}
\put(13,36){\line(0,-1){15}}
\bezier{183}(5,21)(5,17)(9,17)
\bezier{183}(13,21)(13,17)(9,17)
\end{picture}
}

\def\ovII{
\unitlength .7mm
\linethickness{0.4pt}
\ifx\plotpoint\undefined\newsavebox{\plotpoint}\fi 
\begin{picture}(29,23)(0,0)
\bezier{130}(8,17)(6,19)(8,21)
\bezier{130}(8,21)(10,23)(12,21)
\put(8,17){\line(1,-1){15}}
\put(12,21){\line(1,-1){15}}
\bezier{130}(27,6)(29,4)(27,2)
\bezier{130}(23,2)(25,0)(27,2)
\end{picture}
}

\def\ovIV{
\unitlength .7mm
\thinlines
\begin{picture}(29,23)(0,0)
\bezier{130}(-8,17)(-6,19)(-8,21)
\bezier{130}(-8,21)(-10,23)(-12,21)
\put(-8,17){\line(-1,-1){15}}
\put(-12,21){\line(-1,-1){15}}
\bezier{130}(-27,6)(-29,4)(-27,2)
\bezier{130}(-23,2)(-25,0)(-27,2)
\end{picture}
}

\def\ovIIImensi{
\unitlength .43mm
\thinlines
\begin{picture}(75,66)(0,0)
\bezier{221}(34,62)(38,66)(42,66)
\bezier{221}(42,66)(46,66)(50,62)
\put(34,62){\line(-1,-1){19}}
\bezier{221}(15,43)(11,39)(11,35)
\bezier{183}(11,35)(11,28)(18,28)
\put(84,0){
\bezier{221}(-15,43)(-11,39)(-11,35)
\bezier{183}(-11,35)(-11,28)(-18,28)
}
\put(18,28){\line(1,0){50}}
\put(50,62){\line(1,-1){19}}
\end{picture}
}

\def\ovVtrochumensi{
\thinlines
\unitlength=1.3pt
\begin{picture}(48,66)(0,0)
\bezier{221}(34,62)(38,66)(42,66)
\put(34,62){\line(-1,-1){19}}
\bezier{221}(15,43)(11,39)(11,35)
\bezier{178}(42,66)(45.5,66)(47,62)
\bezier{185}(47,62)(48,59.5)(45,55)
\put(45,55){\line(-2,-3){21.33}}
\bezier{211}(23.5,23)(20,17)(15,18)
\put(11,35){\line(0,-1){12}}
\bezier{183}(11,25)(11,19)(15,18)
\end{picture}
}

\def\ovVmensi{
\thinlines
\unitlength=1.2pt
\begin{picture}(48,66)(0,0)
\bezier{221}(34,62)(38,66)(42,66)
\put(34,62){\line(-1,-1){19}}
\bezier{221}(15,43)(11,39)(11,35)
\bezier{178}(42,66)(45.5,66)(47,62)
\bezier{185}(47,62)(48,59.5)(45,55)
\put(45,55){\line(-2,-3){21.33}}
\bezier{211}(23.5,23)(20,17)(15,18)
\put(11,35){\line(0,-1){12}}
\bezier{183}(11,25)(11,19)(15,18)
\end{picture}
}

\def\ovVI
{
\unitlength 4.2mm
\thinlines
\begin{picture}(16,20)(0,0)
\put(9,19){\oval(2,2)[lt]}
\put(8,19){\line(0,-1){6}}
\put(8,13){\line(0,1){0}}
\put(9,20){\line(1,0){6}}
\put(15,19){\oval(2,2)[rt]}
\put(9,13){\oval(2,2)[lb]}
\put(9,12){\line(1,0){6}}
\put(16,18){\line(0,1){0}}
\put(16,18){\line(0,1){1}}
\put(16,18){\line(0,-1){5}}
\put(15,13){\oval(2,2)[rb]}
\end{picture}

}
\def\kozulkanova{
{
\unitlength=1.6pt
\begin{picture}(60.00,34.00)(0.00,17.00)
\thicklines
\put(55.00,15.00){\makebox(0.00,0.00)[lb]{$h$}}
\put(7.00,11.00){\makebox(0.00,0.00)[rb]{$g$}}
\put(10.00,-5.00){\makebox(0.00,0.00)[t]{$l$}}
\put(34.00,34.00){\makebox(0.00,0.00)[lb]{$f$}}
\put(50.00,10.00){\makebox(0.00,0.00){$\bullet$}}
\put(10.00,10.00){\makebox(0.00,0.00){$\bullet$}}
\put(10.00,0.00){\makebox(0.00,0.00){$\bullet$}}
\put(30.00,30.00){\makebox(0.00,0.00){$\bullet$}}
\put(50.00,10.00){\line(1,-2){12.00}}
\put(50.00,10.00){\line(-1,-2){12.00}}
\put(30.00,30.00){\line(1,-1){20.00}}
\put(10.00,10.00){\line(1,-2){12.00}}
\put(10.00,10.00){\line(0,-1){10.00}}
\put(10.00,10.00){\line(-1,-2){12.00}}
\put(30.00,30.00){\line(-1,-1){20.00}}
\put(30.00,50.00){\line(0,-1){20.00}}
\end{picture}}
}

\def\kozulkanovaI{
{
\unitlength=1.6pt
\begin{picture}(60.00,30.00)(0.00,17.00)
\thinlines
\put(-.8,-13.5){\ovI}
\put(-5,-17){\ovVmensi}
\put(-62,-102){\ovVI}
\thicklines
\put(53.00,13.00){\makebox(0.00,0.00)[lb]{$h$}}
\put(6,11.50){\makebox(0.00,0.00)[rb]{$g$}}
\put(10.00,-5.00){\makebox(0.00,0.00)[t]{$l$}}
\put(34.00,34.00){\makebox(0.00,0.00)[lb]{$f$}}
\put(50.00,10.00){\makebox(0.00,0.00){$\bullet$}}
\put(10.00,10.00){\makebox(0.00,0.00){$\bullet$}}
\put(10.00,0.00){\makebox(0.00,0.00){$\bullet$}}
\put(30.00,30.00){\makebox(0.00,0.00){$\bullet$}}
\put(50.00,10.00){\line(1,-2){12.00}}
\put(50.00,10.00){\line(-1,-2){12.00}}
\put(30.00,30.00){\line(1,-1){20.00}}
\put(10.00,10.00){\line(1,-2){12.00}}
\put(10.00,10.00){\line(0,-1){10.00}}
\put(10.00,10.00){\line(-1,-2){12.00}}
\put(30.00,30.00){\line(-1,-1){20.00}}
\put(30.00,50.00){\line(0,-1){20.00}}
\end{picture}}
}

\def\kozulkanovaII{
{
\unitlength=1.6pt
\begin{picture}(60.00,30.00)(0.00,17.00)
\thinlines
\put(-.8,-13.5){\ovI}
\put(13,6){\ovII}
\put(-62,-102){\ovVI}
\thicklines
\put(54.00,15.00){\makebox(0.00,0.00)[lb]{$h$}}
\put(7.00,11.00){\makebox(0.00,0.00)[rb]{$g$}}
\put(10.00,-5.00){\makebox(0.00,0.00)[t]{$l$}}
\put(34.00,34.00){\makebox(0.00,0.00)[lb]{$f$}}
\put(50.00,10.00){\makebox(0.00,0.00){$\bullet$}}
\put(10.00,10.00){\makebox(0.00,0.00){$\bullet$}}
\put(10.00,0.00){\makebox(0.00,0.00){$\bullet$}}
\put(30.00,30.00){\makebox(0.00,0.00){$\bullet$}}
\put(50.00,10.00){\line(1,-2){12.00}}
\put(50.00,10.00){\line(-1,-2){12.00}}
\put(30.00,30.00){\line(1,-1){20.00}}
\put(10.00,10.00){\line(1,-2){12.00}}
\put(10.00,10.00){\line(0,-1){10.00}}
\put(10.00,10.00){\line(-1,-2){12.00}}
\put(30.00,30.00){\line(-1,-1){20.00}}
\put(30.00,50.00){\line(0,-1){20.00}}
\end{picture}}
}

\def\kozulkanovaIII{
{
\unitlength=1.6pt
\begin{picture}(60.00,30.00)(0.00,17.00)
\thinlines
\put(13,6){\ovII}
\put(-4,-15.5){\ovIIImensi}
\put(-62,-102){\ovVI}
\thicklines
\put(54.00,15.00){\makebox(0.00,0.00)[lb]{$h$}}
\put(6.00,11.00){\makebox(0.00,0.00)[rb]{$g$}}
\put(10.00,-5.00){\makebox(0.00,0.00)[t]{$l$}}
\put(34.00,35.00){\makebox(0.00,0.00)[lb]{$f$}}
\put(50.00,10.00){\makebox(0.00,0.00){$\bullet$}}
\put(10.00,10.00){\makebox(0.00,0.00){$\bullet$}}
\put(10.00,0.00){\makebox(0.00,0.00){$\bullet$}}
\put(30.00,30.00){\makebox(0.00,0.00){$\bullet$}}
\put(50.00,10.00){\line(1,-2){12.00}}
\put(50.00,10.00){\line(-1,-2){12.00}}
\put(30.00,30.00){\line(1,-1){20.00}}
\put(10.00,10.00){\line(1,-2){12.00}}
\put(10.00,10.00){\line(0,-1){10.00}}
\put(10.00,10.00){\line(-1,-2){12.00}}
\put(30.00,30.00){\line(-1,-1){20.00}}
\put(30.00,50.00){\line(0,-1){20.00}}
\end{picture}}
}
\def\kozulkanovaIV{
{
\unitlength=1.6pt
\begin{picture}(60.00,30.00)(0.00,17.00)
\thinlines
\put(39,6){\ovIV}
\put(-6,-18){\ovVtrochumensi}
\put(-62,-102){\ovVI}
\thicklines
\put(53.00,13.00){\makebox(0.00,0.00)[lb]{$h$}}
\put(6.00,11.50){\makebox(0.00,0.00)[rb]{$g$}}
\put(10.00,-5.00){\makebox(0.00,0.00)[t]{$l$}}
\put(34.00,35.00){\makebox(0.00,0.00)[lb]{$f$}}
\put(50.00,10.00){\makebox(0.00,0.00){$\bullet$}}
\put(10.00,10.00){\makebox(0.00,0.00){$\bullet$}}
\put(10.00,0.00){\makebox(0.00,0.00){$\bullet$}}
\put(30.00,30.00){\makebox(0.00,0.00){$\bullet$}}
\put(50.00,10.00){\line(1,-2){12.00}}
\put(50.00,10.00){\line(-1,-2){12.00}}
\put(30.00,30.00){\line(1,-1){20.00}}
\put(10.00,10.00){\line(1,-2){12.00}}
\put(10.00,10.00){\line(0,-1){10.00}}
\put(10.00,10.00){\line(-1,-2){12.00}}
\put(30.00,30.00){\line(-1,-1){20.00}}
\put(30.00,50.00){\line(0,-1){20.00}}
\end{picture}}
}
\def\kozulkanovaV{
{
\unitlength=1.6pt
\begin{picture}(60.00,30.00)(0.00,17.00)
\thinlines
\put(-6,-15.5){\ovIIImensi}
\put(39,6){\ovIV}
\put(-62,-102){\ovVI}
\thicklines
\put(54.50,11.00){\makebox(0.00,0.00)[lb]{$h$}}
\put(5.50,15.00){\makebox(0.00,0.00)[rb]{$g$}}
\put(10.00,-5.00){\makebox(0.00,0.00)[t]{$l$}}
\put(34.00,34.00){\makebox(0.00,0.00)[lb]{$f$}}
\put(50.00,10.00){\makebox(0.00,0.00){$\bullet$}}
\put(10.00,10.00){\makebox(0.00,0.00){$\bullet$}}
\put(10.00,0.00){\makebox(0.00,0.00){$\bullet$}}
\put(30.00,30.00){\makebox(0.00,0.00){$\bullet$}}
\put(50.00,10.00){\line(1,-2){12.00}}
\put(50.00,10.00){\line(-1,-2){12.00}}
\put(30.00,30.00){\line(1,-1){20.00}}
\put(10.00,10.00){\line(1,-2){12.00}}
\put(10.00,10.00){\line(0,-1){10.00}}
\put(10.00,10.00){\line(-1,-2){12.00}}
\put(30.00,30.00){\line(-1,-1){20.00}}
\put(30.00,50.00){\line(0,-1){20.00}}
\end{picture}}
}

\def\kozulkalabeled{
{
\unitlength=1.200000pt
\begin{picture}(60.00,40.00)(0.00,10.00)
\thicklines
\put(60.00,20.00){\makebox(0.00,0.00){$h$}}
\put(0.00,20.00){\makebox(0.00,0.00){$g$}}
\put(40.00,40.00){\makebox(0.00,0.00){$f$}}
\put(10,-20){\makebox(0.00,0.00)[t]{$l$}}
\put(29,22){\makebox(0.00,0.00)[r]{\scriptsize $1$}}
\put(40,28){\makebox(0.00,0.00)[r]{\scriptsize $2$}}
\put(50.00,10.00){\makebox(0.00,0.00){$\bullet$}}
\put(10.00,10.00){\makebox(0.00,0.00){$\bullet$}}
\put(30.00,30.00){\makebox(0.00,0.00){$\bullet$}}
\put(50.00,10.00){\line(1,-2){15}}
\put(50.00,10.00){\line(-1,-2){15}}
\put(30.00,30.00){\line(1,-1){20.00}}
\put(10.00,10.00){\line(1,-2){15.00}}
\put(10.00,10.00){\line(0,-1){25.00}}
\put(10.00,10.00){\line(-1,-2){15.00}}
\put(5.00,-3.00){\makebox(0.00,0.00)[l]{\scriptsize $1$}}
\put(12.00,-5.00){\makebox(0.00,0.00)[l]{\scriptsize $2$}}
\put(16.00,5.00){\makebox(0.00,0.00)[l]{\scriptsize $3$}}
\put(30.00,30.00){\line(-1,-1){20.00}}
\put(30.00,50.00){\line(0,-1){20.00}}
\put(20,-20){
\put(32,22){\makebox(0.00,0.00)[r]{\scriptsize $1$}}
\put(38.5,26.5){\makebox(0.00,0.00)[r]{\scriptsize $2$}}
}
\put(-5,-25){\makebox(0.00,0.00)[t]{$1$}}
\put(10,-15){\makebox(0.00,0.00){$\bullet$}}
\put(25,-25){\makebox(0.00,0.00)[t]{$2$}}
\put(35,-25){\makebox(0.00,0.00)[t]{$3$}}
\put(65,-25){\makebox(0.00,0.00)[t]{$4$}}
\end{picture}}
}

\def\kozulkalabeledtau{
{
\unitlength=1.200000pt
\begin{picture}(60.00,40.00)(-70.00,10.00)
\thicklines
\put(-60.00,20.00){\makebox(0.00,0.00){$h$}}
\put(0.00,20.00){\makebox(0.00,0.00){$g$}}
\put(-20.00,40.00){\makebox(0.00,0.00){$f\tau$}}
\put(-10,-20){\makebox(0.00,0.00)[t]{$l$}}
\put(-63,0){
\put(29,22){\makebox(0.00,0.00)[l]{\scriptsize $1$}}
\put(40,28){\makebox(0.00,0.00)[l]{\scriptsize $2$}}
}
\put(-50.00,10.00){\makebox(0.00,0.00){$\bullet$}}
\put(-10.00,10.00){\makebox(0.00,0.00){$\bullet$}}
\put(-30.00,30.00){\makebox(0.00,0.00){$\bullet$}}
\put(-50.00,10.00){\line(-1,-2){15}}
\put(-50.00,10.00){\line(1,-2){15}}
\put(-30.00,30.00){\line(-1,-1){20.00}}
\put(-10.00,10.00){\line(-1,-2){15.00}}
\put(-10.00,10.00){\line(0,-1){25.00}}
\put(-10.00,10.00){\line(1,-2){15.00}}
\put(-16,0){
\put(5.00,-3.00){\makebox(0.00,0.00)[r]{\scriptsize $1$}}
\put(12.00,-5.00){\makebox(0.00,0.00)[r]{\scriptsize $2$}}
\put(16.00,5.00){\makebox(0.00,0.00)[r]{\scriptsize $3$}}
}
\put(-30.00,30.00){\line(1,-1){20.00}}
\put(-30.00,50.00){\line(0,-1){20.00}}
\put(-84,-20){
\put(32,22){\makebox(0.00,0.00)[l]{\scriptsize $1$}}
\put(38.5,26.5){\makebox(0.00,0.00)[l]{\scriptsize $2$}}
}
\put(5,-25){\makebox(0.00,0.00)[t]{$2$}}
\put(-10,-15){\makebox(0.00,0.00){$\bullet$}}
\put(-25,-25){\makebox(0.00,0.00)[t]{$1$}}
\put(-35,-25){\makebox(0.00,0.00)[t]{$4$}}
\put(-65,-25){\makebox(0.00,0.00)[t]{$3$}}
\end{picture}}
}

\def\kozulkalabeledtaubis{
{
\unitlength=1.200000pt
\begin{picture}(60.00,40.00)(-70.00,10.00)
\thicklines
\put(-60.00,20.00){\makebox(0.00,0.00){$h$}}
\put(0.00,20.00){\makebox(0.00,0.00){$g$}}
\put(-20.00,40.00){\makebox(0.00,0.00){$f$}}
\put(-10,-20){\makebox(0.00,0.00)[t]{$l$}}
\put(-63,0){
\put(29,22){\makebox(0.00,0.00)[l]{\scriptsize $2$}}
\put(40,28){\makebox(0.00,0.00)[l]{\scriptsize $1$}}
}
\put(-50.00,10.00){\makebox(0.00,0.00){$\bullet$}}
\put(-10.00,10.00){\makebox(0.00,0.00){$\bullet$}}
\put(-30.00,30.00){\makebox(0.00,0.00){$\bullet$}}
\put(-50.00,10.00){\line(-1,-2){15}}
\put(-50.00,10.00){\line(1,-2){15}}
\put(-30.00,30.00){\line(-1,-1){20.00}}
\put(-10.00,10.00){\line(-1,-2){15.00}}
\put(-10.00,10.00){\line(0,-1){25.00}}
\put(-10.00,10.00){\line(1,-2){15.00}}
\put(-16,0){
\put(5.00,-3.00){\makebox(0.00,0.00)[r]{\scriptsize $1$}}
\put(12.00,-5.00){\makebox(0.00,0.00)[r]{\scriptsize $2$}}
\put(16.00,5.00){\makebox(0.00,0.00)[r]{\scriptsize $3$}}
}
\put(-30.00,30.00){\line(1,-1){20.00}}
\put(-30.00,50.00){\line(0,-1){20.00}}
\put(-84,-20){
\put(32,22){\makebox(0.00,0.00)[l]{\scriptsize $1$}}
\put(38.5,26.5){\makebox(0.00,0.00)[l]{\scriptsize $2$}}
}
\put(5,-25){\makebox(0.00,0.00)[t]{$2$}}
\put(-10,-15){\makebox(0.00,0.00){$\bullet$}}
\put(-25,-25){\makebox(0.00,0.00)[t]{$1$}}
\put(-35,-25){\makebox(0.00,0.00)[t]{$4$}}
\put(-65,-25){\makebox(0.00,0.00)[t]{$3$}}
\end{picture}}
}

\begin{document}
\bibliographystyle{plain}
\baselineskip 1.3em

\title{Operads and PROPs}
\author[M. Markl]{Martin MARKL}
\thanks{The author was supported by the grant GA \v CR 201/05/2117 and by
   the Academy of Sciences of the Czech Republic,
   Institutional Research Plan No.~AV0Z10190503.}
\date{January  25, 2006}
\email{markl@math.cas.cz}

\address{Mathematical Institute of the Academy, {\v Z}itn{\'a} 25,
         115 67 Prague 1, The Czech Republic}

\begin{abstract}
We review definitions and basic properties of operads, \PROP{s} and
algebras over these structures.
\end{abstract}

\maketitle

\begin{center}
{\em Dedicated to the memory of Jakub Jan Ryba (1765--1815)\/}
\end{center}

\hyphenation{au-toch-tho-nous}


Operads involve an abstraction of the family $\{{\it
Map}(X^n,X)\}_{n\geq 0}$ of composable functions of several variables
together with an action of permutations of variables. As such, they
were originally studied as a tool in homotopy theory, specifically for
iterated loop spaces and homotopy invariant structures, but the theory
of operads has recently received new inspiration from homological
algebra, category theory, algebraic geometry and mathematical
physics. The name {\it operad} and the formal definition appear first
in the early 1970's in J.P.~May's book~\cite{may:1972}, but a year or
more earlier, M.~Boardman and R.~Vogt~\cite{boardman-vogt:BAMS68}
described the same concept under the name categories of {\em operators
in standard form\/}, inspired by {\small PROP}s and {\small PACT}s of
Adams and Mac~Lane~\cite{maclane:RiceUniv.Studies63}.  As pointed out
in~\cite{leinster:TAC04}, also Lambek's definition of
multicategory~\cite{lambek:69}
(late 1960s) was almost equivalent to what is called today a
colored or many-sorted operad. Another important precursor was the
associahedron $K$ that appeared in J.D.~Stasheff's 1963
paper~\cite{stasheff:TAMS63} on homotopy associativity of H-spaces.
We do not, however, aspire to write an account on the \underline{histor}y of
operads and their applications here -- we refer to the introduction
of~\cite{markl-shnider-stasheff:book}, to~~\cite{may:202over},
\cite{vogt:hist}, or to the report~\cite{stasheff:report} instead.

Operads are important not in and of themselves but, like {\small
PROP}s, through their representations, more commonly called {\em
algebras\/} over operads or {\em operad algebras\/}. If an operad
is thought of as a kind of algebraic theory, then an algebra over
an operad is a model of that theory.  Algebras over operads
involve most of `classical' algebras (associative, Lie,
commutative associative, Poisson,~\&c.), loop spaces, moduli
spaces of algebraic curves, vertex operator algebras, \&c. {\em
Colored\/} or {\em many-sorted\/} operads then describe diagrams
of homomorphisms of these objects, homotopies between
homomorphisms, modules,~\&c.

\PROP{s} generalize operads in the sense that they admit operations
with several inputs and {\em several\/} outputs. Therefore various
bialgebras (associative, Lie, infinitesimal) are \PROP{ic}
algebras. \PROP{s} were also used to encode `profiles' of structures
in formal differential geometry~\cite{merkulov:PROP,merkulov:defquant}.

By the \underline{renaissance} of operads we mean the first half of the
nineties of the last century when several papers which stimulated the
rebirth of interest in operads
appeared~\cite{getzler-jones:preprint,ginzburg-kapranov:DMJ94,%
       hinich-schechtman:AdvinSovMath1993,Huang:CMP94,hl,%
       kimura-stasheff-voronov:1993,markl:ws93}.
Let us mention the most important new ideas that emerged during this period.

First of all, operads were recognized as the underlying combinatorial
structure of the moduli space of stable algebraic curves in complex
geometry, and of compactifications of configuration spaces of points
of affine spaces in real geometry.  In mathematical physics, several
very important concepts such as vertex operator algebras or various
string theories were interpreted as algebras over operads.
On the algebraic side, the notion of {\em Koszulness\/}
of operads was introduced and studied, and the relation between resolutions
of operads and deformations of their algebras was recognized.
See~\cite{loday:bourbaki} for an autochthonous account of the renaissance.
Other papers that later became influential
then followed in a rapid
succession~\cite{getzler:genus0,getzler-kapranov:CompM98,getzler-kapranov:CPLNGT95,%
       markl:JPAA96,markl:zebrulka}.

Let us list some most important outcomes of the renaissance of
operads.  The choice of the material for this incomplete catalog has
been of course influenced by the author's personal expertise and
inclination towards algebra, geometry and topics that are commonly
called mathematical physics. We will therefore not be able to pay as
much attention to other aspects of operads, such as topology, category
theory and homotopy theory, as they deserve.

{\it Complex geometry.} Applications involve moduli spaces of stable
complex algebraic curves of genus zero~\cite{ginzburg-kapranov:DMJ94},
enumerative geometry, Frobenius
manifolds, quantum cohomology and cohomological
field theory~\cite{kontsevich-manin:CMP94,manin:frob}.
The moduli space of genus zero curves
exhibits an additional symmetry that leads to a generalization called
{\em cyclic\/} operads~\cite{getzler-kapranov:CPLNGT95}. {\em Modular\/}
operads~\cite{getzler-kapranov:CompM98} then describe
the combinatorial structure of the space of curves of arbitrary genus.

{\it Real geometry.} Compactifications of configuration spaces of
points in real smooth manifolds are operads in the category of smooth
manifolds with corners or modules over these
operads~\cite{markl:cf}. This fact is crucial for the theory of
configuration spaces with summable labels~\cite{salvatore}. The cacti
operad~\cite{voronov:01} lies behind the Chas-Sullivan product on the free
loop space of a smooth manifold~\cite{chas-sullivan}, see
also~\cite{cohen-voronov}.  Tamarkin's proof of the formality of
Hochschild cochains of the algebra of functions on smooth
manifolds~\cite{tamarkin:deligne} explained in~\cite{hinich:FM03} uses
obstruction theory for operad algebras and the affirmative answer to the
Deligne conjecture~\cite{deligne:letter,kontsevich-soibelman}.

{\it Mathematical physics.} The formality mentioned in the previous
item implies the existence of the deformation quantization of Poisson
manifolds~\cite{kontsevich:LMPh03}. We must not forget to mention the
operadic interpretation of vertex operator algebras~\cite{huang:book},
string theory~\cite{kimura-stasheff-voronov:1993} and Connes-Kreimer's
approach to renormalization~\cite{connes-kreimer:I}. Operads and
multicategories are important also for Beilinson-Drinfeld's theory of
chiral algebras~\cite{beilinson-drinfeld:chiral_algebras}.

{\it Algebra.}  Operadic
cohomology~\cite{balavoine:jpaa98,fox-markl:ContM97,getzler-jones:preprint,%
                 ginzburg-kapranov:DMJ94,markl-shnider-stasheff:book}
provides a uniform treatment of all
`classical' cohomology theories, such as the Hochschild cohomology
of associative algebras, Harrison cohomology of associative
commutative algebras, Chevalley-Eilenberg cohomology of Lie
algebras,~\&c.  Minimal models for operads~\cite{markl:zebrulka} offer
a conceptual understanding of strong homotopy algebras, their
homomorphisms and homotopy invariance~\cite{markl:ha}. Operads serve
as a natural language for various types of
`multialgebras'~\cite{loday:digebres,LFCG}.  Relation between
Koszulness of operads and properties of posets was studied
in~\cite{fresse:02}. Also the concept of the operadic distributive law
turned out to be useful~\cite{fox-markl:ContM97,markl:dl}.

{\it Model structures.}  It turned
out~\cite{berger-moerdijk:coloured,%
getzler-jones:preprint,hinich:CA97,schwanzl-vogt:Archiv91} 
that algebras over a
reasonable (possibly colored) operad form a model category that generalizes the
classical model structures of the categories of dg commutative
associative algebras and dg Lie
algebras~\cite{quillen:Ann.ofMath.69,sullivan:Publ.IHES77}.
Operads, in a reasonable monoidal model category, themselves form
a model category~\cite{berger-moerdijk:02,getzler-jones:preprint} 
such that algebras
over cofibrant operads are homotopy invariant, see also~\cite{spitzweck}.
Minimal operads mentioned in the previous 
item are particular cases of cofibrant
dg-operads and the classical
$W$-construction~\cite{boardman-vogt:BAMS68} is a functorial
cofibrant replacement in the category of topological
operads~\cite{vogt:TopAppl03}. The above model structures are
important for various constructions in the homology theory of
(free or based) loop spaces~\cite{cohen-voronov,hu} and
formulations of `higher' Deligne conjecture~\cite{hu-kriz}.

{\it Topology.} Operads as gadgets organizing homotopy
coherent structures are important in the brave new algebra approach to
topological Hochschild cohomology and algebraic
\hbox{$K$-theory}, see~\cite{elmendorf-kriz-mandell-may:stable,%
elmendorf-kriz-mandell-may:book,may:1998,vogt:brave},
or~\cite{elmendorf:04} for a historical background.
A description of a localized category of integral and $p$-adic
homotopy types by $E_\infty$-operads was given
in~\cite{mandell:TOP01,mandell:integral}.
An operadic approach to partial algebras and their completions was
applied in~\cite{kriz-may} to mixed Tate motives over the
rationals. See also an overview~\cite{may:202}.

{\it Category theory.}  
Operads and multicategories were used as a
language in which to propose a definition of weak
$\omega$-category~\cite{batanin98,batanin:eckmann,batanin,leinster:book}.
Operads themselves can be viewed as special kinds of algebraic theory
(as can multicategories, if one allows many-sorted theories),
see~\cite{mauri}.  There are also `categorical' generalizations of
operads, e.g.~the globular operads of~\cite{batanin:globular} and
$T$-categories of~\cite{burroni}. An interesting presentation of \PROP-like
structures in enriched monoidal categories can be found
in~\cite{may:caterads}. 

{\it Graph Complexes.}  Each cyclic operad $\calP$ determines a graph
complex~\cite{getzler-kapranov:CompM98,markl:ws98}.  As observed
earlier by M.~Kontsevich~\cite{kontsevich:93}, these graph complexes
are, for some specific choices of $\calP$, closely related to some
very interesting objects such as moduli spaces of Riemann surfaces,
automorphisms of free groups or primitives in the homology of
certain infinite-dimensional Lie
algebras, see also~\cite[II.5.5]{markl-shnider-stasheff:book}.  In the
same vein, complexes of directed graphs are related to {\small
PROP}s~\cite{mv,vallette:thesis,vallette:extract,vallette:CMR04} and
directed graphs with back-in-time
edges are tied to wheeled PROPs introduced
in~\cite{merkulov:defquant}.

{\it Deformation theory and homotopy invariant structures in algebra.}
A concept of homotopy invariant structures in algebra parallel to the
classical one in topology~\cite{boardman-vogt:BAMS68,boardman-vogt:73}
was developed in~\cite{markl:ha}. It was explained
in~\cite{kontsevich-soibelman,markl:JPAA96,markl:ba} how
cofibrant resolutions of operads or \PROP{s} determine a cohomology
theory governing deformations of related
algebras. In~\cite{markl:ib}, deformations were identified with
solutions of the Maurer-Cartan equation of a certain strongly homotopy Lie
algebra constructed in a very explicit way from a cofibrant resolution
of the underlying operad or \PROP.

\begin{center}
-- -- -- -- --
\end{center}

\noindent
{\bf Terminology.}
As we already observed, operads are abstractions of families of
composable functions. Given functions $f : X^{\times n} \to X$ and
$g_i : X^{\times k_i} \to X$, $1 \leq i \leq n$, one may
consider the simultaneous composition
\hfill\break
(I) \rule{0em}{1.5em} \hfill
$f(\Rada g1n) : X^{\times {(k_1 + \cdots + k_n)}} \to X.$ \hfill \break
\rule{0em}{1.5em}%
One may also consider, for $f : X^{\times n} \to X$,
$g : X^{\times m} \to X$ and $1 \leq i \leq n$, the individual compositions
\hfill\break
(II)\rule{0em}{1.5em} \hfill
$f(\id,\cdots,\id,g,\id,\cdots,\id) : X^{\times (m+m-1)} \to X$,
\hfill\break
\rule{0em}{1.5em}%
with  $g$ at the $i$th place. While May's original definition of
an operad~\cite{may:1972} was an abstraction of type~(I)
compositions, there exist an alternative approach based on
type~(II) compositions. This second point of view was formalized
in the 1963~papers by Gerstenhaber~\cite{gerstenhaber:AM63} and
Stasheff~\cite{stasheff:TAMS63}. A definition that included the
symmetric group action was formulated much later in the author's
paper~\cite{markl:zebrulka} in which the two approaches were also
compared.

In the presence of operadic units, these approaches are equivalent.
There are, however, situations where one needs also
non-unital versions, and then the two approaches lead to
\underline{different} structures -- a non-unital
structure of the second type always determines a non-unital
structure of the first type, but {\em not vice versa\/}! It turns out that
more common are non-unital structures of the second type; they describe,
for example, the underlying combinatorial structure of the moduli
space of stable complex curves.

We will therefore call the non-unital
versions of the first type of operads {\em non-unital May's
operads\/}, while the second version simply {\em non-unital
operads\/}. We opted for this terminology, which was used already in
the first version of~\cite{markl:zebrulka}, after a long hesitation,
being aware that it might not be universally welcome. Note
that non-unital operads are sometimes called {\em (Markl's)
pseudo-operads\/}~\cite{markl:zebrulka,markl-shnider-stasheff:book}.

\vskip 1cm
\noindent
{\sc Outline of the paper:}
     \ref{sec1}.~Operads -- page~\pageref{sec1} \hfill\break
\hphantom{{\sc Outline of the paper:}}
     \ref{sec1bis}.~Non-unital operads -- page~\pageref{sec1bis}
\hfill\break
\hphantom{{\sc Outline of the paper:}}
     \ref{sec2}.~Operad algebras -- page~\pageref{sec2}
\hfill\break
\hphantom{{\sc Outline of the paper:}}
     \ref{sec3}.~Free operads and trees -- page~\pageref{sec3}
\hfill\break
\hphantom{{\sc Outline of the paper:}}
     \ref{sec4}.~Unbiased definitions -- page~\pageref{sec4}
\hfill\break
\hphantom{{\sc Outline of the paper:}}
     \ref{sec5}.~Cyclic operads -- page~\pageref{sec5}
\hfill\break
\hphantom{{\sc Outline of the paper:}}
     \ref{sec6}.~Modular operads -- page~\pageref{sec6}
\hfill\break
\hphantom{{\sc Outline of the paper:}}
     \ref{sec7}.~\PROP{s} -- page~\pageref{sec7}
\hfill\break
\hphantom{{\sc Outline of the paper:}}
     \ref{sec8}.~Properads, dioperads and $\frac 12$PROPs -- page~\pageref{sec8}
\hfill\break
\hphantom{{\sc Outline of the paper:}}
     \hphantom{9.}~References -- page~\pageref{refs}

\vskip .5cm

In the first three sections we review basic definitions of (unital
and non-unital) operads and operad algebras, and give examples
that illustrate these notions. The fourth section describes free
operads and their relation to rooted trees. In the fifth section we
explain that operads can be defined as algebras over the monad of rooted
trees. In the following two sections we show that, replacing rooted
trees by other types of trees, one obtains two important
generalizations  -- cyclic and modular
operads. In the last two sections, \PROP{s} and their versions are
recalled; this article is the first expository text where these structures are
systematically treated.

Sections~\ref{sec1}--\ref{sec2} are based on the classical
book~\cite{may:1972} by J.P.~May and 
the author's article~\cite{markl:zebrulka}.
Sections~\ref{sec3}--\ref{sec6} follow the seminal
paper~\cite{ginzburg-kapranov:DMJ94} by V.~Ginzburg and M.M.~Kapranov,
and papers~\cite{getzler-kapranov:CPLNGT95,getzler-kapranov:CompM98}
by E.~Getzler and M.M.~Kapranov. The last two sections are based on
the preprint~\cite{mv} of  A.A.~Voronov and the author, and on an e-mail
message~\cite{kontsevich:message} from M.~Kontsevich.
We were also influenced by T.~Leinster's
concept of biased versus un-biased
definitions~\cite{leinster:book}. At some places, our exposition
follows the monograph~\cite{markl-shnider-stasheff:book} by
S.~Shnider, J.D.~Stasheff and the author.

\vskip .5em \noindent {\bf Acknowledgment.} I would like to
express my gratitude to Ezra Getzler, Tom Leinster,  Peter May, 
Jim Stasheff, Bruno Vallette and Rainer
Vogt for careful reading the manuscript and many useful
suggestions. I would also like to thank the Institut des Hautes
\'Etudes Scientifiques, Bures-sur-Yvette, France, for the
hospitality during the period when this article was completed.

\section{Operads}
\label{sec1}

Although operads, operad algebras and most of related structures can be
defined in an arbitrary symmetric monoidal category with countable
coproducts, we decided to follow the choice
of~\cite{kriz-may} and formulate precise definitions only for the category
${\tt Mod}_\bfk = ({\tt Mod}_\bfk,\otimes)$ of modules over 
a~commutative unital ring ${\bfk}$, with the monoidal structure given by
the tensor product $\otimes = \ot_\bfk$ over $\bfk$.  The reason for
such a decision was to give, in Section~\ref{sec3}, a clean
construction of free operads. 
In a general monoidal category, this construction involves the unordered
$\odot$-product~\cite[Definition~II.1.38]{markl-shnider-stasheff:book}
so the free operad is then a double colimit,
see~\cite[Section~II.1.9]{markl-shnider-stasheff:book}.
Our choice also allows us to
write formulas involving maps in terms of
elements, which is sometimes a~welcome simplification.
We believe that the reader can easily reformulate our definitions
into other monoidal categories or
consult~\cite{markl-shnider-stasheff:book,may:202defs}. 

Let $\bfk[\Sigma_n]$
denote the $\bfk$-group ring of the symmetric group $\Sigma_n$.

\begin{definition}[May's operad]
\label{a}
An {\em operad\/} in the category of $\bfk$-modules is a collection
$\calP =\{\calP(n)\}_{n\geq 0}$ of right $\bfk[\Sigma_n]$-modules, together
with $\bfk$-linear maps (operadic compositions)
\begin{equation}
\label{pocasi_je_na_nic}
\gamma : \calP(n) \ot \calP(k_1) \ot \cdots \ot \calP(k_n)
\to
\calP(k_1 + \cdots + k_n),
\end{equation}
for $n \geq 1$ and $k_1,\ldots,k_n \geq 0$, and a unit map $\eta :
\bfk \to \calP(1)$. These data fulfill the following axioms.

\noindent
{\em Associativity.\/} Let $n \geq 1$ and let $m_1,\ldots,m_n$ and
$k_1,\ldots,k_m$, where $m := m_1 + \cdots + m_n$, be non-negative
integers. Then the following diagram, in which $g_s := m_1 + \cdots
+ m_{s-1}$ and $h_s = k_{g_s + 1} \cdots + k_{g_{s+1}}$, for $1 \leq s
\leq n$,  commutes.
\begin{center}
\setlength{\unitlength}{1em}
\begin{picture}(20,15)(-2,-1)
\put(0,12){\makebox(0,0){$
    \displaystyle \left(
    \calP(n) \ot \bigotimes_{s=1}^n
    \calP(m_s)\right) \ot \bigotimes_{r=1}^m
    \calP(k_r)$}}
\put(8,12){\vector(1,0){7.5}}
\put(11,12.4){\makebox(0,0)[b]{\scriptsize $\gamma \ot \id$}}
\put(20,12){\makebox(0,0){$ \displaystyle\calP(m) \ot  \bigotimes_{r=1}^m
    \calP(k_r)$}}
\put(0,0){\makebox(0,0){$
    \displaystyle
    \calP(n) \ot \bigotimes_{s=1}^n  \left(
    \calP(m_s) \ot\bigotimes_{q=1}^{m_s}
    \calP(k_{g_s + q})\right)
     $}}
\put(20,0){\makebox(0,0){$ \displaystyle\calP(n) \ot  \bigotimes_{s=1}^n
    \calP(h_s)$}}
\put(8.5,0){\vector(1,0){7}}
\put(12,.4){\makebox(0,0)[b]{\scriptsize $\id \ot  (\bigotimes_{s= \rm
      1}^n \gamma)$}}
\put(0,10){\vector(0,-1){8}}
\put(-.4,6){\makebox(0,0)[r]{\scriptsize  shuffle}}
\put(20,6){\makebox(0,0){$\calP(k_1 + \cdots + k_m)$}}
\put(20,10){\vector(0,-1){3}}
\put(20.4,8.5){\makebox(0,0)[l]{\scriptsize $ \gamma$}}
\put(20,2){\vector(0,1){3}}
\put(20.4,3.5){\makebox(0,0)[l]{\scriptsize $ \gamma$}}
\end{picture}
\end{center}

\noindent
{\em Equivariance.\/} Let $n \geq 1$, let $k_1,\ldots,k_n$ be
non-negative integers and $\sigma \in \Sigma_n$, $\tau_1 \in
\Sigma_{k_1},\ldots, \tau_n \in \Sigma_{k_n}$ permutations.  Let
$\sigma(k_1,\ldots,k_n) \in \Sigma_{k_1 + \cdots + k_n}$ denote the
permutation that permutes $n$ blocks
$(1,\ldots,k_1),\ldots,(k_{n-1}+1,\ldots,k_n)$ as $\sigma$ permutes
$(1,\ldots,n)$ and let $\tau_1 \oplus \cdots \oplus \tau_n \in
\Sigma_{k_1 + \cdots + k_n}$ be the block sum of permutations. Then
the following diagrams commute.
\begin{center}
\setlength{\unitlength}{1em}
\begin{picture}(20,14.5)(0,-.5)
\put(0,8){
\put(0,5){\makebox(0,0){$\calP(n) \ot \calP(k_1) \ot \cdots \ot
    \calP(k_n)$}}
\put(20,5){\makebox(0,0){$\calP(n) \ot \calP(k_{\sigma(1)}) \ot \cdots
    \ot \calP(k_{\sigma(n)})$}}
\put(0,0){\makebox(0,0){$\calP(k_1 + \cdots + k_n)$}}
\put(20,0){\makebox(0,0){$\calP(k_{\sigma(1)} + \cdots + k_{\sigma(n)})$}}
\put(10,5.5){\makebox(0,0)[b]{\scriptsize $\sigma \ot \sigma^{-1}$}}
\put(6.5,5){\vector(1,0){6}}
\put(10,.5){\makebox(0,0)[b]{\scriptsize $\sigma(k_{\sigma(1)},\ldots,k_{\sigma(n)})$}}
\put(4,0){\vector(1,0){11}}
\put(-.5,2.5){\makebox(0,0)[r]{\scriptsize $\gamma$}}
\put(0,4){\vector(0,-1){3}}
\put(20.5,2.5){\makebox(0,0)[l]{\scriptsize $\gamma$}}
\put(20,4){\vector(0,-1){3}}
}
\put(0,5){\makebox(0,0){$\calP(n)\ot \calP(k_1) \ot \cdots \ot \calP(k_n)$}}
\put(20,5){\makebox(0,0){$\calP(n) \ot \calP(k_1) \ot \cdots
    \ot \calP(k_n)$}}
\put(0,0){\makebox(0,0){$\calP(k_1 + \cdots + k_n)$}}
\put(20,0){\makebox(0,0){$\calP(k_1 + \cdots + k_n)$}}
\put(10,5.5){\makebox(0,0)[b]{\scriptsize $\id \ot \tau_1 \ot \cdots \ot  \tau_n$}}
\put(6.5,5){\vector(1,0){7}}
\put(10,.5){\makebox(0,0)[b]{\scriptsize $\tau_1 \oplus \cdots \oplus \tau_n$}}
\put(4,0){\vector(1,0){12}}
\put(-.5,2.5){\makebox(0,0)[r]{\scriptsize $\gamma$}}
\put(0,4){\vector(0,-1){3}}
\put(20.5,2.5){\makebox(0,0)[l]{\scriptsize $\gamma$}}
\put(20,4){\vector(0,-1){3}}
\end{picture}
\end{center}

\noindent
{\em Unitality.\/} For each $n \geq 1$, the following diagrams commute.
\begin{center}
\setlength{\unitlength}{1em}
\begin{picture}(25,6)
\put(-.25,5){\makebox(0,0){$\calP(n) \ot \bfk^{\ot n}$}}
\put(-.4,2.5){\makebox(0,0)[r]{\scriptsize $\id \ot \eta^{\ot n}$}}
\put(0,4){\vector(0,-1){3}}
\put(4.5,5.5){\makebox(0,0)[b]{\scriptsize $\cong$}}
\put(3,5){\vector(1,0){3.5}}
\put(4,2.9){\makebox(0,0)[b]{\scriptsize $\gamma$}}
\put(2,1){\vector(3,2){5}}
\put(8,5){\makebox(0,0){$\calP(n)$}}
\put(.5,0){\makebox(0,0){$\calP(n) \ot {\calP(1)}^{\ot n}$}}
%
\put(20,0){
\put(.7,5){\makebox(0,0){$\bfk \ot \calP(n)$}}
\put(-.4,2.5){\makebox(0,0)[r]{\scriptsize $\eta \ot \id$}}
\put(0,4){\vector(0,-1){3}}
\put(4.5,5.5){\makebox(0,0)[b]{\scriptsize $\cong$}}
\put(3,5){\vector(1,0){3.5}}
\put(4,2.9){\makebox(0,0)[b]{\scriptsize $\gamma$}}
\put(2,1){\vector(3,2){5}}
\put(8,5){\makebox(0,0){$\calP(n)$}}
\put(0.1,0){\makebox(0,0){$\calP(1) \ot \calP(n)$}}
}
\end{picture}
\end{center}
\end{definition}

A straightforward modification of the above definition makes sense in
any symmetric monoidal category $(\calM,\odot,{\bf 1})$ such as the
category of differential graded modules, simplicial sets, topological
spaces, \&c, 
see~\cite[Definition~II.1.4]{markl-shnider-stasheff:book}
or ~\cite[Definition~1]{may:202defs}. We then
speak about {\em differential graded\/} operads, {\em simplicial\/}
operads, {\em topological\/} operads,~\&c.

\begin{example}
\label{end}
\rm
All properties axiomatized by Definition~\ref{a} can be read
from the {\em endomorphism operad\/}
$\End_V = \{\End_V(n)\}_{n \geq 0}$ of a $\bfk$-module~$V$. It
is defined by setting
$\End_V(n)$ to be the space of $\bfk$-linear maps
$V^{\ot n} \to V$.
The operadic composition of
$f \in \End_V(n)$ with $g_1 \in \End_V(k_1),\ldots,g_n \in \End_V(k_n)$
is given by the usual composition of multilinear maps as
\[
\gamma(f,g_1,\ldots,g_n) := f(g_1 \ot \cdots \ot g_n),
\]
the symmetric group acts by
\[
\gamma \sigma(f,g_1,\ldots,g_n) := f(g_{\sigma^{-1}(1)} \ot \cdots \ot
g_{\sigma^{-1}(n)}),\
\sigma \in \Sigma_n,
\]
and  the unit map
$\eta : \bfk \to \End_V(1)$ is given by $\eta(1) := \id_V : V \to V$.
The endomorphism operad can be constructed over an object of
an arbitrary symmetric monoidal category with an internal hom-functor,
as it was done in~\cite[Definition.~II.1.7]{markl-shnider-stasheff:book}.
\end{example}

One often considers operads $\calA$ such that $\calA(0) = 0$ (the
trivial $\bfk$-module). We will indicate that $\calA$ is of this type
by writing $\calA = \coll\calA1$.

\begin{example}
\label{9}
\rm Let us denote by $\Associative = \coll\Associative1$ the operad
with $\Associative(n) := \bfk[\Sigma_n]$, $n \geq 1$, and the operadic
composition defined as follows. Let $\id_n \in \Sigma_n$, $\id_{k_1}
\in \Sigma_{k_1},\ldots,\id_{k_n} \in \Sigma_{k_n}$ be the identity
permutations. Then
\[
\gamma(\id_n,\id_{k_1}, \ldots, \id_{k_n}) := \id_{k_1+ \cdots + k_n}
\in \Sigma_{k_1 + \cdots + k_n}.
\]
The above formula determines $\gamma(\sigma,\tau_1,\ldots\tau_n)$ for
general $\sigma \in \Sigma_n$, $\tau_1 \in \Sigma_{k_1},\ldots,\tau_n
\in \Sigma_{k_n}$ by the equivariance axiom. The unit map $\eta : \bfk
\to \Associative(1)$ is given by $\eta(1):= \id_1$.
\end{example}

\begin{example}
\label{tuto_sobotu_zacinaji_zavody}
{\rm
Let us give an example of a topological operad. For $k \geq 1$, the {\em little
$k$-discs operad\/} $\calD_k = \{\calD_k(n)\}_{n \geq 0}$ is defined as
follows~\cite[Section~II.4.1]{markl-shnider-stasheff:book}.
Let
\[
\DD^k := \{(x_1,\ldots,x_k) \in {\mathbb R}^k;\ x_1^2 + \cdots
+ x_k^2 \leq 1\}
\]
be the standard closed disc in ${\mathbb R}^k$. A
little $k$-disc is then a linear embedding $d : \DD^k \hookrightarrow \DD^k$
which is the restriction of a linear map ${\mathbb R}^k \to {\mathbb R}^k$
with parallel axes. The $n$-th space $\calD_k(n)$ of the little
$k$-disc operad is the space of all $n$-tuples $(d_1,\ldots,d_n)$ of
little $k$-discs such that the images of $d_1$,\ldots,$d_n$ have
mutually disjoint interiors. The operad structure is obvious -- the
symmetric group $\Sigma_n$ acts on $\calD_k(n)$ by permuting the
labels of the little discs and the structure map $\gamma$ is given by
composition of embeddings. The unit is the identity embedding
$\id : \DD^k \hookrightarrow \DD^k$.
}\end{example}

\begin{example}
\label{projdu_prohlidkou?}
{\rm
The collection of normalized singular chains
$C_*(\calT) = \{C_*(\calT(n))\}_{n \geq 0}$ of a topological operad
$\calT = \coll \calT0$  is an operad in
the category of differential graded ${\mathbb Z}$-modules. For a
ring $R$, the singular homology $H_*(\calT(n);R) =
H_*(C_*(\calT(n))\ot_{\mathbb Z} R)$ forms an operad $H_*(\calT;R)$ in
the category of graded $R$-modules,
see~\cite[Section~I.5]{kriz-may} for details.
}
\end{example}

\begin{definition}
\label{pristi_tyden_prohlidka}
Let $\calP = \coll \calP0$ and $\calQ = \coll\calQ0$ be two operads. A
{\em homomorphism\/} $f : \calP \to \calQ$ is a sequence $f
= \{f(n) : \calP(n) \to \calQ(n)\}_{n\geq 0}$ of equivariant maps
which commute with the operadic compositions and preserve the units.

An operad $\calR = \coll\calR0$ is a {\em suboperad\/} of\/ $\calP$
if $\calR(n)$ is, for each $n \geq 0$, a $\Sigma_n$-submodule of
$\calP(n)$ and if all structure operations of\/ $\calR$ are the
restrictions of those of\/ $\calP$. Finally, an {\em ideal\/} in the
operad $\calP$ is the collection $\calI = \coll\calI0$ of
$\Sigma_n$-invariant subspaces $\calI(n) \subset \calP(n)$ such that
\[
\gamma_\calP(f,\Rada g1n) \in \calI(k_1 + \cdots + k_n)
\]
if either $f \in \calI(n)$ or $g_i \in \calI(k_i)$ for some $1 \leq i
\leq n$.
\end{definition}

\begin{example}
\label{zitra_strasne_vedro}
{\rm
Given an operad $\calP = \coll \calP0$, let $\ccalP =
\coll\ccalP0$ be the collection defined by $\ccalP(n) := \calP(n)$ for
$n \geq 1$ and $\ccalP(0) := 0$. Then $\ccalP$ is a suboperad of
$\calP$.  The correspondence $\calP \mapsto \ccalP$ is a full
embedding of the category of operads $\calP$ with $\calP(0) \cong
\bfk$ into the category of operads~$\calA$ with $\calA(0) =
0$. Operads satisfying $\calP(0) \cong \bfk$ have been called {\em
unital\/} while operads with $\calA(0) = 0$ {\em
non-unital\/} operads. We will not use this terminology because
non-unital operads will mean something different in this article, see
Section~\ref{sec1bis}.

An example of an operad $\calA$ which is not of the
form $\ccalP$ for some operad $\calP$ with $\calP(0) \cong \bfk$ can
be constructed as follows.
Observe first that operads $\calP$ with the property that
\[
\mbox {$\calP(0) \cong \bfk$ and $\calP(n) =
0$ for $n \geq 2$}
\]
are the same as augmented associative algebras.
Indeed, the space $\calP(1)$ with the operation
$\circ_1 : \calP(1) \ot \calP(1) \to \calP(1)$
is clearly a unital associative algebra,  augmented by
the composition
\[
\calP(1) \stackrel{\cong}{\longrightarrow} \calP(1)
\ot \bfk \stackrel{\cong}{\longrightarrow}
\calP(1) \ot \calP(0) \stackrel{\circ_1}{\longrightarrow}
\calP(0) \cong \bfk.
\]
Now take an arbitrary unital associative algebra $A$ and define the operad
$\calA = \coll \calA1$ by
\[
\calA(n) := \cases{A}{for $n = 1$ and}0{for $n \not= 1$,}
\]
with $\circ_1 : \calA(1) \ot \calA(1) \to \calA(1)$ the
multiplication of $A$.
It follows from the above considerations that $\calA = \ccalP$ for some
operad $\calP$ with $\calP(0) \cong \bfk$ if and only if $A$ admits an
augmentation. Therefore any unital associative algebra that does not admit an
augmentation produces the desired example.
}
\end{example}

\begin{example}
{\rm
Kernels, images,~\&c., of homomorphisms between operads in
the category of $\bfk$-modules are defined componentwise. For example,
if $f : \calP \to \calQ$ is such homomorphism, then
$\Ker(f) = \coll{\Ker(f)}0$ is the collection with
\[
Ker(f)(n)
:= \Ker\left(\rr f : \calP(n) \to \calQ(n)
\right),\ n \geq 0.
\]
It is clear that $\Ker(f)$ is an ideal in $\calP$.

Also quotients are defined componentwise. If $\calI$ is an ideal in
$\calP$, then the collection $\calP/\calI = \coll{(\calP/\calI)}0$
with $(\calP/\calI)(n) := \calP(n)/\calI(n)$ for $n \geq 0$, has a
natural operad structure induced by the structure of $\calP$. The
canonical projection $\calP \to \calP/\calI$ has the expected
universal property. The kernel of this projection equals $\calI$.
}
\end{example}

Sometimes it suffices to consider operads without the symmetric group
action. This notion is formalized by:

\begin{definition}[May's non-$\Sigma$ operad]
\label{b1}
A {\em non-$\Sigma$ operad\/} in the category of\/ $\bfk$-modules is a
collection $\ucalP =\{\ucalP(n)\}_{n\geq 0}$ of\/ $\bfk$-modules,
together with operadic compositions
\[
\ugamma : \ucalP(n) \ot \ucalP(k_1) \ot \cdots \ot \ucalP(k_n)
\to
\ucalP(k_1 + \cdots + k_n),
\]
for $n \geq 1$ and $k_1,\ldots,k_n \geq 0$, and a unit map $\ueta :
\bfk \to \ucalP(1)$ that fulfill the associativity and unitality
axioms of Definition~\ref{a}.
\end{definition}

Each operad can be considered as a non-$\Sigma$ operad by forgetting
the $\Sigma_n$-actions. On the other hand, given a \ns\ operad
${\ucalP}$, there is an associated operad $\Sigma[\ucalP]$ with
$\Sigma[\ucalP](n) := \ucalP(n)\ot \bfk[\Sigma_n]$, $n \geq 0$,
with the structure operations induced by the structure operations of
$\underline{\calP}$. Operads of this form are sometimes called {\em
regular\/} operads.

\begin{example}
\label{wx}
{\rm
Consider the operad $\Com = \coll\Com1$ such that $\Com(n) :=
\bfk$ with the trivial $\Sigma_n$-action, $n \geq 1$, and the operadic
compositions~(\ref{pocasi_je_na_nic}) given by the canonical
identifications
\[
\Com(n) \ot \Com(k_1) \ot \cdots \ot \Com(k_n) \cong \bfk^{\ot (n+1)}
\stackrel{\cong}\longrightarrow \bfk \cong
\Com(k_1 + \cdots + k_n).
\]
The operad $\Com$ is obviously not regular. Observe also that $\Com
\cong \widehat{\End}_\bfk$, where $\widehat{\End}_\bfk$ is the endomorphism
operad of the ground ring without the initial component, see
Example~\ref{zitra_strasne_vedro} for the notation.

Let $\uAss$ denote the operad $\Com$ considered as a \ns\
operad. Its symmetrization $\Sigma[\uAss]$ then equals the operad
$\Associative$ introduced in Example~\ref{9}.
}
\end{example}

As we already observed, there is an alternative approach to
operads. For the purposes of comparison, in the rest of this
section and in the following section we will refer to operads
viewed from this alternative perspective as to Markl's operads.
See also the remarks on the terminology in the introduction.

\begin{definition}
\label{b3}
A {\em Markl's operad\/} in the category of\/ $\bfk$-modules is a collection
$\calS =\{\calS(n)\}_{n\geq 0}$ of right $\bfk[\Sigma_n]$-modules,
together with $\bfk$-linear maps ($\circ_i$-compositions)
\[
\circ_i : \calS(m) \ot \calS(n) \to \calS(m +n -1),
\]
for $1 \leq i \leq m$ and $n \geq 0$.
These data fulfill the following axioms.

\noindent
{\em Associativity.\/} For each $1 \leq j \leq a$, $b,c \geq 0$, $f\in
\calS(a)$, $g \in \calS(b)$ and $h \in \calS(c)$,
\[
(f \circ_j g)\circ_i h =
\left\{
\begin{array}{ll}
(f \circ_i h) \circ_{j+c-1}g,& \mbox{for } 1\leq i< j,
\\
f \circ_j(g \circ_{i-j+1} h),& \mbox{for }
j\leq i~< b+j, \mbox{ and}
\\
(f \circ_{i-b+1}h) \circ_j g,& \mbox{for }
j+b\leq i\leq a+b-1,
\end{array}
\right.
\]
see  Figure~\ref{axiom}.
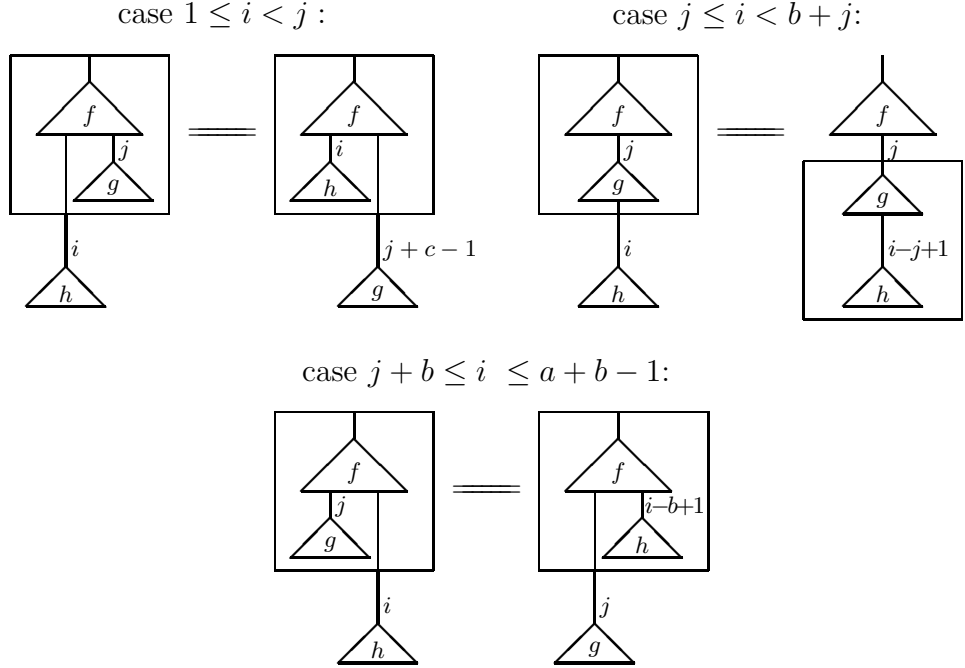
\begin{figure}[t]
\begin{center}
\begin{picture}(230.00,240.00)(40.00,0.00)

\put(60,245){\makebox(0,0){\rm case $1 \leq i < j:$}}
\put(255,245){\makebox(0,0){\rm case $j \leq i < b+j$:}}
\put(-100,0){
\put(160.00,200.00){\makebox(0.00,0.00)%
          {$=\hskip-2mm =\hskip-2mm =\hskip-2mm =\hskip-2mm =\hskip-2mm =$}}


\thicklines
\put(130.00,200.00){\line(-1,1){20.00}}
\put(90.00,200.00){\line(1,0){40.00}}
\put(110.00,220.00){\line(-1,-1){20.00}}
\put(110.00,230.00){\line(0,-1){10.00}}
\put(110.00,211.00){\makebox(0.00,0.00)[t]{\scriptsize $f$}}

\put(9,0){
\put(110.00,200.00){\line(0,-1){10.00}}
\put(110.00,183.00){\makebox(0.00,0.00)[t]{\scriptsize $g$}}
\put(112.00,198.00){\makebox(0.00,0.00)[lt]{\scriptsize $j$}}
\put(95.00,175.00){\line(1,0){30.00}}
\put(125.00,175.00){\line(-1,1){15.00}}
\put(110.00,190.00){\line(-1,-1){15.00}}
}


\thinlines
\put(130.00,170.00){\line(-1,0){40.00}}
\put(90.00,230.00){\line(1,0){40.00}}
\put(80.00,170.00){\line(1,0){10.00}}
\put(80.00,230.00){\line(0,-1){60.00}}
\put(90.00,230.00){\line(-1,0){10.00}}
\put(140.00,170.00){\line(-1,0){10.00}}
\put(140.00,230.00){\line(0,-1){60.00}}
\put(130.00,230.00){\line(1,0){10.00}}

\thicklines

\put(-9,0){
\put(110.00,143.00){\makebox(0.00,0.00)[t]{\scriptsize $h$}}
\put(112.00,160.00){\makebox(0.00,0.00)[tl]{\scriptsize $i$}}
\put(125.00,135.00){\line(-1,1){15.00}}
\put(95.00,135.00){\line(1,0){30.00}}
\put(110.00,150.00){\line(-1,-1){15.00}}
\put(110.00,170.00){\line(0,-1){20.00}}
\thinlines
\put(110.00,170.00){\line(0,1){30.00}}
}
}


\thicklines
\put(130.00,200.00){\line(-1,1){20.00}}
\put(90.00,200.00){\line(1,0){40.00}}
\put(110.00,220.00){\line(-1,-1){20.00}}
\put(110.00,230.00){\line(0,-1){10.00}}
\put(110.00,211.00){\makebox(0.00,0.00)[t]{\scriptsize $f$}}

\put(-9,0){
\put(110.00,200.00){\line(0,-1){10.00}}
\put(110.00,183.00){\makebox(0.00,0.00)[t]{\scriptsize $h$}}
\put(112.00,198.00){\makebox(0.00,0.00)[lt]{\scriptsize $i$}}
\put(95.00,175.00){\line(1,0){30.00}}
\put(125.00,175.00){\line(-1,1){15.00}}
\put(110.00,190.00){\line(-1,-1){15.00}}
}


\thinlines
\put(130.00,170.00){\line(-1,0){40.00}}
\put(90.00,230.00){\line(1,0){40.00}}
\put(80.00,170.00){\line(1,0){10.00}}
\put(80.00,230.00){\line(0,-1){60.00}}
\put(90.00,230.00){\line(-1,0){10.00}}
\put(140.00,170.00){\line(-1,0){10.00}}
\put(140.00,230.00){\line(0,-1){60.00}}
\put(130.00,230.00){\line(1,0){10.00}}

\thicklines

\put(9,0){
\put(110.00,143.00){\makebox(0.00,0.00)[t]{\scriptsize $g$}}
\put(112.00,160.00){\makebox(0.00,0.00)[tl]{\scriptsize $j+c-1$}}
\put(125.00,135.00){\line(-1,1){15.00}}
\put(95.00,135.00){\line(1,0){30.00}}
\put(110.00,150.00){\line(-1,-1){15.00}}
\put(110.00,170.00){\line(0,-1){20.00}}
\thinlines
\put(110.00,170.00){\line(0,1){30.00}}
}

\put(100,0){


\put(160.00,200.00){\makebox(0.00,0.00)%
          {$=\hskip-2mm =\hskip-2mm =\hskip-2mm =\hskip-2mm =\hskip-2mm =$}}


\thicklines
\put(130.00,200.00){\line(-1,1){20.00}}
\put(90.00,200.00){\line(1,0){40.00}}
\put(110.00,220.00){\line(-1,-1){20.00}}
\put(110.00,230.00){\line(0,-1){10.00}}
\put(110.00,211.00){\makebox(0.00,0.00)[t]{\scriptsize $f$}}


\put(110.00,200.00){\line(0,-1){10.00}}
\put(110.00,183.00){\makebox(0.00,0.00)[t]{\scriptsize $g$}}
\put(112.00,198.00){\makebox(0.00,0.00)[lt]{\scriptsize $j$}}
\put(95.00,175.00){\line(1,0){30.00}}
\put(125.00,175.00){\line(-1,1){15.00}}
\put(110.00,190.00){\line(-1,-1){15.00}}


\thinlines
\put(130.00,170.00){\line(-1,0){40.00}}
\put(90.00,230.00){\line(1,0){40.00}}
\put(80.00,170.00){\line(1,0){10.00}}
\put(80.00,230.00){\line(0,-1){60.00}}
\put(90.00,230.00){\line(-1,0){10.00}}
\put(140.00,170.00){\line(-1,0){10.00}}
\put(140.00,230.00){\line(0,-1){60.00}}
\put(130.00,230.00){\line(1,0){10.00}}

\thicklines

\put(110.00,143.00){\makebox(0.00,0.00)[t]{\scriptsize $h$}}
\put(112.00,160.00){\makebox(0.00,0.00)[tl]{\scriptsize $i$}}
\put(125.00,135.00){\line(-1,1){15.00}}
\put(95.00,135.00){\line(1,0){30.00}}
\put(110.00,150.00){\line(-1,-1){15.00}}
\put(110.00,170.00){\line(0,-1){20.00}}
\thinlines
\put(110.00,170.00){\line(0,1){5.00}}
}

\put(200,0){



\thicklines
\put(130.00,200.00){\line(-1,1){20.00}}
\put(90.00,200.00){\line(1,0){40.00}}
\put(110.00,220.00){\line(-1,-1){20.00}}
\put(110.00,230.00){\line(0,-1){10.00}}
\put(110.00,211.00){\makebox(0.00,0.00)[t]{\scriptsize $f$}}


\put(0,-5){
\put(110.00,205.00){\line(0,-1){15.00}}
\put(110.00,183.00){\makebox(0.00,0.00)[t]{\scriptsize $g$}}
\put(112.00,203.00){\makebox(0.00,0.00)[lt]{\scriptsize $j$}}
\put(95.00,175.00){\line(1,0){30.00}}
\put(125.00,175.00){\line(-1,1){15.00}}
\put(110.00,190.00){\line(-1,-1){15.00}}
}


\put(0,-40){
\thinlines
\put(130.00,170.00){\line(-1,0){40.00}}
\put(90.00,230.00){\line(1,0){40.00}}
\put(80.00,170.00){\line(1,0){10.00}}
\put(80.00,230.00){\line(0,-1){60.00}}
\put(90.00,230.00){\line(-1,0){10.00}}
\put(140.00,170.00){\line(-1,0){10.00}}
\put(140.00,230.00){\line(0,-1){60.00}}
\put(130.00,230.00){\line(1,0){10.00}}
}

\thicklines

\put(110.00,143.00){\makebox(0.00,0.00)[t]{\scriptsize $h$}}
\put(112.00,160.00){\makebox(0.00,0.00)[tl]{%
  \scriptsize $i\hskip -1mm -\hskip -1mm j\hskip -1mm + \hskip -1mm 1$}}
\put(125.00,135.00){\line(-1,1){15.00}}
\put(95.00,135.00){\line(1,0){30.00}}
\put(110.00,150.00){\line(-1,-1){15.00}}
\put(110.00,170.00){\line(0,-1){20.00}}
}

\put(160,110){\makebox(0.00,0.00){\rm case $j + b \leq i~\leq a+b- 1$:}}

\put(0,-135){

\put(160.00,200.00){\makebox(0.00,0.00)%
          {$=\hskip-2mm =\hskip-2mm =\hskip-2mm =\hskip-2mm =\hskip-2mm =$}}


\thicklines
\put(130.00,200.00){\line(-1,1){20.00}}
\put(90.00,200.00){\line(1,0){40.00}}
\put(110.00,220.00){\line(-1,-1){20.00}}
\put(110.00,230.00){\line(0,-1){10.00}}
\put(110.00,211.00){\makebox(0.00,0.00)[t]{\scriptsize $f$}}

\put(-9,0){
\put(110.00,200.00){\line(0,-1){10.00}}
\put(110.00,183.00){\makebox(0.00,0.00)[t]{\scriptsize $g$}}
\put(112.00,198.00){\makebox(0.00,0.00)[lt]{\scriptsize $j$}}
\put(95.00,175.00){\line(1,0){30.00}}
\put(125.00,175.00){\line(-1,1){15.00}}
\put(110.00,190.00){\line(-1,-1){15.00}}
}


\thinlines
\put(130.00,170.00){\line(-1,0){40.00}}
\put(90.00,230.00){\line(1,0){40.00}}
\put(80.00,170.00){\line(1,0){10.00}}
\put(80.00,230.00){\line(0,-1){60.00}}
\put(90.00,230.00){\line(-1,0){10.00}}
\put(140.00,170.00){\line(-1,0){10.00}}
\put(140.00,230.00){\line(0,-1){60.00}}
\put(130.00,230.00){\line(1,0){10.00}}

\thicklines

\put(9,0){
\put(110.00,143.00){\makebox(0.00,0.00)[t]{\scriptsize $h$}}
\put(112.00,160.00){\makebox(0.00,0.00)[tl]{\scriptsize $i$}}
\put(125.00,135.00){\line(-1,1){15.00}}
\put(95.00,135.00){\line(1,0){30.00}}
\put(110.00,150.00){\line(-1,-1){15.00}}
\put(110.00,170.00){\line(0,-1){20.00}}
\thinlines
\put(110.00,170.00){\line(0,1){30.00}}
}
}

\put(100,-135){

\thicklines
\put(130.00,200.00){\line(-1,1){20.00}}
\put(90.00,200.00){\line(1,0){40.00}}
\put(110.00,220.00){\line(-1,-1){20.00}}
\put(110.00,230.00){\line(0,-1){10.00}}
\put(110.00,211.00){\makebox(0.00,0.00)[t]{\scriptsize $f$}}

\put(9,0){
\put(110.00,200.00){\line(0,-1){10.00}}
\put(110.00,183.00){\makebox(0.00,0.00)[t]{\scriptsize $h$}}
\put(111.00,198.00){\makebox(0.00,0.00)[lt]{%
       \scriptsize $i \hskip -1mm  - \hskip -1mm b\hskip -1mm +\hskip -1mm 1$}}
\put(95.00,175.00){\line(1,0){30.00}}
\put(125.00,175.00){\line(-1,1){15.00}}
\put(110.00,190.00){\line(-1,-1){15.00}}
}


\thinlines
\put(130.00,170.00){\line(-1,0){40.00}}
\put(90.00,230.00){\line(1,0){40.00}}
\put(80.00,170.00){\line(1,0){10.00}}
\put(80.00,230.00){\line(0,-1){60.00}}
\put(90.00,230.00){\line(-1,0){10.00}}
\put(144.00,170.00){\line(-1,0){14.00}}
\put(144.00,230.00){\line(0,-1){60.00}}
\put(130.00,230.00){\line(1,0){14.00}}

\thicklines

\put(-9,0){
\put(110.00,143.00){\makebox(0.00,0.00)[t]{\scriptsize $g$}}
\put(112.00,160.00){\makebox(0.00,0.00)[tl]{\scriptsize $j$}}
\put(125.00,135.00){\line(-1,1){15.00}}
\put(95.00,135.00){\line(1,0){30.00}}
\put(110.00,150.00){\line(-1,-1){15.00}}
\put(110.00,170.00){\line(0,-1){20.00}}
\thinlines
\put(110.00,170.00){\line(0,1){30.00}}
}
}
\end{picture}
\end{center}
\caption{%
\label{axiom}
Flow charts explaining the associativity in Markl's operads.
}
\end{figure}

\noindent
{\em Equivariance.\/}
For each  $1 \leq i \leq m$, $n \geq 0$, $\tau \in \Sigma_m$ and
$\sigma \in \Sigma_n$,  let $\tau \circ_i
\sigma \in \Sigma_{m+n -1}$ be given by inserting the permutation
$\sigma$ at the $i$th place in $\tau$. Let $f \in \calS(m)$ and $g \in
\calS(n)$. Then
\[
(f\tau)\circ_i(g\sigma) = (f\circ_{\tau(i)} g)(\tau \circ_i
\sigma).
\]
\noindent
{\em Unitality.\/}
There exists
$e \in \calS(1)$ such that
\begin{equation}
\label{Jitka_S_zase_nepise}
f\circ_i e = e\ \mbox { and }\ e \circ_1 g = g
\end{equation}
for each  $1\leq i \leq m$, $n \geq 0$, $f \in \calS(m)$ and $g \in \calS(n)$.
\end{definition}

\begin{example}
{\rm
All axioms in Definition~\ref{b3} can be read
from the endomorphism operad
$\End_V = \{\End_V(n)\}_{n \geq 0}$ of a $\bfk$-module~$V$ reviewed in
Example~\ref{end}, with $\circ_i$-operations given by
\[
f \circ_i g := f ( \id^{\otimes i-1}_V \ot g \ot \id^{\otimes m-1}_V),
\]
for $f \in \End_V(m)$, $g \in \End_V(n)$, $1 \leq i \leq m$ and $n
\geq 0$.
}
\end{example}

The following proposition shows that Definition~\ref{a} describes the
same objects as Definition~\ref{b3}.

\begin{proposition}
\label{cc}
The category of May's operads is
isomorphic to the category of Markl's operads.
\end{proposition}

\noindent
{\bf Proof.}
Given a Markl's operad $\calS = \coll{\calS}0$ as in
Definition~\ref{b3}, define a May's
operad $\calP=\May(\calS)$ by $\calP(n) := \calS(n)$ for $n \geq 0$,
with the $\gamma$-operations given by
\begin{equation}
\label{a1}
\gamma(f,\Rada g1n) :=
(\cdots((f\circ_n g_n)\circ_{n-1} g_{n-1}) \cdots) \circ_1g_1
\end{equation}
where $f \in \calP(n)$, $g_i \in \calP(k_i)$, $1 \leq i \leq n$,
$\Rada k1n \geq 0$. The unit morphism $\eta : \bfk \to \calP(1)$ is
defined by $\eta(1) := e$.
It is easy to verify that $\May(-)$ extends to a functor
from the category of Markl's operads  the
category of May's operads.

On the other hand, given a May's operad $\calP$, one
can define a Markl's operad $\calS =\Mar(\calP)$ by $\calS(n) :=
\calP(n)$ for $n\geq 0$, with the $\circ_i$-operations:
\begin{equation}
\label{a2}
f\circ_i g :=
\gamma(f,\underbrace{e,\ldots,e}_{i-1},g,\underbrace{e,\ldots,e}_{m-i}),
\end{equation}
for $f \in \calS(m)$, $g \in \calS(n)$, $m \geq 1$, $n \geq 0$,
where $e := \eta(1) \in \calP(1)$.
It is again obvious that $\Mar(-)$ extends to a functor
that the functors
$\May(-)$ and $\Mar(-)$
are mutually inverse isomorphisms between the category of
Markl's operads and the category of May's operads.%
\qed

The equivalence between May's and Markl's operads implies
that an operad can be defined by specifying $\circ_i$-operations and
a unit. This is sometimes simpler that to define the
$\gamma$-operations directly, as illustrated by:

\begin{example}
\label{Jitka_se_ozvala}
{\rm
Let $\Sigma$ be a Riemann sphere, that
is, a nonsingular complex projective curve of genus $0$. By a
puncture or a parametrized hole
we mean a point $p$ of $\Sigma$ together with a
holomorphic
embedding of the standard closed
disc $U = \{z \in {\mathbb C}\,;\ |z| \leq 1\}$ to
$\Sigma$ centered at the point. Thus a puncture is
a holomorphic embedding $u: \tilde U \to \Sigma$,
where $\tilde U \subset {\mathbb C}$ is an open neighborhood of $U$ and
$u(0) = p$. We say that two punctures $u_1: \tilde U_1 \to \Sigma$  and
$u_2 : \tilde U_2 \to \Sigma$ are disjoint, if
\[
u_1(\intU)
\cap u_2(\intU) = \emptyset,
\]
where $\intU : = \{z \in {\mathbb C}\,;\ |z| < 1\}$ is the interior of $U$.

Let $\hatcalM_0(n)$ be the moduli space of Riemann spheres
$\Sigma$ with $n+1$ disjoint punctures $u_i : \tilde U_i \to \Sigma$,
$0 \leq i \leq n$, modulo the action of complex projective
automorphisms.  The topology of $\hatcalM_0(n)$ is a very
subtle thing and we are not going to discuss this issue here;
see~\cite{huang:book}. The constructions below will be made only `up
to topology.'

Renumbering the holes $\Rada u1n$ defines on each $\hatcalM_0(n)$ a
natural right $\Sigma_n$-action and the $\Sigma$-module
$\hatcalM_0 = \coll{\hatcalM_0}0$ forms a
topological operad under sewing Riemannian spheres at punctures.  Let
us describe this operadic structure using the
$\circ_i$-formalism. Thus, let $\Sigma$ represent an element $x \in
\hatcalM_0(m)$ and $\Delta$ represent an element $y \in
\hatcalM_0(n)$. For $ 1 \leq i \leq m$, let $u_i: \tilde
U_i \to \Sigma$ be the $i$th puncture of $\Sigma$ and let $u_0:\tilde
U_0 \to \Delta$ be the $0$th puncture of $\Delta$.

There certainly exists some $0 < r < 1$ such that both $\tilde U_0$
and $\tilde U_i$ contain the disc $U_{1/r} : =
\{\hbox{$z \in {\mathbb C}$}\,;\
|z| <1/r\}$. Let now $\Sigma_r: = \Sigma \setminus u_i(U_{r})$ and
$\Delta_r: = \Delta \setminus u_0(U_{r})$. Define finally
\[
\Xi:= (\Sigma_r \bigsqcup \Delta_r)/\sim,
\]
where the relation $\sim$  is given by
\[
\Sigma_r \ni u_i(\xi) \sim u_0(1/\xi) \in \Delta_r,
\]
for $r < |\xi| < 1/r$. It is immediate to see that $\Xi$ is a
well-defined punctured Riemannian sphere, with $n+m-1$ punctures
induced in the obvious manner from those of $\Sigma$ and $\Delta$, and
that the class of the punctured surface $\Xi$ in the moduli space
$\hatcalM_0(m+n-1)$ does not depend on the representatives
$\Sigma$, $\Delta$ and on $r$.  We define $x \circ _i y$ to be the
class of $\Xi$.

The unit $e \in\hatcalM_0(1)$ can be defined as follows.
Let $\CPj$ be the complex projective line with homogeneous
coordinates $[z,w]$, $z,w \in {\mathbb C}$,~\cite[Example~I.1.6]{wells:book}.
Let $0 := [0,1]\in \CPj$ and $\infty:=[1,0]
\in \CPj$. Recall that we have two canonical isomorphisms
$p_\infty : \CPj \setminus \infty \to {\mathbb C}$ and
$p_0 : \CPj \setminus 0 \to {\mathbb C}$ given~by
\[
p_\infty ([z,w]) := z/w \mbox { and } p_0 ([z,w]) := w/z.
\]
Then $p_\infty^{-1} : {\mathbb C} \to \CPj$ (respectively $p_0^{-1} :
{\mathbb C} \to \CPj$) is a puncture at $0$ (respectively at $\infty$).
We define $e \in \hatcalM_0(1)$ to be the class of
$(\CPj,p_0^{-1}, p_{\infty}^{-1})$.

It is not hard to verify that the above
constructions make the collection $\hatcalM_0 =
\coll{\hatcalM_0}0$ a Markl's operad. By
Proposition~\ref{cc}, $\hatcalM_0$ is a also
May's operad.
}
\end{example}

In the rest of this article, we will consider
May's and Markl's operads as two versions of
the same object which we will call
simply a (unital) operad.

\section{Non-unital operads}
\label{sec1bis}

It turns out that the combinatorial structure of
the moduli space of stable genus zero curves is captured by a certain
non-unital version of operad.  Let $\calM_{0,n+1}$
be the moduli space of \hbox{$(n+1)$}-tuples $(\Rada x0n)$ of distinct
numbered points on the complex projective line $\CPj$ modulo
projective automorphisms, that is, transformations of the form
\[
\CPj \ni [{\xi}_1,{\xi}_2]
\mapsto [a{\xi}_1+b{\xi}_2,c{\xi}_1 + d{\xi}_2] \in \CPj,
\]
where $a,b,c,d \in {\mathbb C}$ with $ad-bc \not=0$.

The moduli space $\calM_{0,n+1}$ has, for $n \geq 2$, a canonical
compactification $\barcalM_0(n) \supset \calM_{0,n+1}$ introduced by
A.~Grothendieck and
F.F.~Knudsen~\cite{deligne:LNM288,knudsen:MS83}. The space
$\barcalM_0(n)$ is the moduli space of stable $(n+1)$-pointed curves of
genus $0$:

\begin{definition}
\label{jItka}
A {\em stable $(n+1)$-pointed curve of genus $0$} is an object
\[
(C;\Rada x0n),
\]
where $C$ is a (possibly reducible) algebraic curve
with at most nodal singularities and $\Rada x0n \in C$ are distinct
smooth points such that
\begin{itemize}
\item[(i)]
each component of $C$ is isomorphic to $\CPj$,
\item[(ii)]
the graph of intersections of components of $C$ (i.e.~the graph whose
vertices correspond to the components of $C$ and edges to the intersection
points of the components) is a tree, and
\item[(iii)]
each component of $C$ has at least three special points, where a special
point means either one of the $x_i$, $0 \leq i \leq n$,
or a singular point of $C$ (the stability).
\end{itemize}
\end{definition}

It can be easily seen that a stable curve $(C;\Rada x0n)$ admits no
infinitesimal automorphisms that fix marked points $\Rada x0n$,
therefore $(C;\Rada x0n)$ is `stable' in the usual sense.
Observe also that $\barcalM_0(0) =\barcalM_0(1)= \emptyset$ (there are
no stable curves with less than three marked points) and that
$\barcalM_0(2) = \mbox {\rm the point}$
corresponding to the three-pointed stable curve $(\CPj; \infty,1,0)$.
The space $\calM_{0,n+1}$ forms an open dense part of $\barcalM_0(n)$
consisting of marked curves $(C;\Rada x0n)$ such that $C$ is
isomorphic to $\CPj$.

Let us try to equip the collection $\barcalM_0 = \{\barcalM_0(n)\}_{n \geq
2}$ with an
operad structure as in Definition~\ref{a}.\label{New_Orleans}
For $C = (C,\Rada x1n) \in
\barcalM_0(n)$ and $C_i = (C_i,\Rada {y^i}1{k_i}) \in \barcalM_0(k_i)$, $1
\leq i \leq n$, let
\begin{equation}
\label{bude_vedro_a_cista_termika}
\gamma(C,C_1,\ldots,C_n) \in \barcalM_0(k_1+\cdots
+ k_n)
\end{equation}
be the stable marked curve obtained from the disjoint union
$C\sqcup C^1 \sqcup \cdots \sqcup C^n$ by identifying, for each $1
\leq i \leq n$, the point $x_i
\in C$ with the point $y_0^i \in C_i$, introducing a nodal
singularity, and relabeling the remaining marked points
accordingly. The symmetric group acts on $\barcalM_0(n)$ by
\[
(C,x_0,x_1,\ldots,x_n) \longmapsto
(C,x_0,x_{\sigma(1)},\ldots,x_{\sigma(n)}),\
\sigma \in \Sigma_n.
\]

We have defined the $\gamma$-compositions and the symmetric group
action, but there is no room for the identity, because
$\barcalM_0(1)$ is empty!  The above structure is, therefore,
a non-unital operad in the sense of the following definition (which is
formulated, as all definitions in this article, for the monoidal
category of $\bfk$-modules).

\begin{definition}
A {\em May's non-unital operad\/} in the category of\/ $\bfk$-modules is a
collection $\calP =\{\calP(n)\}_{n\geq 0}$ of\/ $\bfk[\Sigma_n]$-modules,
together with operadic compositions
\[
\gamma : \calP(n) \ot \calP(k_1) \ot \cdots \ot \calP(k_n)
\to
\calP(k_1 + \cdots + k_n),
\]
for $n \geq 1$ and $k_1,\ldots,k_n \geq 0$, that fulfill the
associativity and equivariance axioms of Definition~\ref{a}.
\end{definition}

We may as well define
on the collection $\barcalM_0 =
\{\barcalM_0(n)\}_{n \geq 2}$  operations
\begin{equation}
\label{dead_or_alive}
\circ_i:
\barcalM_0(m) \times \barcalM_0(n)
\to \barcalM_0(m+n-1)
\end{equation}
for $m,n \geq  2$, $1 \leq i \leq m$, by
\[
(C^1;y_0,\ldots,y_m) \times
(C^2;\rada{x_0}{x_{n}})
\longmapsto
(C;\Rada y0{i-1},\Rada x0n,\Rada y{i+1}m)
\]
where $C$ is the quotient of the disjoint union $C^1 \bigsqcup
C^2$ given by identifying $x_0$ with $y_i$ at a nodal singularity,
see Figure~\ref{zitra_vecer_je_prilet}.
\begin{figure}
\begin{center}
{
\unitlength=.83pt
\begin{picture}(210.00,120.00)(0.00,0.00)
\put(160.00,40.00){\makebox(0.00,0.00)[tl]{$y_i=x_0$}}
\put(151.00,48.50){\makebox(0.00,0.00){$\bullet$}}
\put(210.00,120.00){\makebox(0.00,0.00)[l]{$C^2$}}
\put(10.00,10.00){\makebox(0.00,0.00){$C^1$}}
\thicklines
\bezier{100}(180.00,120.00)(190.00,110.00)(210.00,100.00)
\bezier{100}(170.00,100.00)(180.00,80.00)(200.00,70.00)
\bezier{98}(150.00,40.00)(150.00,20.00)(170.00,0.00)
\bezier{98}(180.00,90.00)(190.00,100.00)(200.00,120.00)
\put(160.00,70.00){\line(1,1){20.00}}
\bezier{98}(150.00,40.00)(150.00,60.00)(160.00,70.00)
\bezier{98}(70.00,90.00)(80.00,70.00)(90.00,70.00)
\bezier{150}(60.00,70.00)(90.00,80.00)(100.00,90.00)
\put(30.00,60.00){\line(3,1){30.00}}
\put(50.00,50.00){\line(2,-3){20.00}}
\bezier{98}(40.00,90.00)(40.00,70.00)(50.00,50.00)
\bezier{98}(0.00,40.00)(20.00,30.00)(30.00,10.00)
\bezier{150}(110.00,40.00)(140.00,50.00)(180.00,50.00)
\bezier{98}(60.00,30.00)(80.00,30.00)(110.00,40.00)
\put(40.00,30.00){\line(1,0){20.00}}
\bezier{98}(0.00,20.00)(20.00,30.00)(40.00,30.00)
\end{picture}}
\end{center}
\caption{\label{zitra_vecer_je_prilet}
The $\circ_i$-compositions in $\barcalM_0 = \{\barcalM_0(n)\}_{n\geq 2}$.}
\end{figure}
The collection $\barcalM_0 = \coll{\barcalM_0}2$ with
$\circ_i$-operations~(\ref{dead_or_alive}) is an example of another
version of non-unital operads, recalled in:

\begin{definition}
A {\em non-unital Markl's operad\/} in the category of\/ $\bfk$-modules is a
collection $\calP =\{\calP(n)\}_{n\geq 0}$ of\/ $\bfk[\Sigma_n]$-modules,
together with operadic compositions
\[
\circ_i : \calS(m) \ot \calS(n) \to \calS(m +n -1),
\]
for $1 \leq i \leq m$ and $n \geq 0$,
that fulfill the
associativity and equivariance axioms of Definition~\ref{b3}.
\end{definition}

A we saw in Proposition~\ref{cc}, in the presence of
operadic units, May's operads are the same as Markl's
operads. Surprisingly,
the non-unital versions of these structures are
{\em radically different\/} --
Markl's operads capture more information than May's operads! This is
made precise in the following:

\begin{proposition}
\label{cd}
The category of non-unital Markl's operads is a subcategory of
the category of non-unital May's operads.
\end{proposition}

\noindent
{\bf Proof.}
It is easy to see that~(\ref{a1}) defines, as in the proof of
Proposition~\ref{cc}, a functor
$\psi\May(-)$
which is an embedding of the category of non-unital Markl's operads into the
category of  non-unital May's operads.%
\qed

Observe that
formula~(\ref{a2}), inverse to~(\ref{a1}), does not make sense without units.
The relation between various versions of operads discussed so far
is summarized in the following diagram of categories
and their inclusions:

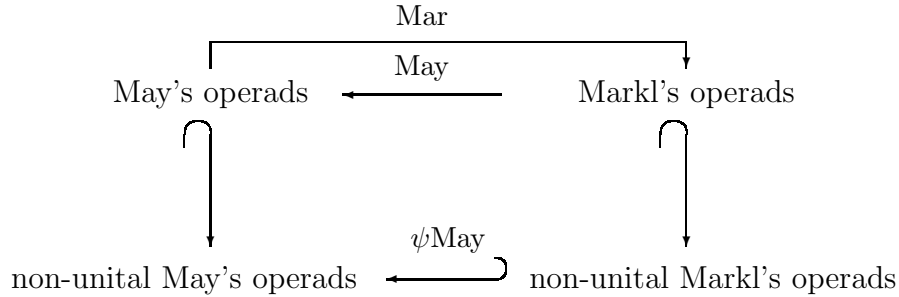
\begin{center}
{
\unitlength=1.000000pt
\begin{picture}(140.00,125.00)(0.00,-10.00)
\put(70.00,100.00){\makebox(0.00,0.00){\small $\Mar$}}
\put(70.00,80.00){\makebox(0.00,0.00){\small $\May$}}
\put(80.00,15.00){\makebox(0.00,0.00){\small $\psi\May$}}
\put(0,5){
\qbezier(165.00,55.00)(170.00,55.00)(170.00,50.00)
\qbezier(165.00,55.00)(160.00,55.00)(160.00,50.00)
\put(160.00,50.00){\line(0,-1){5}}
\put(-180,0){
\qbezier(165.00,55.00)(170.00,55.00)(170.00,50.00)
\qbezier(165.00,55.00)(160.00,55.00)(160.00,50.00)
\put(160.00,50.00){\line(0,-1){5}}
}
}
\put(-10.00,90.00){\line(0,-1){10.00}}
\put(170.00,90.00){\line(-1,0){180.00}}
\put(170.00,90.00){\vector(0,-1){10.00}}
\put(100.00,70.00){\vector(-1,0){60.00}}
\put(97.5,00){$\longhookleftarrow$}
\put(170.00,55.00){\vector(0,-1){42.50}}
\put(-10.00,55.00){\vector(0,-1){42.50}}
\put(180.00,0.00){\makebox(0.00,0.00){non-unital Markl's operads}}
\put(-20,0.00){\makebox(0.00,0.00){non-unital May's operads}}
\put(170.00,70.00){\makebox(0.00,0.00){Markl's operads}}
\put(-10,70.00){\makebox(0.00,0.00){May's operads}}
\end{picture}}
\end{center}

The following example shows that non-unital Markl's operads form a
proper sub-category of the category of non-unital May's operads.

\begin{example}
{\rm
We describe a non-unital May's operad $\calV =
\coll \calV0$ which is not of the form $\psi\May(\calS)$ for some
non-unital Markl's operad $\calS$.  Let
\[
\calV(n) := \cases{\bfk}{for $n =2$ or $4$, and}0{otherwise.}
\]
The only non-trivial $\gamma$-composition is $\gamma: \calV(2) \ot
\calV(2) \ot\calV(2) \to \calV(4)$, given as the canonical isomorphism
\[
\calV(2) \ot \calV(2) \ot\calV(2) \cong \bfk^{\ot 3}
\stackrel{\cong}{\longrightarrow} \bfk \cong  \calV(4).
\]

Suppose that $\calV = \May(\calS)$ for some non-unital Markl's
operad $\calS$. Then,
according to~(\ref{a1}), for $f,g_1,g_2\in \calV(2)$,
\[
\gamma(f,g_1,g_2) = (f \circ_2 g_2) \circ_1 g_1.
\]
Since $(f \circ_2 g) \in \calV(3) = 0$, this would imply that $\gamma$ is
trivial, which is not true.
}
\end{example}

Proposition~\ref{toto_pisu_na_zavodech_v_Moravske_Trebove_2005} below
shows that Markl's rather than May's non-unital operads are true
non-unital versions of operads. We will need the following definition
in which $\calK = \coll\calK1$ is the trivial (unital) operad with $\calK
(1) := \bfk$ and $\calK(n) = 0$, for $n \not= 1$.

\begin{definition}
An {\em augmentation\/} of an operad $\calP$ in the category of\/
$\bfk$-modules is a homomorphism $\epsilon : \calP \to \calK$. Operads
with an augmentation are called {\em augmented\/} operads. The kernel
\[
\barcalP : = \Ker\left( \epsilon : \calP \to \calK\right)
\]
is called the {\em augmentation ideal\/}.
\end{definition}

The following proposition was proved in~\cite{markl:zebrulka}.

\begin{proposition}
\label{toto_pisu_na_zavodech_v_Moravske_Trebove_2005}
The correspondence $\calP \mapsto \barcalP$ is an isomorphism between
the category of augmented operads and the category of Markl's
non-unital operads.
\end{proposition}

\noindent
{\bf Proof.}
The $\circ_i$-operations of $\calP$ obviously restrict to $\barcalP$,
making it a non-unital Markl's operad.
It is simple to describe a functorial inverse $\calS \mapsto {\wcalS}$ of the
correspondence $\calP \mapsto \barcalP$.  For a Markl's non-unital
operad $\calS$,
denote by $\wcalS$ the collection\label{uz_je_po_zavodech}
\begin{equation}
\label{Tomicek}
\wcalS (n) := \cases{\calS(n)}{for $n \not =1$, and}
                  {\calS(1) \oplus \bfk}{for $n=1$.}
\end{equation}
The $\circ_i$-operations of $\wcalS$ are uniquely determined by requiring
that they extend the $\circ_i$-operations of $\calS$ and
satisfy~(\ref{Jitka_S_zase_nepise}), with the unit
$e := 0 \oplus 1_\bfk  \in  \calS(1) \oplus \bfk = \wcalS(1)$.
Informally, $\wcalS$ is obtained from the
Markl's non-unital operad $\calS$ by adjoining a unit.%
\qed

Observe that if $\calS$ were a May's, not Markl's, non-unital operad,
the construction of $\wcalS$ described in the above proof would not
make sense, because we would not know how to define
\[
\gamma(f,\underbrace{e,\ldots,e}_{i-1},g,\underbrace{e,\ldots,e}_{m-i})
\]
for $f \in \calS(m)$, $g \in \calS(n)$, $m \geq 2$, $n \geq 0$, $1
\leq i \leq m$.
Proposition~\ref{toto_pisu_na_zavodech_v_Moravske_Trebove_2005} should
be compared to the obvious statement that the category of augmented
unital associative algebras is isomorphic to the category of
(non-unital) associative algebras.  In the following proposition,
$\catOper$ denotes the category of $\bfk$-linear operads and
$\psi\catOper$ the category of $\bfk$-linear Markl's non-unital operads.

\begin{proposition}
Let $\calP$ be an augmented operad and $\calQ$ an arbitrary operad in
the category of $\bfk$-modules.
Then there exists a natural isomorphism
\begin{equation}
\label{probihaji_zavody_a_jsem_v_Praze_protoze_je_strasne_pocasi!!!}
\Mor_{\catOper}(\calP,\calQ)
\cong \Mor_{\psi\catOper}(\barcalP,\psi \May(\calQ)).
\end{equation}
\end{proposition}

The proof is simple and we leave it the reader.
Combining~(\ref{probihaji_zavody_a_jsem_v_Praze_protoze_je_strasne_pocasi!!!})
with the isomorphism of
Proposition~\ref{toto_pisu_na_zavodech_v_Moravske_Trebove_2005}
one obtains a natural isomorphism
\begin{equation}
\label{probihaji_zavody_a_jsem_v_Praze_protoze_je_strasne_pocasi!!}
\Mor_{\catOper}(\wcalS,\calQ)
\cong \Mor_{\psi\catOper}(\calS,\psi \May(\calQ))
\end{equation}
which holds for each Markl's non-unital operad $\calS$ and operad $\calQ$.
Isomorphism~(\ref{probihaji_zavody_a_jsem_v_Praze_protoze_je_strasne_pocasi!!})
means that $\widetilde {\rule {1em}{0em}} \hskip -.5em : \psi\catOper
\to \catOper$ and $\psi\May : \catOper \to \psi\catOper$
are adjoint functors. This
adjunction
will be used in the
construction of free operads in Section~\ref{sec3}.

In the rest of this article, non-unital Markl's operads will be called
simply non-unital operads.
This will not lead to confusion, since all non-unital operads referred
to in the rest of this article will be Markl's.

\section{Operad algebras}
\label{sec2}

As we already remarked, operads are important
through their representations  called operad algebras or simply algebras.

\begin{definition}
\label{jitky}
Let $V$ be a $\bfk$-module and $\End_V$ the endomorphism operad
of $V$ recalled in Example~\ref{end}. A {\em $\calP$-algebra\/} is
a homomorphism of operads $\rho :\calP \to \End_V$.
\end{definition}

The above definition admits an obvious generalization into an
arbitrary symmetric monoidal category with
an internal hom-functor.
The last assumption is necessary for the existence of the
`internal' endomorphism operad,
see~\cite[Definition~II.1.20]{markl-shnider-stasheff:book}.
Definition~\ref{jitky} can be however unwrapped into the form given
in~\cite[Definition~2.1]{kriz-may}
that makes sense in an
arbitrary symmetric monoidal category without the internal hom-functor
assumption:

\begin{proposition}
\label{pro:alg}
Let $\calP$ be an operad.  A
$\calP$-algebra is the same as a $\bfk$-module
$V$ together with maps
\begin{equation}
\label{algebra}
\alpha : \calP(n) \ot V^{\otimes n} \to V,\ n\geq 0,
\end{equation}
that satisfy the following axioms.

\noindent
{\em Associativity.\/}
For each $n \geq 1$ and non-negative integers $k_1,\ldots,k_n$, the following diagram commutes.
\begin{center}
\setlength{\unitlength}{1em}
\begin{picture}(20,13)(0,-1)
\put(0,10){\makebox(0,0){$\displaystyle
\left(
\calP(n) \ot \bigotimes_{s=1}^n \calP(k_s) \right)
\ot \bigotimes_{s=1}^n V^{\otimes k_s}
$}}
\put(7.5,10){\vector(1,0){4}}
\put(9.5,10.4){\makebox(0,0)[b]{\scriptsize $\gamma \ot \id$}}
\put(19,10){\makebox(0,0){$\calP(k_1 + \cdots + k_n) \ot V^{\ot (k_1 + \cdots + k_n)}$}}
\put(0,0){\makebox(0,0){$\displaystyle\calP(n) \ot \bigotimes_{s=1}^n \left(
    \calP(k_s) \ot V^{\ot k_s} \rule{0em}{1.5em} \right)$}}
\put(20,0){\makebox(0,0){$\calP(n) \ot V^{\ot n}$}}
\put(6.7,0){\vector(1,0){10}}
\put(11.5,.4){\makebox(0,0)[b]{\scriptsize $ \id \ot (\bigotimes_{s=\rm 1}^n \alpha)$}}
\put(0,8.5){\vector(0,-1){6.5}}
\put(-.4,5){\makebox(0,0)[r]{\scriptsize  shuffle}}
\put(20,5){\makebox(0,0){$V$}}
\put(20,9){\vector(0,-1){3}}
\put(20.4,7.5){\makebox(0,0)[l]{\scriptsize $ \alpha$}}
\put(20,1){\vector(0,1){3}}
\put(20.4,3.5){\makebox(0,0)[l]{\scriptsize $ \alpha$}}
\end{picture}
\end{center}

\noindent
{\em Equivariance.\/} For each $n \geq 1$ and $\sigma \in \Sigma_n$,
the following diagram commutes.
\begin{center}
\setlength{\unitlength}{1em}
\begin{picture}(15,6)(0,.9)
\put(2,5){\makebox(0,0){$\calP(n) \ot V^{\ot n}$}}
\put(13,5){\makebox(0,0){$\calP(n) \ot V^{\ot n}$}}
\put(5,5){\vector(1,0){5}}
\put(7.5,5.4){\makebox(0,0)[b]{\scriptsize $\sigma \ot \sigma^{-1}$}}
\put(7.5,1.3){\makebox(0,0){$V$}}
\put(3,4){\vector(3,-2){3.5}}
\put(12,4){\vector(-3,-2){3.5}}
\put(4,2.5){\makebox(0,0)[tr]{\scriptsize $\alpha$}}
\put(11,2.5){\makebox(0,0)[tl]{\scriptsize $\alpha$}}
\end{picture}
\end{center}

\noindent
{\em Unitality.\/}
For each $n \geq 1$, the following diagram commutes.
\begin{center}

\setlength{\unitlength}{1em}
\begin{picture}(10,6)
\put(0,5){\makebox(0,0){$\bfk \ot V$}}
\put(-.4,2.5){\makebox(0,0)[r]{\scriptsize $\eta \ot \id$}}
\put(0,4){\vector(0,-1){3}}
\put(4.5,5.5){\makebox(0,0)[b]{\scriptsize $\cong$}}
\put(2,5){\vector(1,0){5}}
\put(4,2.9){\makebox(0,0)[b]{\scriptsize $\alpha$}}
\put(2,1){\vector(3,2){5}}
\put(8,5){\makebox(0,0){$V$}}
\put(-.6,0){\makebox(0,0){$\calP(1) \ot V$}}
\end{picture}
\end{center}
\end{proposition}

We leave as an exercise to formulate a version of
Proposition~\ref{pro:alg} that would use $\circ_i$-operations instead of
$\gamma$-operations.

\begin{example}
\label{Jitce_budu_volat_az_v_patek}
{\rm
In this example we verify, using Proposition~\ref{pro:alg},
that algebras over the operad $\Com =
\coll\Com1$ recalled in Example~\ref{wx} are ordinary commutative
associative algebras.  To simplify the exposition, let us agree that
$v$'s with various subscripts denote elements of $V$.  Since $\Com(n) =
\bfk$ for $n \geq 1$, the structure map~(\ref{algebra}) determines,
for each $n \geq 1$, a linear map $\mu_n : \otexp Vn \to V$ by
\[
\mu_n(\Rada v1n) := \alpha(1_n, \Rada v1n),
\]
where $1_n$ denotes in this example the unit $1_n \in \bfk
=\Com(n)$. The associativity of Proposition~\ref{pro:alg} says that
\begin{equation}
\label{qqq1}
\mu_n\left(\rr\mu_{k_1}(\Rada v1{k_1}),\ldots,
\mu_{k_n}(\Rada v{k_1+\cdots+k_{n-1}+1}{k_1+\cdots+k_n})\right) =
\mu_{k_1+\cdots+k_n}(\Rada v1{k_1+\cdots+k_n}),
\end{equation}
for each $n,\Rada k1n \geq 1$.  The
equivariance of Proposition~\ref{pro:alg} means that each $\mu_n$ is
fully symmetric
\begin{equation}
\label{qqq2}
\mu_n(\Rada v1n) = \mu_n(\Rada v{\sigma(1)}{\sigma(n)}),\ \sigma \in \Sigma_n,
\end{equation}
and the unitality implies that $\mu_1$ is the identity map,
\begin{equation}
\label{qqq3}
\mu_1(v) =v.
\end{equation}

The above structure can be identified with a commutative associative
multiplication on $V$. Indeed, the bilinear map $\cdot : = \mu_2: V \ot
V\to V$ is clearly associative:
\begin{equation}
\label{qqq4}
(v_1 \cdot v_2) \cdot v_3 = v_1 \cdot (v_2 \cdot v_3)
\end{equation}
and commutative:
\begin{equation}
\label{qqq5}
v_1 \cdot v_2 = v_2 \cdot v_1.
\end{equation}
On the other hand, $\mu_1(v) := v$ and
\[
\mu_n(\Rada v1n) :=  (\cdot\! \cdot\! \cdot (v_1 \cdot v_2)
\ \cdots \ v_{n-1})\cdot v_n\ \mbox { for $n \geq 2$}
\]
defines multilinear maps $\{\mu_n : \otexp Vn \to V\}$
satisfying~(\ref{qqq1})--(\ref{qqq3}).
It is equally easy to verify that algebras over the operad
$\Associative$ introduced in Example~\ref{9} are ordinary associative
algebras.
}\end{example}

Following Leinster~\cite{leinster:book}, one could say
that~(\ref{qqq1})--(\ref{qqq3}) is an {\em unbiased\/} definition
of associative commutative algebras, while~(\ref{qqq4})--(\ref{qqq5})
is a definition of the same object {\em biased\/} towards bilinear
operations. Operads therefore provide unbiased definitions of algebras.

\begin{example}
{\rm
Let us denote by $U\Com$ the endomorphism operad $\End_\bfk$ of
the ground ring $\bfk$. It is easy to verify that $U\Com$-algebras are
{\em unital\/} commutative associative algebras. We leave it to
the reader to describe the operad $U\Associative$ governing unital
associative operads.
}
\end{example}

Algebras over a \ns\ operad $\ucalP$ are defined as algebras, in the
sense of Definition~\ref{jitky}, over the symmetrization
$\Sigma[\ucalP]$ of $\ucalP$ .  Algebras over non-unital
operads discussed in Section~\ref{sec1bis}
are defined by appropriate obvious modifications of
Definition~\ref{jitky}.

\begin{example}
{\rm
Let $Y$ be a topological space with a base point $*$ and
${\mathbb S}^k$ the $k$-dimensional sphere, $k \geq 1$. The $k$-fold
loop space $\Omega^k Y$ is the space of all continuous maps
${\mathbb S}^k \to Y$ that send the south pole of ${\mathbb S}^k$ to
the base point of $Y$. Equivalently, $\Omega^k Y$ is the space of all
continuous maps $\lambda : (\DD^k,\SS^{k-1}) \to (Y,*)$ from the
standard closed $k$-dimensional disc $\DD^k$ to $Y$ that map the
boundary $\SS^{k-1}$ of $\DD^k$ to the base point of $Y$.  Let us show,
following Boardman and Vogt~\cite{boardman-vogt:73}, that $\Omega^k Y$ is a
natural topological algebra over the little $k$-discs operad $\calD_k =
\{\calD_k(n)\}_{n \geq 0}$ recalled in
Example~\ref{tuto_sobotu_zacinaji_zavody}.

The action $\alpha:\calD_k(n) \times (\Omega^k Y)^{\times n} \to \Omega^k Y$
is, for $n \geq 0$, defined as follows.  Given
$\lambda_i:(\DD^k,\SS^{k-1}) \to (Y,*)\in \Omega^k Y$, $1\leq i \leq
n$, and little $k$-discs $d =(d_1,\dots, d_n) \in \calD_k(n)$ as in
Example~\ref{tuto_sobotu_zacinaji_zavody}, then
\[
\alpha(d,\lambda_1,\dots,\lambda_n) :
(\DD^k,\SS^{k-1}) \to (Y,*) \in \Omega^k Y
\]
is the map defined to be $\lambda_i: \DD^k \to Y$
(suitably rescaled) on the image of $d_i$, and to be $*$ on the
complement of the images of the maps $d_i$, $1 \leq i \leq n$.

Therefore each $k$-fold loop space is a $\calD_k$-space.  The
following classical theorem is a certain form of the inverse
statement.

\begin{theorem}
(Boardman-Vogt~\protect\cite{boardman-vogt:73}, May~\protect\cite{may:1972})
A path-connected $\calD_k$-algebra $X$ has the weak homotopy type of a
$k$-fold loop space.
\end{theorem}

The connectedness assumption in the above theorem can be weakened by
assuming that the $\calD_k$-action makes the set $\pi_0(X)$ of path
components of $X$ a group.
}
\end{example}

\begin{example}
{\rm
The non-unital operad $\barcalM_0$ of stable pointed curves of genus~$0$
(also called the {\em configuration (non-unital) operad\/}) recalled on
page~\pageref{New_Orleans} is a non-unital operad in the category of smooth
complex projective varieties. It therefore makes sense, as explained in
Example~\ref{projdu_prohlidkou?}, to consider its
homology operad $H_*(\barcalM_0,\bfk) =
\{H_*(\barcalM_0(n),\bfk)\}_{n \geq 2}$.

An algebra over this non-unital operad is called a (tree level) {\em
cohomological conformal field theory\/} or a {\em hyper-commutative
algebra\/}~\cite{kontsevich-manin:CMP94}.  It consist of  a  family
$\left\{(\rada --): V^{\otimes n} \to V\right\}_{n \geq 2}$  of linear
operations
which are totally symmetric, that is
\[
(\rada{v_{\sigma(1)}}{v_{\sigma(n)}}) = (\Rada v1n),\
\]
for each permutation $\sigma \in \Sigma_n$. Moreover, we require the
following form of associativity:
\begin{equation}
\label{houstne_to}
\sum_{(S,T)} \((u,v,x_i;\ i~\in S),w,x_j;\ j \in T)
=
\sum_{(S,T)} (u,(v,w,x_i;\ i~\in S),x_j; j \in T),
\end{equation}
where $u,v,w,\Rada x1n \in V$ and $(S,T)$ runs over disjoint
decompositions $S \sqcup T = \{\rada 1n\}$.
For $n=0$,~(\ref{houstne_to}) means the (usual) associativity of the
bilinear operation
$(-,-)$, i.e.~$\((u,v),w)= (u,(v,w\))$. For $n=1$ we get
\[
\((u,v),w,x) + \((u,v,x),w)= (u,(v,w,x\)) + (u,(v,w),x).
\]
}
\end{example}

\begin{example}
{\rm
In this example, $\bfk$ is a field of characteristic~$0$.
The non-unital operad $\barcalM_0(\RR) =
\{\barcalM_0(\RR)(n)\}_{n \geq 2}$ of real points in the configuration
operad $\barcalM_0$ is called the {\em mosaic
non-unital operad\/}~\cite{devadoss:tes}. Algebras over the homology
$H_*(\barcalM_0(\RR),\bfk) = \{H_*(\barcalM_0(\RR)(n),\bfk)\}_{n \geq 2}$
of this operad were recently identified~\cite{EHKR} with {\em
$2$-Gerstenhaber algebras\/}, which are structures $(V,\mu,\tau)$
consisting of a commutative associative product $\mu : V\ot V \to V$ and
an anti-symmetric degree $+1$ ternary operation $\tau : V \ot V
\ot V \to V$  which satisfies the generalized Jacobi identity
\[
\sum_\sigma
{\it sgn\/}(\sigma) \cdot \tau(\tau(x_{\sigma(1)},
x_{\sigma(2)},x_{\sigma(3)}),
x_{\sigma(4)},x_{\sigma(5)}) =0,
\]
where the summation runs over all $(3,2)$-unshuffles $\sigma(1) <
\sigma(2)<\sigma(3)$, $\sigma(4) < \sigma(5)$. Moreover, the ternary
operation $\tau$
is tied to the multiplication $\mu$ by the distributive law
\[
\tau(\mu(s,t),u,v) =   \mu(\tau(s,u,v),t) +
(-1)^{(1 +|u| + |v|)|s|} \cdot   \mu(s,\tau(t,u,v)),\ s,t,u,v \in V,
\]
saying that the assignment $s \mapsto \tau(s,u,v)$ is a degree $(1+|u|
+ |v|)$-derivation of the associative commutative algebra $(V,\mu)$, for each
$u,v \in V$.
}
\end{example}

\section{Free operads and trees}
\label{sec3}

The purpose of this section is three-fold. First, we want to study
free operads because each operad is a quotient of a free one. The
second reason why we are interested in free operads is that their
construction involves trees. Indeed, it turns out that rooted trees
provide `pasting schemes' for operads and that, replacing trees by
other types of graphs, one can introduce several important
generalizations of operads, such as cyclic operads, modular operads,
and \PROP{s}. The last reason is that the free operad functor defines
a monad which provides an unbiased definition of operads as algebras
over this monad. Everything in this section is written for
$\bfk$-linear operads, but the constructions can be generalized
into an arbitrary symmetric monoidal category with countable
coproducts $(\calM,\odot,{\bf 1})$ whose monoidal product $\odot$ is
distributive over coproducts,
see~\cite[Section~II.1.9]{markl-shnider-stasheff:book}.

Recall that a {\em $\Sigma$-module\/} is a collection $E =
\coll E0$ in which each $E(n)$ is a right $\bfk[\Sigma_n]$-module.
There is an obvious forgetful functor\/ $\forget :\catOper \to
\sigmamod$ from the category  $\catOper$ of $\bfk$-linear operads
to the category $\sigmamod$ of $\Sigma$-modules.

\begin{definition}
\label{Ja_tu_jitku_snad_opravdu_miluji}
The {\em free operad functor\/} is a left
adjoint~\cite[\S~II.7]{hilton-stammbach} $\Gamma : \sigmamod
\to \catOper$ to the forgetful functor\/ $\forget : \catOper
\to \sigmamod$. This means that there exists a functorial
isomorphism
\[
\Mor_{\catOper}(\Gamma(E),\calP) \cong
\Mor_{\ssigmamod}(E,\forget(\calP))
\]
for an arbitrary $\Sigma$-module $E$ and operad $\calP$. The operad
$\Gamma(E)$ is the {\em free operad\/} generated by the $\Sigma$-module
$E$. Similarly,
the {\em free non-unital operad functor\/} is a left adjoint $\Psi :
\sigmamod \to \psi\catOper$ of the obvious forgetful functor
$\forget_\psi :\psi\catOper \to \sigmamod$, that is
\[
\Mor_{\psi\catOper}(\Psi(E),\calS) \cong
\Mor_{\ssigmamod}(E,\forget_\psi(\calS)),
\]
where $E$ is a $\Sigma$-module and $\calS$ a non-unital operad. The
non-unital operad $\Psi(E)$ is the {\em free non-unital operad\/} generated by
the  $\Sigma$-module $E$.
\end{definition}

Let $\widetilde {\rule {1em}{0em}} \hskip -.5em : \psi\catOper
\to \catOper$ be the functor of `adjoining the unit' considered
in the proof of
Proposition~\ref{toto_pisu_na_zavodech_v_Moravske_Trebove_2005}
on page~\pageref{uz_je_po_zavodech}.  Functorial
isomorphism~(\ref{probihaji_zavody_a_jsem_v_Praze_protoze_je_strasne_pocasi!!})
implies that one may take
\begin{equation}
\label{dnes_jsem_boural_auto}
\Gamma := \widetilde {\Psi},
\end{equation}
which means that the free operad $\Gamma(E)$ can be obtained from the
free non-unital operad $\Psi(E)$ by formally adjoining the unit.

Let us indicate how to construct the free non-unital operad $\Psi(E)$,
a precise description will be given later in this section. The free
non-unital operad $\Psi(E)$ must be built up from all formal
$\circ_i$-compositions of elements of $E$ modulo the axioms listed in
Definition~\ref{b3}. For instance, given $f \in E(2)$, $g\in E(3)$,
$h \in E(2)$ and $l \in E(0)$, the component $\Psi(E)(5)$ must
contain the following five compositions
\begin{equation}
\label{Byl_jsem_s_Jituskou!}
\def\arraystretch{1.4}
\begin{array}{c}
(f \c 1(g\c 2 l)) \c3 h,\
(f\c2 h) \c1 (g \c2 l),\
((f \c2 h)\c1 g)\c 2 l,
\\
((f \c1 g)\c2 l) \c3 h\ \mbox { and }\
((f \c1g)\c4 h)\c2 l.
\end{array}
\end{equation}

\begin{figure}[t]
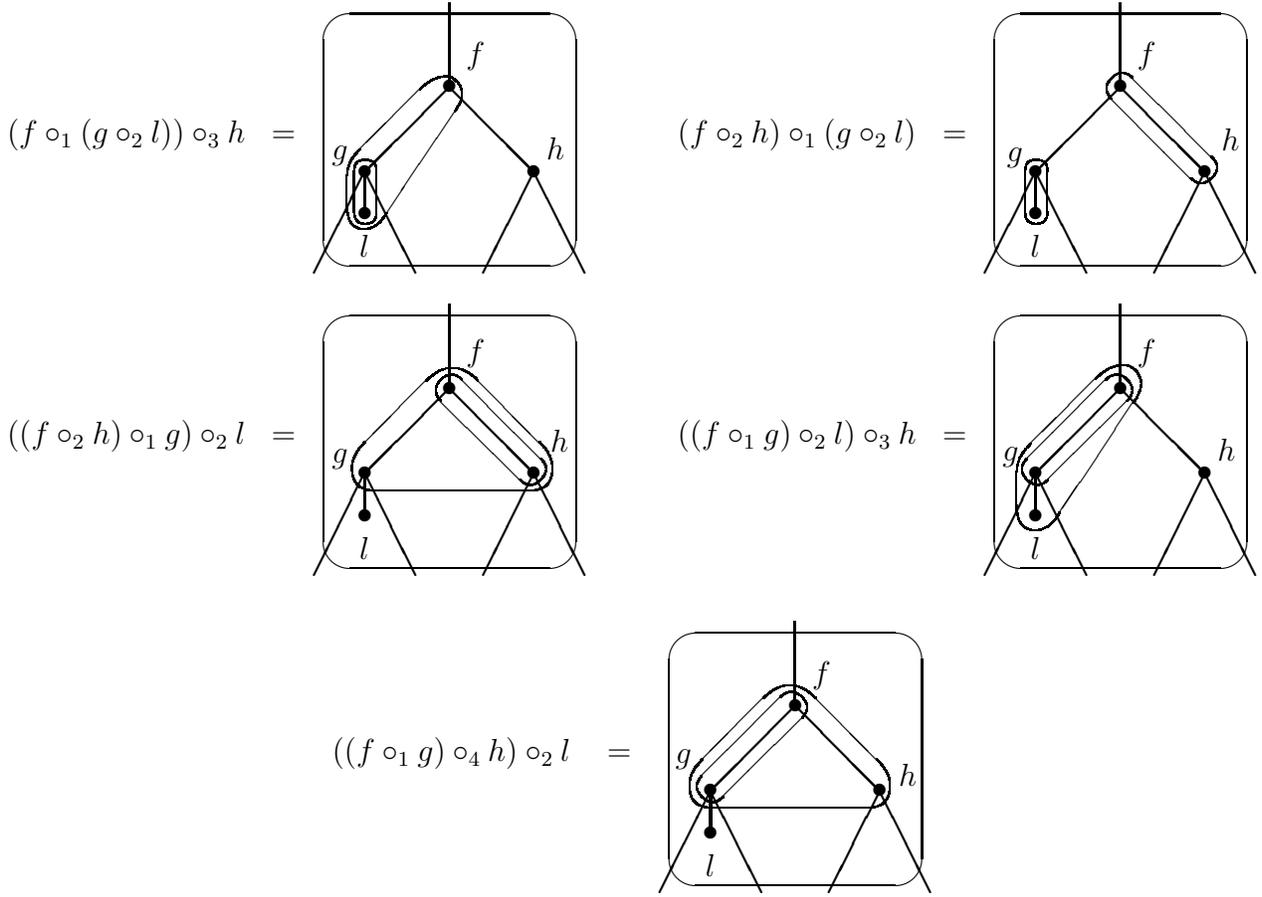

\[
\begin{array}{rclrcl}
(f \c 1(g\c 2 l)) \c3 h &  = & \kozulkanovaI \hskip 1cm &
(f\c2 h) \c1 (g \c2 l)& =& \kozulkanovaII
\\
((f \c2 h)\c1 g)\c 2 l& = &\kozulkanovaIII\rule{0pt}{110pt}&
((f \c1 g)\c2 l) \c3 h &= & \kozulkanovaIV
\end{array}
\]
\[
((f \c1g)\c4 h)\c2 l \hskip .9em  = \hskip .9em\kozulkanovaV\rule{0pt}{110pt}
\]
\caption{\label{porad_na_ni_myslim}Flow diagrams in non-unital
  operads.
\protect\rule{0pt}{60pt}}
\end{figure}

The elements in~(\ref{Byl_jsem_s_Jituskou!}) can be depicted by the
`flow diagrams' of Figure~\ref{porad_na_ni_myslim} on
page~\pageref{porad_na_ni_myslim}.
Nodes of these diagrams are decorated by elements $f,g,h$ and $l$ of
$E$ in such a way that an element of $E(n)$ decorates a node with $n$
input lines, $n \geq 0$.  Thin `amoebas' indicate the nesting which
specifies the order in
which the $\circ_i$-operations are performed.  The associativity of
Definition~\ref{b3} however says that the result of the composition
does not depend on the order, therefore the amoebas can be erased and
the common value of the compositions represented by
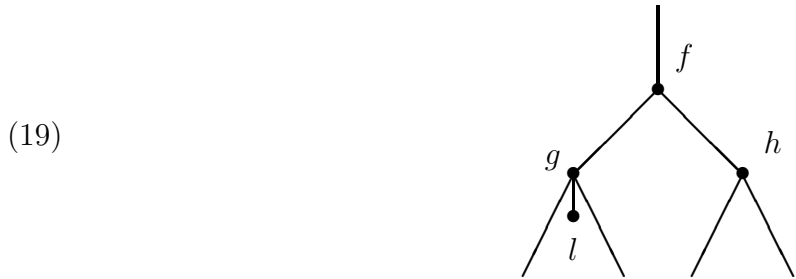
\begin{equation}
\label{kozulka}
\kozulkanova
\end{equation}
\rule{0pt}{39pt}

Let us look more closely how diagram~(\ref{kozulka}) determines an
element of the (still hypothetical) free non-unital operad $\Psi(E)$. The
crucial\label{prohlidku_mam_za_sebou!}
fact is that the underlying graph of~(\ref{kozulka}) is a
planar rooted tree.  Recall that a {\em tree\/} is a finite connected
simply connected graph without loops and multiple edges. For a tree
$T$ we denote, as usual, by $\Vert(T)$ the set of vertices and
$\Edg(T)$ the set of edges of $T$. The number of
edges adjacent to a vertex $v \in \Vert(T)$  is called the {\em valence\/} of
$v$ and denoted $\val(v)$.
We assume that one is given a subset
\[
\ext(T) \subset \{ v \in \Vert(T);\ \val(v) = 1\}
\]
of {\em external\/} vertices, the remaining vertices are  {\em
internal\/}.
Let us denote
\[
\vert(T) : = \Vert(T) \setminus \ext(T)
\]
the set of all internal vertices.  Henceforth, we will assume that our trees have
at least one internal vertex. This excludes at this stage the {\em exceptional
tree\/} consisting of two external vertices connected by an edge.

Edges adjacent to external
vertices are the {\em legs\/} of $T$.  A tree is {\em rooted\/} if one
of its legs, called the {\em root\/}, is marked and all other edges
are oriented, pointing to the root.  The legs different from the root
are the {\em leaves\/} of $T$.
For example, the tree in~(\ref{kozulka}) has 4~internal vertices
decorated $f$, $g$, $h$ and $l$, and 4~leaves.
Finally, the {\em planarity\/} means
that an embeddings of $T$ into the plane is specified. In our pictures, the
root will always be placed on the top. By a vertex we will always mean
an internal one.

The planarity and a choice of the root of the underlying tree
of~(\ref{kozulka}) specifies a total order of the set $\In(v)$ of
input edges of each vertex $v \in \vert(T)$ as well as a total order
of the set $\Leaf(T)$ of the leaves of $T$, by numbering
from the left to the right:
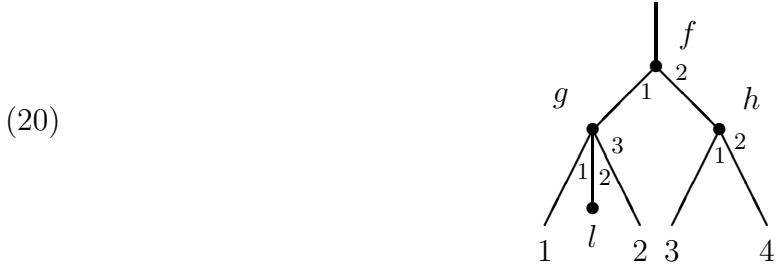
\begin{equation}
\label{g}
\kozulkalabeled
\end{equation}
\rule{0pt}{35pt}

\noindent
This tells us that $l$ should be inserted into the second input of
$g$, $g$ into the first input of $f$
and $h$ into the second input of $f$.
Using `abstract variables' $v_1,v_2,v_3$ and $v_4$, the element
represented by~(\ref{g}) can also be written as the `composition'
$f(g(v_1,l,v_2),h(v_3,v_4))$.

Now we need to take into account also the symmetric group action. If
$\tau$ is the generator of $\Sigma_2$, then the obvious equality
\[
f(g(v_1,l,v_2),h(v_3,v_4)) = f\tau (h(v_3,v_4),g(v_1,l,v_2))
\]
of `abstract compositions' coming from the equivariance of
Definition~\ref{b3} translates into the following equality of flow
diagrams:
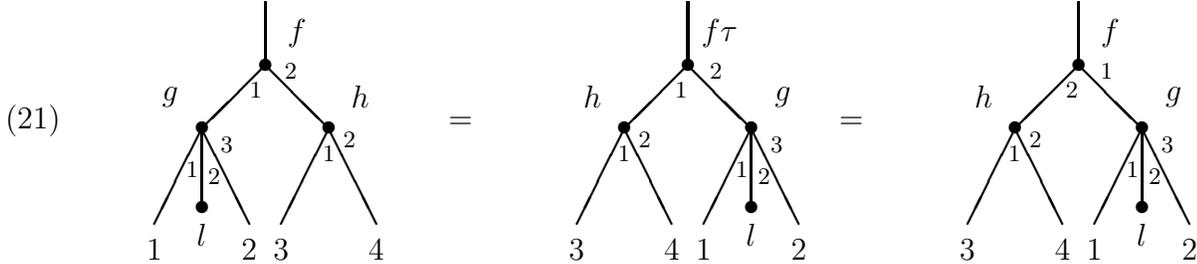
\begin{equation}
\label{kralicek}
\kozulkalabeled \hskip 30pt =  \hskip 30pt\kozulkalabeledtau
 \hskip 30pt =  \hskip 30pt\kozulkalabeledtaubis
\end{equation}
\rule{0pt}{45pt}

Relation~(\ref{kralicek}) shows that the equivariance of
Definition~\ref{b3} violates the linear orders induced by the planar
embedding of $T$. This leads us to the conclusion
that the flow diagrams describing
elements of free non-unital operads are (abstract, non-planar) rooted,
leaf-labeled decorated trees.

Let us describe, after these motivations, a precise construction of
$\Psi(E)$. The first subtlety one needs to understand is how to decorate
vertices of non-planar trees.  To this end, we need to explain how each
$\Sigma$-module $E = \coll E0$ naturally extends into a functor
(denoted again~$E$) from the category $\Set$ of finite sets and
their bijections to the category of $\bfk$-modules.
If $X$ and $Y$ are finite sets, denote by
\begin{equation}
\label{maminka_se_zase_nacpala}
\Bij(Y,X) := \{\vartheta : X \stackrel{\cong}{\longrightarrow} Y\}
\end{equation}
the set of all isomorphisms between $X$ and $Y$ (notice the
unexpected direction of the arrow!). It is clear that $\Bij(Y,X)$
is a natural left $\Aut_Y$- right $\Aut_X$-bimodule, where
$\Aut_X := \Bij(X,X)$ and $\Aut_Y := \Bij(Y,Y)$ are the sets of
automorphisms with group structure given by composition. For a
finite set $S \in \Set$ of cardinality $n$ and a $\Sigma$-module
$E = \coll E0$ define $E(S)$ to be
\begin{equation}
\label{m-m}
E(S) := E(n) \times_{\Sigma_n} \Bij([n],S)
\end{equation}
where, as usual, $[n] := \{\rada 1n\}$ and, of course, $\Sigma_n =
\Aut_{\en}$.

Let us recall that a {\em (leaf-) labeled rooted $n$-tree\/} is a rooted tree
$T$ together with a specified bijection $\ell : \Leaf(T)
\stackrel{\sim}{\rightarrow} [n]$. Let $\TTT_n$ be the category of
labeled rooted $n$-trees and their bijections. For $T \in \TTT_n$
define
\begin{equation}
\label{ale_jeste_mne_cekaji_odbery}
E(T) := \bigotimes_{v \in \vert(T)} E(\In(v))
\end{equation}
where $\In(v)$ is, as before, the set of all input edges of a
vertex $v \in \vert(T)$.
It is easy to verify that $E \mapsto E(T)$ defines a functor
from the category $\TTT_n$ to the category of $\bfk$-modules.

Recall that the colimit of a covariant functor
$F: \calD \to \Modk$ is the quotient
\[
\colim{x \in \calD} F(x) = \bigoplus_{x \in \calD} F(x)/\sim,
\]
where $\sim$ is the equivalence generated by
\[
F(y) \ni a \sim F(f)(a) \in F(z),
\]
for each $a \in F(y)$, $y,z \in \calD$ and $f \in {\it
Mor\/}_\calD(y,z)$.  Define finally
\begin{equation}
\label{msss}
\Psi(E)(n) := \colim{T \in \TTT_n} E(T),\ n \geq 0.
\end{equation}
The following theorem was proved in~\cite[II.1.9]{markl-shnider-stasheff:book}.

\begin{theorem}
There exists a natural non-unital operad structure on the $\Sigma$-module
\[
\Psi(E) = \coll {\Psi(E)}0,
\]
with the $\circ_i$-operations given by the grafting of trees and the
symmetric group re-labeling the leaves, such that $\Psi(E)$ is the
free non-unital operad generated by the $\Sigma$-module $E$.
\end{theorem}

One could simplify~(\ref{msss}) by introducing $\tree(n)$ as the
set of {\em isomorphism\/} classes of $n$-trees from $\TTT_n$ and
defining $\Psi(E)$ by the formula
\begin{equation}
\label{draw}
\Psi(E)(n) = \bigoplus_{[T]\in \tree(n)}E(T),\ n \geq 0,
\end{equation}
which does not involve the colimit. The drawback of~(\ref{draw}) is
that it assumes a choice of a representative $[T]$ of each isomorphism
class in $\tree(n)$, while~(\ref{msss}) is functorial and admits
simple generalizations to other types of operads and {\small
PROP}s. See~\cite[Section~II.1.9]{markl-shnider-stasheff:book} for
other representations of the free non-unital operad functor.

Having constructed the free non-unital operad $\Psi(E)$, we may
use~(\ref{dnes_jsem_boural_auto}) to define the free operad
$\Gamma(E)$. This is obviously equivalent to
enlarging, in~(\ref{msss}) for $n=1$, the category $\TTT_n$ by
the {\em exceptional roted tree\/} \hskip .5em
\raisebox{-.2em}{\rule{.8pt}{1.1em}}  \hskip .5em
with one leg and no internal
vertex. If we denote this enlarged category of trees and their
isomorphisms (which however differs from $\TTT_n$ only at $n=1$) by
$\UTTT_n$, we may represent the free operad as
\begin{equation}
\label{freeop}
\Gamma(E)(n) := \colim{T \in \UTTT_n} E(T),\ n \geq 0.
\end{equation}

If $E$ is a $\Sigma$-module such that $E(0) = E(1) = 0$,
then~(\ref{draw}) reduces to a summation over {\em reduced
  trees\/}, that is trees whose all vertices
have at least two input edges. By simple combinatorics, the number
of isomorphism classes of reduced trees in $\TTT_n$ is finite for each
$n \geq 0$.
This implies the following proposition that says that operads are
relatively small objects.

\begin{proposition}
\label{to_by_mne_zajimalo_jestli_jsem_to_s_ni_pokazil}
Let $E = \coll E0$ be a $\Sigma$-module such that
\[
E(0) = E(1) = 0
\]
and that $E(n)$ are
finite-dimensional for $n \geq 2$. Then the spaces $\Psi(E)(n)$
and  $\Gamma(E)(n)$ are finite-dimensional for each $n \geq 0$.
\end{proposition}

We close this section by showing how the free operad functor
can be used to define operads. It follows from general
principles that any operad $\calP$ is a quotient
$\calP =\Gamma(E)/(R)$, where $E$ and $R$ are $\Sigma$-modules and $(R)$ is
the operadic ideal (see Definition~\ref{pristi_tyden_prohlidka})
generated by $R$ in $\Gamma(E)$.

\begin{example}
\label{000}
{\rm
The commutative associative operad $\Com$ recalled in Example~\ref{wx}
is generated by the $\Sigma$-module
\[
E_{\sCom}(n) :=
\cases{\bfk\cdot\mu}{if $n = 2$ and}0{if $n \not= 2$.}
\]
where $\bfk\cdot\mu$ is the trivial representation of $\Sigma_2$.  The
ideal of relations is generated by
\[
R_{\sCom} : = \Span_\bfk\{\mu (\mu \ot \id) - \mu(\id \ot \mu)\}
\subset \Gamma(E_{{\sCom}})(3),
\]
where $\mu (\mu \ot \id) - \mu(\id \ot \mu)$ is the obvious shorthand
for $\gamma(\mu,\mu,e) - \gamma(\mu,e,\mu)$, with $e$ the unit of
$\Gamma(E_{\sCom})$.

Similarly, the operad $\Ass$ for associative algebras reviewed in
Example~\ref{9} is generated by
the $\Sigma$-module $E_{\sAss}$ such that
\[
E_{\sAss}(n) := \cases{\bfk[\Sigma_2]}{if $n = 2$ and}0{if $n \not= 2$.}
\]
The ideal of relations is generated by the
$\bfk[\Sigma_3]$-closure $R_{\sAss}$ of the associativity
\begin{equation}
\label{associ}
\alpha (\alpha \ot \id)
- \alpha(\id \ot \alpha) \in \Gamma(E_{{\sAss}})(3),
\end{equation}
where $\alpha$ is a generator of the regular representation
$E_{\sAss}(2) = \bfk[\Sigma_2]$.
}
\end{example}

\begin{example}
\label{Lie}{\rm
The operad $\Lie$ governing Lie algebras is the quotient
$\Lie := \Gamma(E_{\sLie})/(R_{\sLie})$, where $E_{\sLie}$ is the
$\Sigma$-module
\[
E_{\sLie}(n) :=\left\{\begin{array}{ll}\bfk\cdot\beta&\mbox{  if }n=2
\mbox { and}\\
0&\mbox{ if } n\neq 2, \end{array}\right.
\]
with $\bfk\cdot\beta$ is the signum representation of $\Sigma_2$.
The ideal of relations $(R_{\sLie})$ is generated by the Jacobi
identity:
\begin{equation}
\label{JJJ}
\beta(\beta \ot \id) + \beta(\beta \ot \id)c + \beta(\beta \ot \id)c^2 = 0,
\end{equation}
in which $c \in \Sigma_3$ is the cyclic permutation $(1,2,3) \mapsto (2,3,1)$.
}
\end{example}

\begin{example}
\label{Jitka_ted_moc_nepise}
{\rm\
We show how to describe the presentations
of the operads $\Ass$ and $\Lie$ given
in Examples~\ref{000} and~\ref{Lie} in a simple graphical language.
The generator $\alpha$ of $E_{\sAss}$ is
an operation with two inputs and one output, so we depict it as $\gen 12$. The
associativity~(\ref{associ}) then reads as
\[
\ZbbZb = \bZbbZ,
\]
therefore
$\Ass = \Gamma(\gen 12)/(\ZbbZb = \bZbbZ)$.
Also the operad for $\Lie$ algebras is generated by one bilinear
operation $\gen 12$, but this time the operation is anti-symmetric
\[
\glab 12 = - \glab 21.
\]
The Jacobi identity~(\ref{JJJ}) reads
\[
\Jac 123 + \Jac 231 + \Jac 312 = 0.
\]
}
\end{example}

The kind of description used in the above examples is
`tautological' in the sense that it just says that the operad $\calP$
governing a certain type of algebras is generated by operations of
these algebras, with an appropriate symmetry, modulo the axioms
satisfied by these operations. It does not say directly anything about
the properties of the individual spaces $\calP(n)$, $n \geq 0$. Describing
these individual components may be a very nontrivial task, see
for example the formula for the $\Sigma_n$-modules $\Lie(n)$ given
in~\cite[page~50]{markl-shnider-stasheff:book}.
Operads in Examples~\ref{000} and~\ref{Lie} are
quadratic in the sense of the following:

\begin{definition}
\label{je_tady_neporadek}
An operad $\calP$ is {\em quadratic\/} if it has a presentation $\calP =
\Gamma(E)/(R)$, where $E=\calP(2)$ and $R \subset \Gamma(E)(3)$.
\end{definition}

Quadratic operads form a very important class of operads.  Each
quadratic operad $\calP$ has a {\em quadratic dual\/}
$\calP^!$~\cite{ginzburg-kapranov:DMJ94},~\cite[Definition~II.3.37]%
{markl-shnider-stasheff:book}
which is a quadratic operad defined, roughly
speaking, by dualizing the generators of $\calP$ and replacing the
relations of $\calP$ by their annihilator in the dual space.  For
example, $\Ass^! = \Ass$, $\Com^! = \Lie$ and $\Lie^! = \Com$.  A
quadratic operad $\calP$ is {\em Koszul\/} if it has the homotopy type
of the bar construction of its quadratic
dual~\cite{ginzburg-kapranov:DMJ94},~\cite[Definition~II.3.40]%
{markl-shnider-stasheff:book}.  For
quadratic Koszul operads, there is a deep understanding of the
derived category of the corresponding algebras. Operads
$\Ass$, $\Com$ and $\Lie$ above, as well as most quadratic operads
one encounters in everyday life, are Koszul.

\section{Unbiased definitions}
\label{sec4}

In this section, we review the definition of a triple (monad) and
give, in Theorem~\ref{Jana_neni_doma}, a~description of unital
and non-unital operads in terms of algebras over a triple. The
relevant triples come from the endofunctors $\Psi$ and $\Gamma$
recalled in Section~\ref{sec3}.  Let ${\it End}({\calC})$ be the
strict symmetric monoidal category of endofunctors on a category
$\calC$ where multiplication is the composition of functors.

\begin{definition}
\index{triple}\label{deftriple}\index{monad}
A {\em triple} (also called a {\em monad\/}) $T$ on a category $\calC$ is an
associative and unital monoid $\left( T,\mu ,\upsilon \right) $ in
${\it End}(\calC)$. The multiplication $\mu : TT \to T$
and unit morphism $\upsilon:  {\it id} \to T$
satisfy the axioms given by commutativity of the diagrams in
Figure~\ref{triple}.
\end{definition}

Triples arise naturally from pairs of adjoint functors. Given an
adjoint pair~\cite[II.7]{hilton-stammbach}
\[
\bijections{{\tt A}}{{\tt B}}{G}{F},
\]
with associated functorial isomorphism
\[
\Mor_{\tt A}(F(X),Y)
\cong \Mor_{\tt B}(X,G(Y)),\ X \in {\tt B},\ Y \in {\tt A},
\]
there is a triple in ${\tt B}$ defined by $T: =GF$.  The unit of the
adjunction $\id\rightarrow GF$ defines the unit $\upsilon$ of the
triple and the counit of the adjunction $FG\rightarrow \id$ induces a
natural transformation $GFGF\rightarrow GF$ which defines the
multiplication $\mu.$ In fact, it is a theorem of Eilenberg and
Moore~\cite{e-m} that all triples arise in this way from adjoint
pairs.  This is exactly the situation with the free operad and free
non-unital operad functors that were described in Section~\ref{sec3}.  We
will show how operads and non-unital operads can actually be defined using
the concept of an algebra over a triple:

\begin{definition}
\label{snad_projdu}
A  {\em $T$-algebra\/} or {\em algebra over the triple $T$\/}
is an object $A$ of\/ $\calC$  together with a structure morphism
$\alpha :T(A)\rightarrow A$
satisfying
\[
\alpha (T(\alpha ))=\alpha(\mu_A)  \mbox{ and }
\alpha \upsilon_A=\id_{A},
\]
see Figure~\ref{Talgebra}.
\end{definition}
\begin{figure}
\begin{center}
{\setlength{\unitlength}{.66cm}
\begin{picture}(15,4)(2,0)
%
\put(2,3){\makebox(0,0){$TTT$}}
\put(6,3){\makebox(0,0){$TT$}}
\put(2,0){\makebox(0,0){$TT$}}
\put(6,0){\makebox(0,0){$T$}}
\put(1.7,1.5){\makebox(0,0)[r]{$\mu T$}}
\put(6.3,1.5){\makebox(0,0)[l]{$\mu$}}
\put(4.1,0.5){\makebox(0,0)[b]{$\mu$}}
\put(4.1,3.5){\makebox(0,0)[b]{$T\mu$}}
\put(3,0){\vector(1,0){2}}
\put(3.2,3){\vector(1,0){1.8}}
\put(2,2.5){\vector(0,-1){2}}
\put(6,2.5){\vector(0,-1){2}}
%
\put(-.5,0){
\put(9.6,3){\makebox(0,0){$T$}}
\put(12.3,3){\makebox(0,0){$TT$}}
\put(10,3){\vector(1,0){1.8}}
\put(11,3.5){\makebox(0,0){$T\upsilon$}}
\put(12,2.5){\vector(-1,-3){0.7}}
\put(10,2.5){\vector(1,-3){0.7}}
\put(9.5,1.5){\makebox(0,0){$\id$}}
\put(12,1.5){\makebox(0,0){$\mu$}}
\put(11,0){\makebox(0,0){$T$}}
}
%
\put(14.6,3){\makebox(0,0){$T$}}
\put(17.3,3){\makebox(0,0){$TT$}}
\put(15,3){\vector(1,0){1.8}}
\put(15.8,3.5){\makebox(0,0){$\upsilon T$}}
\put(17,2.5){\vector(-1,-3){0.7}}
\put(15,2.5){\vector(1,-3){0.7}}
\put(14.5,1.5){\makebox(0,0){$\id$}}
\put(17,1.5){\makebox(0,0){$\mu$}}
\put(16,0){\makebox(0,0){$T$}}
\end{picture}}
\end{center}
\caption{\label{triple}
Associativity and unit axioms for a triple.}
\end{figure}

\begin{figure}
\begin{center}
{\setlength{\unitlength}{.66cm}
\begin{picture}(10,4)(2,0)
%
\put(2,3){\makebox(0,0){$T(T(A))$}}
\put(6,3){\makebox(0,0){$T(A)$}}
\put(2,0){\makebox(0,0){$T(A)$}}
\put(6,0){\makebox(0,0){$A$}}
\put(1.7,1.5){\makebox(0,0)[r]{$\mu$}}
\put(6.3,1.5){\makebox(0,0)[l]{$\alpha$}}
\put(4,0.5){\makebox(0,0)[b]{$\alpha$}}
\put(4.2,3.5){\makebox(0,0)[b]{$T(\alpha)$}}
\put(3,0){\vector(1,0){2}}
\put(3.2,3){\vector(1,0){1.8}}
\put(2,2.5){\vector(0,-1){2}}
\put(6,2.5){\vector(0,-1){2}}
%
\put(9.6,3){\makebox(0,0){$A$}}
\put(12.5,3){\makebox(0,0){$T(A)$}}
\put(10,3){\vector(1,0){1.6}}
\put(10.9,3.5){\makebox(0,0){$\upsilon_A$}}
\put(12,2.5){\vector(-1,-3){0.7}}
\put(10,2.5){\vector(1,-3){0.7}}
\put(9.5,1.5){\makebox(0,0){$\id$}}
\put(12.2,1.5){\makebox(0,0){$\alpha$}}
\put(11,0){\makebox(0,0){$A$}}
\end{picture}}
\end{center}
\caption{\label{Talgebra}
$T$-algebra structure.}
\end{figure}

The category of $T$-algebras in $\calC$ will be denoted ${\tt
Alg}_T(\calC).$ Since the free non-unital operad functor $\Psi$ and the
free operad functor $\Gamma$ described in Section~\ref{sec3} are left
adjoints to \ $\forget_\psi :\psi\catOper \to \sigmamod$ and \
$\forget :\catOper \to \sigmamod$, respectively, the functors \
$\forget_\psi\Psi$ (denoted simply $\Psi$) and \ $\forget\Gamma$
(denoted~$\Gamma$) define triples on $\sigmamod$.

\begin{theorem}
\label{Jana_neni_doma}
A $\Sigma$-module $\calS$ is a $\Psi$-algebra if and only if it is a
non-unital operad and it is a $\Gamma$-algebra if and only if it is an
operad.  In shorthand:
\[
{\tt Alg}_{\Psi}(\sigmamod) \cong \psi\catOper\
\mbox { and }\
{\tt Alg}_{\Gamma}(\sigmamod) \cong \catOper.
\]
\end{theorem}

\noindent
{\bf Proof.}
We outline first the proof of the implication in the direction from algebra to
non-unital operad.
Let $\calS$ be a $\Psi$-algebra. The restriction of the structure morphism
$\alpha:\Psi(\calS)\longrightarrow \calS$
to the components of $\Psi(\calS)$ supported on
trees with one internal edge defines the non-unital operad composition maps
$\circ_i$, as indicated by:
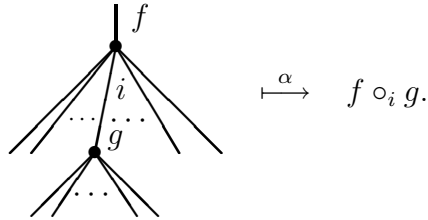
\begin{center}
{
\unitlength=.8pt
\begin{picture}(170.00,90.00)(0.00,10.00)
\thicklines
\put(118.00,60.00){\makebox(0.00,0.00)[l]%
{$\stackrel{\alpha}{\longmapsto}\hskip10pt f \circ_i g$.}}
\put(39.50,10.00){\makebox(0.00,0.00){$\cdots$}}
\put(57.00,45.00){\makebox(0.00,0.00){$\cdots$}}
\put(35.00,45.50){\makebox(0.00,0.00){\scriptsize $\cdots$}}
\put(46.00,31.00){\makebox(0.00,0.00)[lb]{$g$}}
\put(57.00,87.00){\makebox(0.00,0.00)[lb]{$f$}}
\put(50.00,60.00){\makebox(0.00,0.00)[l]{$i$}}
\put(40.00,30.00){\makebox(0.00,0.00){$\bullet$}}
\put(50.00,80.00){\makebox(0.00,0.00){$\bullet$}}
\put(40.00,30.00){\line(2,-3){20.00}}
\put(50.00,80.00){\line(3,-5){30.00}}
\put(40.00,30.00){\line(1,-1){30.00}}
\put(40.00,30.00){\line(-2,-3){20.00}}
\put(40.00,30.00){\line(-1,-1){30.00}}
\put(50.00,80.00){\line(1,-1){50.00}}
\put(50.00,80.00){\line(-4,-5){40.00}}
\put(50.00,80.00){\line(-1,-5){10.00}}
\put(50.00,80.00){\line(-1,-1){50.00}}
\put(50.00,100.00){\line(0,-1){20.00}}
\end{picture}}
\end{center}

In the opposite direction, for a non-unital operad $\calS$, the
$\Psi$-algebra structure $\alpha : \Psi(\calS) \to \calS$ is
the contraction along the edges of underlying trees, using the
$\circ_i$-operations. The proof that $\Gamma$-algebras are operads is
similar.%
\qed

Let us change our perspective and consider formula~(\ref{msss}) as
defining an endofunctor $\Psi : \sigmamod \to \sigmamod$,
ignoring that we already know that it represents free
non-unital operads. We are going to construct maps
\[
\mu  : \Psi \Psi \to \Psi
\mbox { and }
\upsilon : \id \to \Psi
\]
making $\Psi$ a triple on the category $\sigmamod$.
Let us start with the triple multiplication $\mu$.
It follows from~(\ref{msss}) that, for each $\Sigma$-module $E$,
\[
\Psi\Psi(E)(n) := \colim{T \in \TTT_n} \Psi(E)(T),\ n \geq 0.
\]
The elements in the right hand side are
represented by rooted trees $T$ with vertices decorated by
elements of $\Psi(E)$, while elements of $\Psi(E)$ are represented by
rooted trees with vertices decorated by $E$. We may therefore imagine
elements\label{vecer_zavolam_Jitce}
of $\Psi\Psi(E)$ as `bracketed' rooted trees, in the sense indicated in
Figure~\ref{prohlidka}. The triple multiplication $\mu_E : \Psi\Psi(E)
\to \Psi(E)$ then simply erases the braces.
\begin{figure}[t]
\begin{center}
{
\unitlength=1.000000pt
\begin{picture}(270.00,100.00)(0.00,0.00)
\put(130.00,50.00)%
     {\makebox(0.00,0.00){$\stackrel{\mu}\longmapsto$}}
\thicklines
\put(190.00,20.00){\makebox(0.00,0.00){$\bullet$}}
\put(220.00,30.00){\makebox(0.00,0.00){$\bullet$}}
\put(220.00,50.00){\makebox(0.00,0.00){$\bullet$}}
\put(50.00,40.00){\makebox(0.00,0.00){$\bullet$}}
\put(40.00,40.00){\makebox(0.00,0.00){$\bullet$}}
\put(50.00,50.00){\makebox(0.00,0.00){$\bullet$}}
\put(190.00,20.00){\line(1,-1){20.00}}
\put(220.00,50.00){\line(0,-1){20.00}}
\put(220.00,50.00){\line(1,-1){50.00}}
\put(220.00,50.00){\line(-1,-1){50.00}}
\put(220.00,100.00){\line(0,-1){50.00}}
\put(70.00,30.00){\line(1,-1){32.50}}
\put(30.00,30.00){\line(-1,-1){32.50}}
\put(50.00,30.00){\line(0,-1){32.50}}
\put(40.00,40.00){\line(1,-1){10.00}}
\put(50.00,50.00){\line(0,-1){10.00}}
\put(50.00,50.00){\line(1,-1){20.00}}
\put(50.00,50.00){\line(-1,-1){20.00}}
\put(50.00,100.00){\line(0,-1){50.00}}
\thinlines
\put(50.00,50.00){\circle{60.00}}
\put(80.00,20.00){\circle{60.00}}
\thicklines
\put(80.00,20.00){\line(-1,-1){22.5}}
\put(80.00,20.00){\makebox(0.00,0.00){$\bullet$}}
\put(90.00,10.00){\line(-1,-1){10}}
\put(90.00,10.00){\makebox(0.00,0.00){$\bullet$}}
\put(80.00,0.00){\line(0,-1){5}}
\put(170,0){
\thicklines
\put(80.00,20.00){\line(-1,-1){20}}
\put(80.00,20.00){\makebox(0.00,0.00){$\bullet$}}
\put(90.00,10.00){\line(-1,-1){10}}
\put(90.00,10.00){\makebox(0.00,0.00){$\bullet$}}
}
\end{picture}}
\end{center}
\caption{\label{prohlidka}%
Bracketed trees. The left picture shows an
element of $\Psi\Psi(E)(5)$ while the right picture shows the same element
interpreted, after erasing the braces indicated by thin cycles,
as an element of $\Psi(E)(5)$.
For simplicity, we did not show the
decoration of vertices by elements of $E$.}
\end{figure}
The triple unit
$\upsilon_E : E \to \Psi(E)$ identifies elements of $E$ with decorated
corollas:
\[
\unitlength 4mm
\linethickness{0.4pt}
\begin{picture}(20,6.5)(10,16.8)
\put(14,20){\makebox(0,0)[cc]{$E(n)\ni e \hskip 10pt \longleftrightarrow$ }}
\put(20,21){\line(0,-1){1}}
\put(20,22){\line(0,-1){1}}
\multiput(20,20)(-.0333333,-.0333333){60}{\line(0,-1){.0333333}}
\multiput(20,20)(-.033333,-.066667){30}{\line(0,-1){.066667}}
\multiput(20,20)(.0333333,-.0333333){60}{\line(0,-1){.0333333}}
\put(20,20){\makebox(0,0)[cc]{$\bullet$}}
\put(20,18){\makebox(0,0)[cc]{$\ldots$}}
\put(20,17){\makebox(0,0)[cc]{%
   $\underbrace{\rule{16mm}{0mm}}_{\mbox{\small $n$ inputs}}$}}
\put(20.5,20.5){\makebox(0,0)[lb]{$e$}}
\put(22,20){\makebox(0,0)[cl]{$\in \Psi(E)(n),\ n \geq 0$.}}
\end{picture}
\]

It is not difficult to verify that the above constructions indeed make
$\Psi$ a triple, compare~\cite[\S~II.1.12]{markl-shnider-stasheff:book}. Now we can
\underline{define}
non-unital operads as algebras over the triple $(\Psi,\mu,\upsilon)$.  The
advantage of this approach is that, by replacing $\TTT_n$
in~(\ref{msss}) by another category of trees or graphs, one may
obtain triples defining other types of operads and their generalizations.

\begin{figure}
\begin{center}
{
\unitlength=1.000000pt
\thinlines
\begin{picture}(180.00,100.00)(0.00,0.00)
\put(0.00,20.00){\line(1,0){180.00}}
\put(0.00,40.00){\line(1,0){180.00}}
\put(0.00,80.00){\line(1,0){180.00}}
\thicklines
\put(150.00,20.00){\makebox(0.00,0.00){$\bullet$}}
\put(130.00,20.00){\makebox(0.00,0.00){$\bullet$}}
\put(110.00,20.00){\makebox(0.00,0.00){$\bullet$}}
\put(60.00,20.00){\makebox(0.00,0.00){$\bullet$}}
\put(30.00,20.00){\makebox(0.00,0.00){$\bullet$}}
\put(130.00,40.00){\makebox(0.00,0.00){$\bullet$}}
\put(90.00,40.00){\makebox(0.00,0.00){$\bullet$}}
\put(50.00,40.00){\makebox(0.00,0.00){$\bullet$}}
\put(90.00,80.00){\makebox(0.00,0.00){$\bullet$}}
\put(90.00,80.00){\line(0,-1){40.00}}
\put(150.00,20.00){\line(0,-1){20.00}}
\put(150.00,20.00){\line(-1,-2){10.00}}
\put(130.00,40.00){\line(0,-1){20.00}}
\put(110.00,20.00){\line(1,-2){10.00}}
\put(110.00,20.00){\line(-1,-2){10.00}}
\put(30.00,20.00){\line(-1,-2){10.00}}
\put(30.00,20.00){\line(0,-1){20.00}}
\put(60.00,20.00){\line(1,-2){10.00}}
\put(60.00,20.00){\line(-1,-2){10.00}}
\put(30.00,20.00){\line(1,-2){10.00}}
\put(130.00,40.00){\line(-1,-1){20.00}}
\put(50.00,40.00){\line(1,-2){10.00}}
\put(90.00,80.00){\line(1,-1){80.00}}
\put(90.00,80.00){\line(-1,-1){80.00}}
\put(90.00,100.00){\line(0,-1){20.00}}
\end{picture}}
\end{center}
\caption{\label{maytree}A May's tree.}
\end{figure}
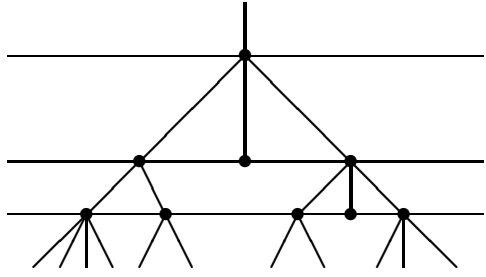

We have already seen in~(\ref{freeop}) that enlarging $\TTT_n$ into
$\UTTT_n$ by adding the exceptional tree, one gets the triple $\Gamma$
describing (unital) operads.  It is not difficult to see that
non-unital May's operads are related to the category
$\MTTT_n$\label{mma} of {\em May's trees\/} which are, by definition,
rooted trees whose vertices can be arranged into levels as in
Figure~\ref{maytree}.
Non-unital May's operads are then
algebras over the triple ${\rm M} : \sigmamod\to \sigmamod$ defined by
\[
{\rm M}(E)(n) := \colim{T \in \MTTT_n} E(T),\ n \geq 0.
\]
These observations are summarized in the first three lines of the
table in Figure~\ref{tables} on page~\pageref{tables}.

\section{Cyclic operads}
\label{sec5}

In the following two sections
we use the approach developed in Section~\ref{sec4} to
introduce cyclic and modular operads.
We recalled, in Example~\ref{Jitka_se_ozvala}, the operad
$\hatcalM_0= \coll {\hatcalM_0}{0}$ of
Riemann spheres with parametrized labeled holes.
Each $\hatcalM_0(n)$ was a
right $\Sigma_n$-space, with the operadic right $\Sigma_n$-action
permuting the labels $\rada 1n$ of the holes $\Rada u1n$.
But each $\hatcalM_0(n)$ obviously
admits a higher type of symmetry which interchanges
labels $\rada 0n$ of {\em all\/} holes, including the label of the
`output' hole $u_0$.  Another example admitting a similar higher symmetry is the
configuration (non-unital)
operad $\barcalM_0 = \{\barcalM_0(n)\}_{n \geq 2}$.

These examples indicate that, for some operads, there is no
clear distinction between `inputs' and the `output.' Cyclic operads,
introduced by E.~Getzler and
M.M.~Kapranov in~\cite{getzler-kapranov:CPLNGT95}, formalize this
phenomenon. They are, roughly speaking, operads with an extra symmetry
that interchanges the output with one of the inputs. Let us recall
some notions necessary to give a~precise definition.

We remind the reader that in this section, as well as everywhere in
this article, main definitions are formulated over the underlying
category of $\bfk$-modules, where $\bfk$ is a commutative associative
unital ring. However, for some constructions, we will require $\bfk$ to
be a {\em field\/}; we will indicate this as usual by speaking about
{\em vector spaces\/} instead of $\bfk$-modules.

Let $\Sigma_n^+$ be the permutation group of the set $\{\rada
0n\}$. The group $\Sigma^+_n$ is, of course, non-canonically
isomorphic to the symmetric group $\Sigma_{n+1}$.  We identify
$\Sigma_n$ with the subgroup of $\Sigma^+_n$ consisting of
permutations $\sigma \in \Sigma^+_n$ such that $\sigma(0)= 0$.  If
$\tau_n \in \Sigmap_n$ denotes the cycle $(\rada 0n)$,
that is, the permutation with $\tau_n(0)=1,\ \tau_n(1)=2,\
\ldots,\tau_n(n)=0$, then $\tau_n$ and $\Sigma_n$ generate
$\Sigmap_n$.

Recall that a {\em cyclic $\Sigma$-module\/}
or a {\em $\Sigma^+$-module\/} is a sequence $W = \{W(n)\}_{n\geq 0}$
such that each $W(n)$ is a (right) $\bfk[\Sigma^+_n]$-module.  Let
$\sigmamodp$ denote the category of cyclic $\Sigma$-modules.
As (ordinary) operads were $\Sigma$-modules with an
additional structure, cyclic operads are $\Sigma^+$-modules with
an additional structure.

We will also need the following `cyclic' analog of~(\ref{m-m}):
if $X$ is a set with $n+1$ elements and $W \in \sigmamodp$, then
\begin{equation}
\label{kkkkk}
W\((X\)) :=
W(n) \times_{\Sigma_n^+} \Bij(\en^+,X),
\end{equation}
where $\en^+ := \{0,\ldots,n\}$, $n \geq 0$.
Double brackets in $W\((X\))$ remind us that the $n$th
piece of the cyclic $\Sigma$-module $W = \coll W0$
is applied on a set with $n+1$ elements, using the extended
$\Sigma_n^+$-symmetry. Therefore
\[
W\((\{\rada 0n\}\)) \cong W(n)\ \mbox { while }\
W(\{\rada 0n\}) \cong W(n+1),\ n \geq 0.
\]

Pasting schemes for cyclic operads are {\em cyclic (leg-) labeled
$n$-trees\/}, by which we mean \underline{un}-\underline{rooted} trees as on
page~\pageref{prohlidku_mam_za_sebou!}, with legs labeled by the set
$\{\rada 0n\}$.  An example of such a tree is given in
Figure~\ref{pozitri_prohlidka}. Since we do not assume a choice of the
root, the edges of a cyclic tree $C$ are not directed and it does not
make sense to speak about inputs and the output of a vertex $v \in
\vert(C)$.  Let $\TTT^+_n$ be the category of cyclic labeled $n$-trees
and their bijections.
\begin{figure}[t]
\begin{center}
{
\unitlength=1.000000pt
\begin{picture}(130.00,110.00)(0.00,0.00)
\thicklines
\put(20.00,80.00){\makebox(0.00,0.00){$8$}}
\put(90.00,110.00){\makebox(0.00,0.00){$9$}}
\put(50.00,0.00){\makebox(0.00,0.00){$0$}}
\put(120.00,10.00){\makebox(0.00,0.00){$7$}}
\put(130.00,60.00){\makebox(0.00,0.00){$6$}}
\put(30.00,20.00){\makebox(0.00,0.00){$5$}}
\put(0.00,60.00){\makebox(0.00,0.00){$4$}}
\put(80.00,0.00){\makebox(0.00,0.00){$3$}}
\put(120.00,80.00){\makebox(0.00,0.00){$2$}}
\put(40.00,110.00){\makebox(0.00,0.00){$1$}}
\put(90.00,40.00){\makebox(0.00,0.00){$\bullet$}}
\put(50.00,40.00){\makebox(0.00,0.00){$\bullet$}}
\put(30.00,50.00){\makebox(0.00,0.00){$\bullet$}}
\put(60.00,90.00){\makebox(0.00,0.00){$\bullet$}}
\put(50.00,80.00){\makebox(0.00,0.00){$\bullet$}}
\put(90.00,80.00){\makebox(0.00,0.00){$\bullet$}}
\put(70.00,60.00){\makebox(0.00,0.00){$\bullet$}}
\put(30.00,50.00){\line(-2,1){20.00}}
\put(30.00,50.00){\line(0,-1){20.00}}
\put(70.00,60.00){\line(-4,-1){40.00}}
\put(70.00,60.00){\line(-2,-5){20.00}}
\put(50.00,80.00){\line(1,1){10.00}}
\put(90.00,40.00){\line(1,-1){20.00}}
\put(90.00,40.00){\line(3,2){30.00}}
\put(90.00,40.00){\line(-1,-3){10.00}}
\put(50.00,80.00){\line(-1,0){20.00}}
\put(40.00,100.00){\line(1,-2){10.00}}
\put(90.00,80.00){\line(1,0){20.00}}
\put(90.00,100.00){\line(0,-1){20.00}}
\put(50.00,40.00){\line(1,1){40.00}}
\put(50.00,80.00){\line(1,-1){40.00}}
\end{picture}}
\end{center}
\caption{\label{pozitri_prohlidka}%
A cyclic labeled tree from $\protect\TTT^+_{9}$.}
\end{figure}

For a cyclic $\Sigma$-module $W$ and a cyclic labeled tree $T$
we have the following cyclic version of the
product~(\ref{ale_jeste_mne_cekaji_odbery})
\[
W\((T\)):=\bigotimes_{v \in \vert(T)}W\((\edge(v)\)).
\]
The conceptual difference between~(\ref{ale_jeste_mne_cekaji_odbery})
and the above formula is that
instead of the set $\In(v)$ of incoming edges of a vertex
$v$ of a rooted tree, here we use the set $\edge(v)$ of {\em all\/} edges
incident with $v$.
Let, finally, $\Psi_+ : \sigmamodp \to \sigmamodp$ be the functor
\begin{equation}
\label{maminka_zase_v_nemocnici}
\Psi_+(W)(n) := \colim{T \in \TTT^+_n} W\((T\)),\ n \geq 0,
\end{equation}
equipped with the triple structure of `forgetting the braces' similar\label{[]}
to that reviewed on page~\pageref{vecer_zavolam_Jitce}. We will use
also the `extended' triple $\Gamma_+ : \sigmamodp \to \sigmamodp$,
\[
\Gamma_+(W)(n) := \colim{T \in \UTTT^+_n} W\((T\)),\ n \geq 0,
\]
where $\UTTT_n^+$ is the obvious extension of the category $\TTT_n^+$ by the
exceptional tree \hskip .5em
\raisebox{-.2em}{\rule{.8pt}{1.1em}}~.

\begin{definition}
A~{\em cyclic\/} (resp.~{\em non-unital cyclic\/}) {\em operad\/} is
an algebra over the triple $\Gamma_+$ (resp.~the triple
$\Psi_+$) introduced above.
\end{definition}

In the following proposition, which
slightly improves~\cite[Theorem~2.2]{getzler-kapranov:CPLNGT95},
$\tau_n \in \Sigmap_n$ denotes the cycle $(\rada 0n)$.

\begin{proposition}
\label{zitra_s_Jitkou_na_alergologii}
A non-unital cyclic operad is the same as a non-unital operad $\eucalC =
\coll{\eucalC}0$ (Definition~\ref{b3}) such that the right
$\Sigma_n$-action on $\eucalC(n)$ extends, for each $n\geq 0$, to an action
of~$\Sigmap_n$ with the property that
for $p\in \eucalC(m)$ and $q\in \eucalC(n)$, $1\leq i \leq m$, $n \geq
0$, the composition maps satisfy
\[
(p \circ_i q) \tau_{m+n-1}
=
\cases{(q \tau_n) \circ_n (p \tau_m)}{if $i=1$, and}
      {(p \tau_m) \circ_{i-1} q}
      {for $2 \leq i \leq m$.}
\]
The above structure is a (unital) cyclic operad
if moreover there exists a $\Sigma_1^+$-invariant operadic
unit $e\in \eucalC(1)$.
\end{proposition}

Proposition~\ref{zitra_s_Jitkou_na_alergologii} gives a biased
definition of cyclic operads whose obvious modification
(see~\cite[Definition~II.5.2]{markl-shnider-stasheff:book})
makes sense in an arbitrary symmetric monoidal
category. We can therefore speak about topological cyclic operads,
differential graded cyclic operads,
simplicial cyclic operads~\&c. Observe that there are no
{\em non-unital cyclic May's operads\/} because it does not make sense
to speak about levels in trees without a choice of the root.

\begin{example}
\label{cyclic-endomorphism-operad}
{\rm\
Let $V$ be a finite dimensional vector space and
$B: V\ot V \to \bfk$ a
nondegenerate symmetric bilinear form.
The form $B$ induces the identification
\[
\Hom{V^{\ot n}}{V} \ni f \longmapsto \widehat B(f):=
B(-,f(-)) \in \Hom{V^{\ot (n+1)}}{\bfk}
\]
of the spaces of linear maps.  The standard right $\Sigmap_n$-action
\[
\widehat B(f)\sigma(\Rada v0n) =
\widehat B(f)(\Rada v{\sigma^{-1}(0)}{\sigma^{-1}(n)}),\
\sigma \in \Sigmap_n,\
\Rada v0n \in V,
\]
defines, via this identification, a right $\Sigmap_n$-action on
$\Hom{V^{\ot n}}{V}$, that is, on the $n$th piece of the endomorphism
operad $\End_V = \coll{\End_V}0$ recalled in Example~\ref{end}.  It is
easy to show that, with the above action, $\End_{V}$ is a cyclic
operad in the monoidal category of vector spaces, called the {\em
cyclic endomorphism operad\/}
of the pair $V = (V,B)$. The biased definition of
cyclic operads given in Proposition~\ref{zitra_s_Jitkou_na_alergologii}
can be read off from this example.
}
\end{example}

\begin{example}
{\rm
We saw in Example~\ref{zitra_strasne_vedro} that a unital operad
$\calA = \coll \calA0$ such that $\calA(n) = 0$ for $n \not=1$ is the
same as a unital associative algebra. Similarly, it can be easily
shown that a cyclic operad 
$\calC = \coll \calC0$ satisfying $\calC(n) = 0$ for $n \not=1$ is the
same as a unital associative algebra $A$ with a linear involutive
antiautomorphism, by which we mean a $\bfk$-linear map 
${}^* : A \to A$ such~that
\[
(ab)^* = b^* a^*,\ 
(a^*)^* = a\
\mbox { and }\ 1^* = 1,
\]
for arbitrary $a,b \in A$.  
}  
\end{example}

Let $\calP = \Gamma(E)/(R)$ be a quadratic operad as in
Definition~\ref{je_tady_neporadek}.
The action of $\Sigma_2$ on
$E$ extends to an action of $\Sigma^+_2$, via the sign
representation
$\sgn : \Sigma^+_2 \to \{ \pm 1\} = \Sigma_2$. It can be easily verified
that this action induces a cyclic operad structure on the free operad
$\Gamma(E)$. In particular, $\Gamma(E)(3)$ is a right
$\Sigma^+_3$-module.

\begin{definition}
\label{I_need_a_day_off}
We say that the operad $\calP$
is a {\em cyclic quadratic operad\/}\index{cyclic quadratic operad}
if, in the above presentation, $R$ is a $\Sigma^+_3$-invariant
subspace of $\Gamma(E)(3)$.
\end{definition}

If the condition of the above definition is satisfied,
$\calP$ has a natural induced cyclic operad structure.

\begin{example}
{\rm
By~\cite[Proposition~3.6]{getzler-kapranov:CPLNGT95},
all quadratic operads generated by a one-dimensional space are
cyclic quadratic, therefore the operads
$\Lie$ and $\Com$ are cyclic quadratic. Also the operads
$\Ass$ and the operad $\Poiss$ for Poisson algebras
are cyclic quadratic~\cite[Proposition~3.11]{getzler-kapranov:CPLNGT95}.
A~surprisingly simple
operad which is cyclic and quadratic, but not cyclic quadratic, is
constructed in~\cite[Remark~15]{markl-remm}.

The operad $\hatcalM_0$ of Riemann spheres with labeled
punctures reviewed in Example~\ref{Jitka_se_ozvala} is a topological
cyclic operad. The configuration operad $\barcalM_0$ recalled on
page~\pageref{New_Orleans} is a non-unital topological cyclic operad.
Important examples of non-cyclic operads are the operad $\pre-Lie$ for
pre-Lie algebras~\cite[Section~3]{markl-remm}
and the operad $\Leib$ for Leibniz
algebras~\cite[\S~3.15]{getzler-kapranov:CPLNGT95}.
}
\end{example}

Let $\eucalC$ be an operad, $\alpha : \eucalC(n) \ot \otexp Vn \to
V$, $n \geq 0$, a $\eucalC$-algebra with the underlying vector space~$V$ as
in Proposition~\ref{pro:alg} and $B : V\ot V \to U$ a bilinear form on
$V$ with values in a vector space $U$. We can form a map
\begin{equation}
\label{az}
{\widetilde B}(\alpha) : \eucalC(n) \ot \otexp V{(n+1)} \to U, \ n
\geq 0,
\end{equation}
by the formula
\[
{\widetilde B}(\alpha)(c \ot v_0 \ot \cdots v_n) :=
B(v_0,\alpha(c \ot v_1 \ot \cdots v_n)),\ c \in \eucalC(n),\
\Rada v0n \in V.
\]
Suppose now that the operad $\eucalC$ is cyclic, in particular, that
each $\eucalC(n)$ is a right $\Sigma_n^+$-module. We say that
the bilinear form $B :
V\ot V \to U$ is {\em
invariant\/}~\cite[Definition~4.1]{getzler-kapranov:CPLNGT95},
if the maps ${\widetilde B}(\alpha)$ in~(\ref{az}) are, for each $n \geq
0$, invariant under the diagonal action of $\Sigma^+_n$ on $\eucalC(n)
\ot \otexp V{(n+1)}$. We leave as an exercise to verify that the
invariance of ${\widetilde B}(\alpha)$ for $n=1$ together with the
existence of the operadic unit implies that $B$ is symmetric,
\[
B(v_0,v_1) = B(v_1,v_0),\ v_0,v_1 \in V.
\]

\begin{definition}
A {\em cyclic algebra\/} over a cyclic operad $\eucalC$ is a
$\eucalC$-algebra structure on a vector space~$V$
together with a nondegenerate invariant
bilinear form $B : V\ot V \to \bfk$.
\end{definition}

By~\cite[Proposition~II.5.14]{markl-shnider-stasheff:book},
a cyclic algebra is the same as
a cyclic operad homomorphism $\eucalC \to \End_V$, where $\End_V$ is
the cyclic endomorphism operad of the pair $(V,B)$ recalled in
Example~\ref{cyclic-endomorphism-operad}.

\begin{example}
{\rm
A cyclic algebra over the cyclic operad $\Com$ is a {\em Frobenius algebra\/},
that is, a structure consisting of a
commutative associative multiplication $\cdot : V \ot V \to V$
as in Example~\ref{Jitce_budu_volat_az_v_patek}
together with a non-degenerate symmetric
bilinear form $B : V \ot V \to \bfk$, invariant in the sense
that
\[
B(a \cdot b,c) = B(a, b \cdot c),\ \mbox { for all }\ a,b,c \in V.
\]
Similarly, a cyclic Lie algebra is given by a Lie bracket $[-,-] : V
\ot V \to V$ and a  non-degenerate symmetric
bilinear form $B : V \ot V \to \bfk$ satisfying
\[
B([a,b],c) = B(a,[b,c]),\ \mbox { for }\ a,b,c \in V.
\]
}\end{example}

For algebras over cyclic operads, one may introduce cyclic
cohomology that generalizes the classical cyclic cohomology of
associative 
algebras~\cite{cartier:bourbaki,loday-quillen:84,tamarkin-tsygan:LMP01}
as the non-abelian derived functor of the universal bilinear
form~\cite{getzler-kapranov:CPLNGT95},~\cite[Proposition~II.5.26]%
{markl-shnider-stasheff:book}.
Let us close this section by mentioning two
examples of operads with other types of higher symmetries. The
symmetry required for {\em anticylic operads\/} differs from the
symmetry of cyclic operads by the
sign~\cite[Definition~II.5.20]{markl-shnider-stasheff:book}.
{\em Dihedral operads\/} exhibit
a symmetry governed by the dihedral
groups~\cite[Definition~16]{markl-remm}.

\section{Modular operads}
\label{sec6}

Let us consider again the $\Sigma^+$-module $\hatcalM_0 =
\coll{{\hatcalM}_0}0$ of Riemann spheres with punctures.  We saw that
the operation $M,N \mapsto M \circ_i N$ of sewing the $0$th hole of
the surface $N$ to the $i$th hole of the surface $M$ defined on
$\hatcalM_0$ a cyclic operad structure.
One may generalize this operation by defining, for $M \in
{\hatcalM}_0(m)$, $N \in {\hatcalM}_0(n)$, $0 \leq i \leq m$, $0 \leq
j \leq n$, the element $M \circij ij  N \in {\hatcalM}_0(m+n-1)$ by
sewing the $j$th hole of $M$ to the $i$th hole of $N$. Under this
notation, $\circ_i = \circij i0$.
In the same manner, one may consider a single surface $M \in
\hatcalM_0(n)$, choose labels $i,j$, $0 \leq i\not=j \leq n$, and sew
the $i$th hole of $M$ along the $j$th hole of the {\em same\/}
surface. The result is a new surface $\xi_{\{i,j\}}(M)$, with $n-2$
holes and genus $1$.

This leads us to the system $\hatcalM =
\{\hatcalM(g,n)\}_{g \geq 0, n \geq -1}$, where $\hatcalM\lzavorkaobyc
g,n \rzavorkaobyc$ denotes now the moduli space of genus $g$ Riemann
surfaces with $n+1$ holes.  Observe that we include $\hatcalM(g,n)$
also for $n = -1$; $\hatcalM(g,-1)$ is the moduli space of Riemann
surfaces of genus $g$.  The operations $\circij ij$ and
$\xi_{\{i,j\}}$ act on $\hatcalM$. Clearly, for $M \in
\hatcalM\lzavorkaobyc g,m\rzavorkaobyc$ and $N \in \hatcalM(h,n)$, $0
\leq i \leq m$, $0 \leq j \leq n$ and $g,h \geq 0$,
\begin{equation}
\label{prave_jsem_mluvil_s_Jitkou}
M\circij ij N \in
\hatcalM\lzavorkaobyc g+h,m+n-1\rzavorkaobyc
\end{equation}
and, for $m \geq 1$ and $g \geq 0$,
\begin{equation}
\label{prave_jsem_mluvil_s_Jitkou.}
\xi_{\{i,j\}}(M) \in \hatcalM\lzavorkaobyc g+1,m-2 \rzavorkaobyc.
\end{equation}
A particular case of~(\ref{prave_jsem_mluvil_s_Jitkou}) is the
non-operadic composition
\begin{equation}
\label{Prave_jsem_mluvil_s_Jitkou}
\circij 00 : \hatcalM(g,0) \times \hatcalM(h,0) \to \hatcalM(g+h,-1),\
g,h \geq 0.
\end{equation}

Modular operads are abstractions of the above structure satisfying a
certain additional stability condition.  The following definitions,
taken from~\cite{getzler-kapranov:CompM98}, are made for the category
of $\bfk$-modules, but they can be easily generalized to an arbitrary
symmetric monoidal category with finite colimits,
whose monoidal product $\odot$ is distributive over colimits.
Let us introduce the underlying category for modular
operads.

A {\em modular $\Sigma$-module\/} is a sequence $\calE = \collmod
{{\calE}}$ of $\bfk$-modules such that each
${{\calE}}\lzavorkaobyc g,n\rzavorkaobyc$ has a right
$\bfk[\Sigma^+_n]$-action.  We say that $\calE$
is {\em stable}\index{stable modular $\Sigma$-module} if
\begin{equation}
\label{Cupid_chain}
\calE\lzavorkaobyc g,n \rzavorkaobyc = 0\ \mbox { for } 2g+n-1 \leq 0
\end{equation}
and denote $\MCol$ the category of stable modular $\Sigma$-modules.

Stability~(\ref{Cupid_chain}) says that $\calE\lzavorkaobyc g,n
\rzavorkaobyc$ is trivial for $(g,n) = (0,-1),\ (1,-1),\ (0,0)$ and
$(0,1)$. We will sometimes express the stability of $\calE$ by
writing $\calE = \stablecoll{\calE}$, where
\[
\frakS := \{(g,n)\, |\ g \geq 0,\ n \geq -1 \mbox { and } 2g+n-1 >0\}.
\]
Recall that a genus $g$ Riemann surface with $k$
marked points is stable if it does not admit infinitesimal
automorphisms. This happens if and only if  $2(g-1) + k > 0$, that is,
excluded is the torus with no marked points and the
sphere with less than three marked points. Thus the stability property
of modular $\Sigma$-modules is analogous to the stability of Riemann
surfaces.

Now we introduce graphs that serve as pasting schemes for modular operads.
The naive notion of a graph as we have used it up to this point
is not subtle enough; we need to replace it by a more sophisticated:

\begin{definition}
\label{graph-def}
A {\em graph\/} $\Gamma$ is a finite set $\Flag(\Gamma)$
(whose elements are called {\em flags\/} or {\em half-edges\/})
together with an involution
$\sigma$ and a partition $\lambda$.
The {\em vertices\/} $\vert(\Gamma)$
of a graph $\Gamma$ are the blocks of the partition
$\lambda$, we assume also that the number of these blocks is finite.
The {\em edges\/}
$\Edg(\Gamma)$ are pairs of flags forming a two-cycle of $\sigma$. The
{\em legs\/} $\leg(\Gamma)$ are the fixed points
of $\sigma$.
\end{definition}

We also denote by $\edge(v)$ the flags belonging to the block $v$ or,
in common speech, half-edges adjacent to the vertex $v$.  We say that
graphs $\Gamma_1$ and $\Gamma_2$ are {\em isomorphic\/} if there
exists a set isomorphism $\varphi : \Flag(\Gamma_1) \to
\Flag(\Gamma_2)$ that preserves the partitions and commutes with the
involutions. We may associate to a graph $\Gamma$ a finite
one-dimensional cell complex $|\Gamma|$, obtained by taking one copy
of $[0,\frac 12]$ for each flag, a point for each block of the
partition, and imposing the following equivalence
relation: The points $0\in [0,\frac 12]$ are identified for all flags
in a block of the partition $\lambda$ with the point corresponding to
the block, and the points $\frac 12 \in
[0,\frac 12]$ are identified for pairs of flags exchanged by the
involution $\sigma$.

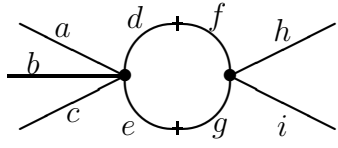
\begin{figure}[t]
\begin{center}
\setlength{\unitlength}{0.00023in}%
\begin{picture}(7524,3100)(2089,-4723)
\thicklines
\put(6001,-3361){\oval(2400,2400)}
\put(7201,-3361){\line( 2, 1){2400}}
\put(7201,-3361){\line( 2,-1){2400}}
\put(4801,-3361){\line(-2, 1){2400}}
\put(4801,-3361){\line(-2,-1){2400}}
\put(4801,-3361){\line(-1, 0){2700}}
\put(6001,-2011){\line( 0,-1){300}}
\put(6001,-4411){\line( 0,-1){300}}
\put(3376,-2311){\makebox(0,0){$a$}}
\put(2701,-3061){\makebox(0,0){$b$}}
\put(3601,-4261){\makebox(0,0){$c$}}
\put(5026,-2000){\makebox(0,0){$d$}}
\put(6901,-2000){\makebox(0,0){$f$}}
\put(4876,-4550){\makebox(0,0){$e$}}
\put(6976,-4550){\makebox(0,0){$g$}}
\put(8401,-2311){\makebox(0,0){$h$}}
\put(8401,-4486){\makebox(0,0){$i$}}
\put(4801,-3361){\makebox(0,0){$\bullet$}}
\put(7201,-3361){\makebox(0,0){$\bullet$}}
\end{picture}
\end{center}
\caption{\label{sput}%
The sputnik.}
\end{figure}

We call $|\Gamma|$ the {\em geometric
realization\/} of $\Gamma$.
Observe that empty blocks of the partition generate isolated vertices
in the geometric realization.
We will usually make no distinction
between the graph and its geometric realization.
As an example (taken from~\cite{getzler-kapranov:CompM98}), consider
the graph with $\{a,b,\ldots,i\}$ as the set of flags, the involution
$\sigma = (df)(eg)$ and the partition $\{a,b,c,d,e\} \cup
\{f,g,h,i\}$. The geometric realization of this graph is the `sputnik' in
Figure~\ref{sput}.

Let us introduce labeled versions of the above notions.  A {\em
(vertex-) labeled graph\/} is a connected graph $\Gamma$ together with
a map $g$ (the {\em genus map\/}) from $\vert(\Gamma)$ to the set
$\{0,1,2,\ldots\}$.  Labeled graphs $\Gamma_1$ and $\Gamma_2$ are
isomorphic if there exists an isomorphism preserving the labels of the
vertices.  The {\em genus\/} $g(\Gamma)$ of a labeled graph $\Gamma$
is defined by
\begin{equation}\
\label{Good_Friday}
g(\Gamma) := b_1(\Gamma)+\sum_{v\in \vert(\Gamma)} g(v),
\end{equation}
where $b_1(\Gamma) := \dim H_1(|\Gamma|)$ is the first Betti number of
the graph $|\Gamma|$, i.e.~the number of independent circuits of
$\Gamma$.  A graph $\Gamma$ is {\em stable\/} if
\[
2(g(v)-1) + |\edge(v)| > 0,
\]
at each vertex $v \in \vert(\Gamma)$.

For $g \geq 0$ and $n \geq -1$, let $\Graphcat gn$ be the groupoid
whose objects are pairs $(\Gamma,\ell)$ consisting of a stable
(vertex-) labeled graph $\Gamma$ of genus $g$ and an isomorphism $\ell
: \Leg(\Gamma) \to \n n$ labeling the legs of $\Gamma$ by elements of
$\n n$.  Morphisms of $\Graphcat gS$ are isomorphisms of
vertex-labeled graphs preserving the labeling of the legs. The
stability implies, via an elementary combinatorial topology that,
for each fixed $g \geq 0$ and $n \geq -1$, there is only a finite
number of isomorphism classes of stable graphs $\Gamma \in \Graphcat
gn$, see~\cite[Lemma~2.16]{getzler-kapranov:CompM98}.

We will also need the following obvious generalization
of~(\ref{kkkkk}): if
$\calE = \{\calE(g,n)\}_{g \geq 0, n \geq -1}$ is a modular
$\Sigma$-module and $X$ a set with
$n+1$ elements, then
\begin{equation}
\label{Jtk}
\calE\((g,X\)) :=
\calE(g,n) \times_{\Sigma_n^+} \Bij(\en^+,X),\  g \geq 0,\ n \geq -1.
\end{equation}
For a modular $\Sigma$-module $\calE = \collmod{\calE}$  and a
labeled graph $\Gamma$, let ${{\calE}}\((\Gamma\))$ be the
product
\begin{equation}
\label{Cosmo}
{{\calE}}\((\Gamma\)): = \bigotimes_{v\in\vert(\Gamma)}
{{\calE}}\((g(v),\edge(v)\)) .
\end{equation}
Evidently, the correspondence $\Gamma \mapsto \calE\((\Gamma\))$ defines a
functor from the category $\Graphcat gn$ to the category of
$\bfk$-modules and their isomorphisms.
We may thus define an endofunctor $\Modtriple$ on
the category
$\MCol$ of stable modular $\Sigma$-modules by the formula
\[
\Modtriple {{\calE}}(g,n) : =
\colim{{\Gamma\in \Graphcat gn}}{{\calE}}\((\Gamma\)),\
g \geq 0,\ n \geq -1.
\]

Choosing a representative for each
isomorphism class in $\Graphcat gn$, one obtains the identification
\begin{equation}
\label{Chrastilek}
\Modtriple {{\calE}}(g,n) \cong
\bigoplus_{{[\Gamma]\in\{\mbox{\scriptsize $\Graphcat gn$}\}}}
{{\calE}}\((\Gamma\))_{{\rm Aut}(\Gamma)},\
g \geq 0,\ n \geq -1,
\end{equation}
where $\{\Graphcat gn\}$ is the set of isomorphism classes of
objects of the groupoid $\Graphcat gn$ and the subscript ${\rm
Aut}(\Gamma)$ denotes the space of
coinvariants. Stability~(\ref{Cupid_chain}) implies that the summation
in the right-hand side of~(\ref{Chrastilek}) is
finite. Formula~(\ref{Chrastilek}) generalizes~(\ref{draw}) which does
not contain coinvariants because there are no nontrivial
automorphisms of leaf-labeled trees. On the other hand,
stable labeled graphs with nontrivial automorphisms are abundant,
an example can be easily constructed from the graph in Figure~\ref{sput}.
The functor $\Modtriple$ carries a triple structure of `erasing
the braces' similar to the one used on
pages~\pageref{vecer_zavolam_Jitce} and~\pageref{[]}.

\begin{definition}
A {\em modular operad\/} is an algebra over the triple $\Modtriple :
\MMod \to \MMod$.
\end{definition}

Therefore a modular operad is a stable modular $\Sigma$-module $\calA
= \stablecoll{\calA}$ equipped with operations that determine coherent
contractions along stable modular graphs.  Observe that the stability
condition is built firmly into the very definition. Very
crucially, modular operads {\em do not have units\/}, because such a unit
ought to be an element of the space $\calA\lzavorkaobyc 0,1
\rzavorkaobyc$ which is empty, by~(\ref{Cupid_chain}).

One can easily introduce un-stable modular operads and their unital versions,
but the main motivating example reviewed below is stable.
We will consider an extension of the
Grothendieck-Knudsen configuration operad $\barcalM_0 = \coll {\barcalM_0}2$
consisting of moduli spaces of stable curves of arbitrary genera
in the sense of the following generalization of Definition~\ref{jItka}:

\begin{definition}
\label{leze_na_mne_chripka}
A {\em stable $(n+1)$-pointed curve\/}, $n \geq 0$, is a connected
complex projective curve $C$ with at most nodal singularities,
together with a `marking' given by a choice $\Rada x0n \in C$ of
smooth points.  The stability means, as usual, that there are no infinitesimal
automorphisms of $C$ fixing the marked points and double points.
\end{definition}

The stability in Definition~\ref{leze_na_mne_chripka} is equivalent to
saying that each smooth component of $C$
isomorphic to the complex projective space ${\mathbb {CP}}^1$ has at least
three special points and that each smooth component isomorphic to the
torus has at least one special point, where by a special point we mean either a
double point or a node.

\begin{figure}
{
\unitlength=.52pt
\begin{picture}(498.50,200.50)(0.00,0.00)
\thicklines
\put(80,0){
\put(498.50,70.50){\makebox(0.00,0.00){$x_0$}}
\put(448.50,70.50){\makebox(0.00,0.00){$x_2$}}
\put(348.50,70.50){\makebox(0.00,0.00){$x_1$}}
\put(250,120){\makebox(0,0){Dual graph $\Gamma(C)$:}}
\put(399.50,80.50){\makebox(0.00,0.00){$\bullet$}}
\put(488.50,130.50){\makebox(0.00,0.00)[rb]{$a_4$}}
\put(448.50,130.50){\makebox(0.00,0.00)[b]{$a_3$}}
\put(398.50,130.50){\makebox(0.00,0.00)[b]{$a_2$}}
\put(358.50,130.50){\makebox(0.00,0.00)[lb]{$a_1$}}
\put(415.50,88.50){\makebox(0.00,0.00){$a_5$}}
\put(498.50,120.50){\makebox(0.00,0.00){$\bullet$}}
\put(448.50,120.50){\makebox(0.00,0.00){$\bullet$}}
\put(398.50,120.50){\makebox(0.00,0.00){$\bullet$}}
\put(348.50,120.50){\makebox(0.00,0.00){$\bullet$}}
\put(448.50,120.50){\line(0,-1){40.00}}
\put(498.50,120.50){\line(0,-1){40.00}}
\put(348.50,120.50){\line(0,-1){40.00}}
\put(418.50,170.50){\line(1,0){10.00}}
\bezier{300}(428.50,170.50)(498.50,170.50)(498.50,120.50)
\bezier{300}(348.50,120.50)(348.50,170.50)(418.50,170.50)
\qbezier(408.50,70.50)(408.50,80.50)(398.50,80.50)
\qbezier(398.50,60.50)(408.50,60.50)(408.50,70.50)
\qbezier(388.50,70.50)(388.50,60.50)(398.50,60.50)
\qbezier(398.50,80.50)(388.50,80.50)(388.50,70.50)
\put(378.50,120.50){\line(1,0){120.00}}
\put(348.50,120.50){\line(0,1){0.00}}
\put(378.50,120.50){\line(-1,0){30.00}}
\put(368.50,120.50){\line(1,0){10.00}}
\put(398.50,120.50){\line(0,-1){40.00}}
}
\put(-20,0){
\put(59.50,44.00){\makebox(0.00,0.00)[t]{$x_2$}}
\put(89.50,160.00){\makebox(0.00,0.00)[b]{$x_1$}}
\put(48.50,110.50){\makebox(0.00,0.00)[l]{$x_0$}}
\put(-30,120){\makebox(0,0){Curve $C$:}}
\put(199.00,160.00){\makebox(0.00,0.00)[b]{$A_5$}}
\put(9.00,21.00){\makebox(0.00,0.00){$A_4$}}
\put(169.00,0.00){\makebox(0.00,0.00)[l]{$A_3$}}
\put(168.50,179.50){\makebox(0.00,0.00)[b]{$A_2$}}
\put(-10.00,204.00){\makebox(0.00,0.00){$A_1$}}
\put(180.00,88.00){\makebox(0.00,0.00){$\bullet$}}
\put(121.50,87.00){\makebox(0.00,0.00){$\bullet$}}
\put(143.00,140){\makebox(0.00,0.00){$\bullet$}}
\put(80.50,149.50){\makebox(0.00,0.00){$\bullet$}}
\put(68.50,51.00){\makebox(0.00,0.00){$\bullet$}}
\put(110.00,35.00){\makebox(0.00,0.00){$\bullet$}}
\put(26.50,60.00){\makebox(0.00,0.00){$\bullet$}}
\put(35.00,114.50){\makebox(0.00,0.00){$\bullet$}}
\put(35.00,171.00){\makebox(0.00,0.00){$\bullet$}}
\bezier{300}(8.50,190.50)(78.50,130.50)(158.50,140.50)
\bezier{300}(18.50,30.50)(48.50,130.50)(28.50,200.50)
\bezier{300}(8.50,60.50)(68.50,60.50)(158.50,10.50)
\bezier{300}(158.50,170.50)(108.50,80.50)(108.50,10.50)
\bezier{300}(178.50,100.50)(178.50,120.50)(188.50,150.50)
\bezier{300}(198.50,60.50)(178.50,70.50)(178.50,100.50)
\bezier{300}(248.50,60.50)(228.50,40.50)(198.50,60.50)
\bezier{300}(248.50,80.50)(253.50,70.50)(248.50,60.50)
\bezier{300}(198.50,90.50)(238.50,100.50)(248.50,80.50)
\bezier{300}(88.50,90.50)(148.50,80.50)(198.50,90.50)
}
\end{picture}}
\caption{A stable curve and its dual graph.
The curve $C$ on the left has five components $A_i$, $1 \leq i \leq
5$, and three marked points
$x_0$, $x_1$ and $x_2$. The dual
graph $\Gamma(C)$ on the right
has five vertices $a_i$, $1 \leq i \leq
5$, corresponding to the components of the curve
and three legs labeled by the marked points.
\label{Domecek}}
\end{figure}
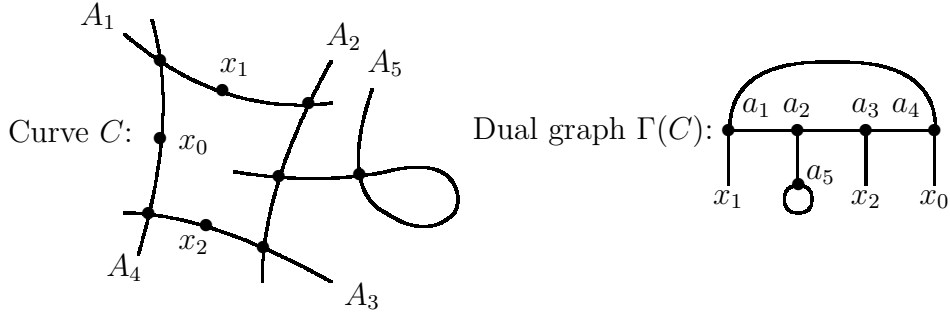

The {\em dual graph\/} $\Gamma =\Gamma(C)$ of a stable $(n+1)$-pointed
curve $C=(C,x_0,\dots,x_n)$ is a labeled graph whose vertices are the
components of $C$, edges are the nodes and its legs are the points
$\{x_i\}_{0\leq i \leq n}$.  An edge $e_y$ corresponding to a nodal
point $y$ joins the vertices corresponding to the components
intersecting at $y$. The vertex $v_K$ corresponding to a branch $K$ is
labeled by the genus of the normalization of $K$. See~\cite[page
23]{hartshorne:book} for the normalization and recall that a curve is
normal if and only if it is nonsingular.
The construction of $\Gamma(C)$ from a
curve $C$ is visualized in Figure~\ref{Domecek}.

Let us denote by $\barcalM_{g,n+1}$ the coarse moduli space~\cite[page
347]{hartshorne:book} of stable $(n+1)$-pointed curves $C$ such that
the dual graph $\Gamma(C)$ has genus $g$, in the sense
of~(\ref{Good_Friday}).  The genus of $\Gamma(C)$ in fact equals the
arithmetic genus of the curve $C$, thus $\barcalM_{g,n+1}$ is the
coarse moduli space of stable curves of arithmetic genus~$g$ with
$n+1$ marked points. By a result of P.~Deligne, F.F.~Knudsen and
D.~Mumford~\cite{deligne-mumford:%
IHES69,knudsen-mumford:MS76,knudsen:MS83}, $\barcalM_{g,n+1}$ is a
projective variety.

Observe that, for a curve $C \in \barcalM_{0,n+1}$, the graph
$\Gamma(C)$ must necessarily be a tree and all components of $C$ must
be smooth of genus $0$, therefore $\barcalM_{0,n+1}$ coincides with
the moduli space $\barcalM_0(n)$\label{barcalMn1} of genus $0$ stable
curves with $n+1$ marked points that we discussed in
Section~\ref{sec1bis}.  Dual graphs of curves $C \in
\barcalM_{g,n+1}$ are stable labeled graphs belonging to $\MGr g{n+1}$.

The symmetric group $\Sigma_n^+$ acts on $\barcalM_{g,n+1}$ by
renumbering the marked points, therefore
\[
\barcalM := \modcoll {\barcalM},
\]
with $\barcalM \lzavorkaobyc g,n \rzavorkaobyc := \barcalM_{g,n+1}$,
is a modular $\Sigma$-module in the category of projective
varieties. Since there are no stable curves of genus $g$ with $n+1$
punctures if $2g+n-1 \leq 0$, $\barcalM$ is a {\em stable\/} modular
$\Sigma$-module.  Let us define the contraction along a stable graph
$\Gamma\in \Graphcat gn$
\begin{equation}
\label{1148}
\alpha_\Gamma:
\barcalM\((\Gamma\)) = \prod_{v\in\vert(\Gamma)}
\barcalM\((g(v),\edge(v)\)) \to\barcalM\lzavorkaobyc g,n \rzavorkaobyc
\end{equation}
by gluing the marked points of curves from
$\barcalM\((g(v),\edge(v)\))$, $v\in\vert(\Gamma)$, according to the graph
$\Gamma$.
To be more precise, let
\[
\prod_{v \in \vert(\Gamma)} C_v,\ \mbox { where }\
C_v \in \barcalM\((g(v),\edge(v)\)),
\]
be an element of $\barcalM\((\Gamma\))$.
Let $e$ be an edge of the graph $\Gamma$
connecting vertices $v_1$ and $v_2$, $e = \{y^e_{v_1},y^e_{v_2}\}$,
where $y^e_{v_i}$ is a marked point of the component $C_{v_i}$,
$i=1,2$, which is also the name of the corresponding
flag of the graph $\Gamma$.
The curve $\alpha_\Gamma(C)$ is then obtained by the identifications
$y^e_{v_1} = y^e_{v_2}$, introducing a nodal singularity, for all $e
\in \Edg(\Gamma)$.
The procedure is the
same as that described for the tree level in Section~\ref{sec1bis}.
As proved in~\cite[\S~6.2]{getzler-kapranov:CompM98},
the contraction maps~(\ref{1148}) define on the stable
modular $\Sigma$-module of
coarse moduli spaces
$\barcalM =\stablecoll{\barcalM}$  a modular operad structure in the
category of complex projective varieties.

Let us look more closely at the structure of the modular triple $\Modtriple$.
Given a (stable or unstable) modular $\Sigma$-module $\calE$,
there is, for each $g \geq 0$ and $n \geq -1$, a natural decomposition
\[
\Modtriple(\calE)(g,n) =
\Modtriple_0(\calE)(g,n) \oplus \Modtriple_1(\calE)(g,n)
\oplus \Modtriple_2(\calE)(g,n) \oplus \cdots,
\]
with $\Modtriple_k(\calE)(g,n)$ the subspace obtained by summing over
graphs $\Gamma$ with $\dim H_1(|\Gamma|) = k$, $k \geq 0$. In
particular, $\Modtriple_0(\calE)(g,n)$ is a summation over simply
connected graphs. It is not difficult to see that
$\Modtriple_0(\calE)$ is a {\em subtriple\/} of $\Modtriple(\calE)$.
This shows that modular operads are
$\Modtriple_0$-algebras with some additional operations (the
`contractions') that raise the genus and generate the higher components
$\Modtriple_k$, $k \geq 1$, of the modular triple $\Modtriple$.

There seems to be a belief expressed in the proof
of~\cite[Lemma~3.4]{getzler-kapranov:CompM98}
and also
in~\cite[Theorem~3.7]{getzler-kapranov:CompM98} that, in the stable case,
the triple $\Modtriple_0$ is equivalent to the non-unital cyclic
operad triple $\Psi_+$, but it is not
so. The triple $\Modtriple_0$ is {\em much bigger\/}, for
example, if $a \in \calE(1,0)$, then $\Modtriple_0(\calE)(2,-1)$
contains a non-operadic element
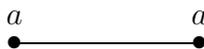
\begin{center}
{
\unitlength=1.000000pt
\begin{picture}(70.00,12.00)(0.00,0.00)
\put(70.00,7.00){\makebox(0.00,0.00)[b]{$a$}}
\put(0.00,7.00){\makebox(0.00,0.00)[b]{$a$}}
\put(70.00,0.00){\makebox(0.00,0.00){$\bullet$}}
\put(0.00,0.00){\makebox(0.00,0.00){$\bullet$}}
\put(0.00,0.00){\line(1,0){70.00}}
\end{picture}}
\end{center}
which can be also written,
using~(\ref{Prave_jsem_mluvil_s_Jitkou}), as $a \circij 00 a$.
The corresponding part
$\Psi_+(\calE)(-1)$ of the cyclic triple is empty.  In the Grothendieck-Knudsen
modular operad $\barcalM$, an element of the above type is realized by two
tori meeting at a nodal point.

On the other hand, the triple $\Modtriple_0$ restricted to the
subcategory of stable modular $\Sigma$-modules $\calE$ such that
$\calE(g,n)=0$ for $g > 0$ indeed coincides with the non-unital cyclic
operad triple $\Psi_+$, as was in fact proved
in~\cite[page~81]{getzler-kapranov:CompM98}. Therefore, given a
modular operad $\calA = \stablecoll{\calA}$, there is an induced
non-unital cyclic operad structure on the cyclic collection
$\calA^\flat := \{\calA(0,n)\}_{n \geq 2}$. We will call $\calA^\flat$
the {\em associated cyclic operad\/}. For example, the
cyclic operad associated to the Grothendieck-Knudsen modular operad
$\barcalM$ equals its genus zero part $\barcalM_0$.

A {\em biased\/} definition of modular operads can be
found in~\cite[Definition~II.5.35]{markl-shnider-stasheff:book}. It is
formulated in terms of operations
\[
\left\{
\circij ij : \calA(g,m) \ot \calA(h,n) \to \calA(g+h,m+n);\
0 \leq i \leq m,\ 0 \leq j \leq n,\ g,h \geq 0
\right\}
\]
together with contractions
\[
\left\{
\xi_{\{i,j\}} : \calA(g,m) \to \calA(g+1,m-2);\
m \geq 1,\ g \geq 0
\right\}
\]
that
generalize~(\ref{prave_jsem_mluvil_s_Jitkou})
and~(\ref{prave_jsem_mluvil_s_Jitkou.}).

\begin{example}
\label{Kacirek_Donald}
{\rm\
Let $V = (V,B)$ be a vector space with a symmetric
inner product $B : V\ot V \to \bfk$.
Denote, for each $g \geq 0$ and $n \geq -1$,
\[
\End_V(g,n) := \otexp V{(n+1)}.
\]
It is clear from definition~(\ref{Cosmo}) that,
for any labeled graph $\Gamma\in \Graphcat gn$,
$\End_V \((\Gamma\)) = \otexp V{\Flag(\Gamma)}$.

Let $\otexp B{\Edg(\Gamma)} : \otexp V{\Flag(\Gamma)} \to \otexp
V{\Leg(\Gamma)}$ be the multilinear form which contracts the factors
of $\otexp V{\Flag(\Gamma)}$ corresponding to the flags which are
paired up as edges of $\Gamma$. Then we define
$\alpha_\Gamma: \End_V \((\Gamma\)) \to \End_V(g,n)$ to be
the map
\[
\alpha_\Gamma :
\End_V\((\Gamma\)) = \otexp V{\Flag(\Gamma)}
\ \stackrel{\otexp B{\Edg(\Gamma)}}{\vlra}\
\otexp
V{\Leg(\Gamma)} \stackrel{\otexp V{\ell}}{\longrightarrow}
\otexp V{(n+1)} = \End_V(g,n),
\]
where $\ell : \Leg(\Gamma) \to \n n$ is the labeling of the legs of $\Gamma$.
It is easy to show that the compositions $\{\alpha_\Gamma;\
\Gamma \in \Graphcat gn \}$
define on $\End_V$ the structure of an un-stable unital modular
operad, see~\cite[\S~2.25]{getzler-kapranov:CompM98}.
}
\end{example}

An {\em algebra\/} over a modular operad $\calA$ is a vector space $V$ with an
inner product $B$, together with a morphism $\rho : \calA \to \End_V$
of modular operads.
Several important structures are algebras over modular
operads. For example, an
algebra over the homology $H_*(\barcalM)$ of the Grothendieck-Knudsen
modular operad is the
same as a cohomological field theory in the sense
of~\cite{kontsevich-manin:CMP94}.  Other physically relevant
algebras over modular operads can be found
in~\cite{getzler-kapranov:CompM98,markl:la,markl-shnider-stasheff:book}.
Relations between modular
operads, chord diagrams and
Vassiliev invariants are studied in~\cite{hinich-vaintrob}.

\section{PROPs}
\label{sec7}

Operads are devices invented to describe structures consisting of
operations with several inputs and {\em one\/} output. There are,
however, important structures with operations having
several inputs and {\em several\/} outputs. Let us recall the most
prominent one:

\begin{example}
\label{za_chvili_Kolin}
{\rm
A (associative) {\em bialgebra\/} is a $\bfk$-module $V$ with
a {\em multiplication\/} $\mu : V \ot V \to V$ and a {\em
comultiplication\/} (also called a {\em diagonal\/}) $\Delta : V \to
V\ot V$. The multiplication is associative:
\[
\mu(\mu \ot \id_V) = \mu(\id_V \ot \mu),
\]
the comultiplication is coassociative:
\[
(\Delta \ot \id_V)\Delta = (\id_V \ot \Delta)\Delta
\]
and the usual compatibility between $\mu$ and $\Delta$ is
assumed:
\begin{equation}
\label{och_ty_Jitky}
\Delta(u \cdot v) = \Delta(u) \cdot \Delta(v)\ \mbox { for }\  u,v \in V,
\end{equation}
where $u \cdot v :=\mu(u,v)$ and the dot $\ \cdot\ $ in the right hand
side denotes the multiplication induced on $V \ot V$ by $\mu$. Loosely
speaking, bialgebras are Hopf algebras without unit, counit and antipode.
}
\end{example}

\PROP{s} (an abbreviation of {\bf pro}duct and {\bf p}ermutation
category) describe structures as in
Example~\ref{za_chvili_Kolin}.
Although \PROP{s} are more general than operads, they appeared
much sooner, in a~1965 Mac~Lane's paper~\cite{maclane:BAMS65}.
This might be explained by the fact that the definition
of \PROP{s} is more compact than that of operads -- compare
Definition~\ref{Jituska} below with Definition~\ref{a} in Section~\ref{sec1}.
\PROP{s} then entered the `renaissance of
operads' in 1996 via~\cite{markl:JPAA96}.

Definition~\ref{Jituska} uses the notion of a symmetric strict
monoidal category which  we consider so basic and commonly known that we
will not recall it, standard citations
are~\cite{maclane:BAMS65,maclane:RiceUniv.Studies63}, see
also~\cite[\S~II.1.1]{markl-shnider-stasheff:book}.
An~example is the category $\Modk$
of $\bfk$-modules, with the monoidal
product $\odot$ given by the tensor product $\otimes = \ot_\bfk$, the
symmetry $S_{U,V} : U \ot V \to V\ot U$ defined as $S_{U,V}(u,v) := v
\ot u$ for $u \in U$ and $v \in V$, and the unit
$\bbjedna$ the ground ring $\bfk$.

\begin{definition}
\label{Jituska}
A (\/$\bfk$-linear) {\em \PROP\/} (called a {\em theory\/}
in~\cite{markl:JPAA96}) is a symmetric strict monoidal
category $\sfP = (\sfP,\odot,S,\bbjedna)$  enriched over\/ $\Modk$  such that
\begin{itemize}
\item[(i)]
the objects are indexed by (or
identified with) the set $\mathbb N = \{0,1,2,\ldots\}$ of natural numbers, and
\item[(ii)]
the product satisfies
$m \odot n=m+n$, for any $m,n \in {\mathbb N}
=\mbox{\rm Ob}(\sfP)$ (hence the unit $\bbjedna$ equals~$0$).
\end{itemize}
\end{definition}

Recall that the $\Modk$-enrichment in the above definition means that
each hom-set $\Mor_\sfP(m,n)$ is a $\bfk$-module and the operations of
the monoidal category $\sfP$ (the composition $\circ$, the product
$\odot$ and the symmetry $S$) are compatible with this $\bfk$-linear
structure.

For a \PROP\ $\sfP$ denote $\sfP(m,n) := \Mor_\sfP(m,n)$. The symmetry
$S$ induces, via the canonical identifications $m \cong
\bbjedna^{\odot m}$ and $n \cong \bbjedna^{\odot n}$, on each
$\sfP(m,n)$ a structure of $(\Sigma_m,\Sigma_n)$-bimodule (left
$\Sigma_m$- right $\Sigma_n$-module such that the left action
commutes with the right one).  Therefore a \PROP\ is a collection
$\sfP = \{\sfP(m,n)\}_{m,n \geq 0}$ of
$(\Sigma_m,\Sigma_n)$-bimodules,
together with two types of compositions, {\em horizontal\/}
\[
\ot : \sfP(m_1,n_1) \ot \cdots \ot \sfP(m_s,n_s) \to
\sfP(m_1+\cdots +m_s,n_1+\cdots +n_s),
\]
induced, for all $m_1,\ldots,m_s,n_1,\ldots,n_s \geq 0$, by the monoidal
product $\odot$ of $\sfP$, and {\em vertical\/}
\[
\circ : \sfP(m,n) \ot \sfP(n,k) \to \sfP(m,k),
\]
given, for all $m,n,k \geq 0$, by the categorial composition.
The monoidal unit is an element $e := \bbjedna \in \sfP(1,1)$. In
Definition~\ref{Jituska}, $\Modk$ can be replaced by an arbitrary
symmetric strict monoidal category.

Let $\sfP = \bicol \sfP 0$ and $\sfQ = \bicol\sfQ0$ be two \PROP{s}. A
{\em homomorphism\/} $f : \sfP \to \sfQ$ is a sequence $f =
\{f(m,n) : \sfP(m,n) \to \sfQ(m,n)\}_{m,n\geq 0}$ of
bi-equivariant maps which commute with both the vertical and horizontal
compositions. An {\em ideal\/} in a \PROP\ $\sfP$ is a
system $\sfI = \bicol\sfI0$ of left $\Sigma_m$- right
$\Sigma_n$-invariant subspaces $\sfI(m,n) \subset \sfP(m,n)$ which is closed,
in the obvious sense, under both the vertical and horizontal
compositions.
Kernels, images,~\&c., of homomorphisms between \PROP{s},
as well as quotients of \PROP{s}
by \PROP{ic} ideals, are defined componentwise,
see~\cite{markl:JPAA96,vallette:thesis,vallette:extract,vallette:CMR04} 
for details.

\begin{example}
{\rm\
The {\em endomorphism\/}  \PROP\ of a $\bfk$-module
$V$ is the system
\[
\End_V = \bicol{\End_V}0
\]
with $\End_V(m,n)$ the space of linear maps $\Lin(\otexp
Vn,\otexp Vm)$ with $n$ `inputs' and $m$ `outputs,' $e \in
\End_V(1,1)$ the identity map, horizontal composition given by the
tensor product of linear maps, and vertical composition by the
ordinary composition of linear maps.
}
\end{example}

Also algebras over \PROP{s} can be introduced in a very
concise way:

\begin{definition}
A {\em $\sfP$-algebra\/} is a strict symmetric monoidal functor
$\lambda: \sfP \to \Modk$ of enriched monoidal categories. The value
$\lambda(1)$ is the {\em underlying space\/} of the algebra $\rho$.
\end{definition}

It is easy to see that a $\sfP$-algebra is the same as a \PROP\
homomorphism $\rho : \sfP \to \End_V$. As in Proposition~\ref{pro:alg},
a $\sfP$-algebra is determined by a system
\[
\alpha: \sfP(m,n) \otimes \otexp Vn \to \otexp Vm,\ m,n, \geq 0,
\]
of linear maps satisfying appropriate axioms.

As before, the first step in formulating an unbiased definition of
\PROP{s} is to specify their underlying category.
A {\em $\Sigma$-bimodule\/} is a system $E= \bicol E0$ such that each
$E(m,n)$ is a left $\bfk[\Sigma_m]$- right $\bfk[\Sigma_n]$-bimodule. Let
$\sigmabimod$ denote the category of $\Sigma$-bimodules.
For $E \in \sigmabimod$ and finite sets $Y,X$ with $m$ resp.~$n$
elements put
\[
E(Y,X) := \Bij(Y,[m]) \times_{\Sigma_m} E(m,n) \times_{\Sigma_n}
\Bij([n],X),\ m,n \geq 0,
\]
where $\Bij(-,-)$ is the same as in~(\ref{maminka_se_zase_nacpala}).
Pasting schemes for \PROP{s} are {\em directed $(m,n)$-graphs\/},
by which we mean finite, not necessary connected, graphs in the sense of
Definition~\ref{graph-def} such that
\begin{itemize}
\item[(i)]
each edge is equipped with a direction
\item[(ii)]
there are no directed cycles and
\item[(iii)] the set of legs is divided into the set of inputs labeled by
$\{1, \dots, n\}$ and the set of outputs labeled by $\{1, \dots, m\}$.
\end{itemize}
An example of a directed graph is given in
Figure~\ref{nemohu_spat}.
We denote by $\Gr mn$ the category of directed $(m,n)$-graphs and
their isomorphisms. The direction of edges determines at
each vertex $v \in \vert(G)$ of a directed graph $G$ a disjoint decomposition
\[
\edge(v) = \In(v) \sqcup \Out(v)
\]
of the set of edges adjacent to $v$ into the set $\In(v)$ of incoming
edges and the set $\Out(v)$ of outgoing edges.  The pair
$(\#(\Out(v)),\#(\In(v))) \in {\mathbb N} \times {\mathbb N}$
is called the {\em
biarity\/} of $v$.
\begin{figure}
\begin{center}
{
\unitlength=1.000000pt
\begin{picture}(250.00,110.00)(-80.00,0.00)
\thicklines
\put(60.00,0.00){\makebox(0.00,0.00){$1$}}
\put(0.00,0.00){\makebox(0.00,0.00){$3$}}
\put(30.00,0.00){\makebox(0.00,0.00){$2$}}
\put(80.00,100.00){\makebox(0.00,0.00){$2$}}
\put(60.00,100.00){\makebox(0.00,0.00){$3$}}
\put(40.00,100.00){\makebox(0.00,0.00){$1$}}
\put(20.00,100.00){\makebox(0.00,0.00){$4$}}
\put(80.00,30.00){\makebox(0.00,0.00){$\bullet$}}
\put(60.00,70.00){\makebox(0.00,0.00){$\bullet$}}
\put(60.00,50.00){\makebox(0.00,0.00){$\bullet$}}
\put(40.00,50.00){\makebox(0.00,0.00){$\bullet$}}
\put(40.00,70.00){\makebox(0.00,0.00){$\bullet$}}
\put(20.00,50.00){\makebox(0.00,0.00){$\bullet$}}
\put(20.00,30.00){\makebox(0.00,0.00){$\bullet$}}
\put(40.00,50.00){\vector(0,1){20.00}}
\put(60.00,10.00){\vector(0,1){40.00}}
\put(30.00,10.00){\vector(-1,2){10.00}}
\put(0.00,10.00){\vector(1,1){20.00}}
\put(20.00,30.00){\vector(2,1){40.00}}
\put(20.00,30.00){\vector(0,1){20.00}}
\put(20.00,50.00){\vector(1,1){20.00}}
\put(80.00,30.00){\vector(0,1){60.00}}
\put(60.00,50.00){\vector(0,1){20.00}}
\put(60.00,50.00){\vector(-1,1){20.00}}
\put(80.00,30.00){\vector(-1,1){20.00}}
\put(40.00,70.00){\vector(-1,1){20.00}}
\put(40.00,70.00){\vector(0,1){20.00}}
\put(40.00,70.00){\vector(1,1){20.00}}
\end{picture}}
\end{center}
\caption{\label{nemohu_spat}%
A directed graph from $\Gr 43$.}
\end{figure}
To incorporate the unit, we need to extend the category $\Gr mn$, for
$m=n$, into the category $\UGr mn$ by allowing the
exceptional graph
\[
\uparrow \ \uparrow \ \uparrow \ \cdots \uparrow\ \in \UGr nn,\ n \geq 1,
\]
with $n$ inputs, $n$ outputs and no vertices.
For a graph $G \in \UGr mn$ and a $\Sigma$-bimodule
$E$, let
\[
E(G) := \bigotimes_{v \in \vert(G)} E(\Out(v),\In(v)).
\]
and
\begin{equation}
\label{jsem_zvedav_jestli_se_ta_druha_Jitka_ozve}
\freePROP(E)(m,n) := \colim{{G \in \UGr mn}}{{E}}(G),\
m,n \geq 0.
\end{equation}
The $\Sigma$-bimodule $\freePROP(E)$
is a \PROP, with the vertical composition given by the disjoint union
of graphs, the horizontal composition by grafting the legs, and the
unit the exceptional graph $\uparrow \hskip .2em \in \freePROP(E)(1,1)$.
The following proposition follows from~\cite{mv}
and~\cite{vallette:thesis,vallette:extract,vallette:CMR04}:

\begin{proposition}
\label{a_pratele_si_znechutim}
The \PROP\ $\freePROP(E)$ is the {\em free \PROP\/} generated by the
$\Sigma$-bimodule $E$.
\end{proposition}

As in the previous sections,
(\ref{jsem_zvedav_jestli_se_ta_druha_Jitka_ozve})
defines a triple $\freePROP : \sigmabimod \to
\sigmabimod$ with the triple multiplication of erasing the braces.
According to general principles~\cite{e-m},
Proposition~\ref{a_pratele_si_znechutim} is almost equivalent to

\begin{proposition}
\PROP{s} are algebras over the triple $\freePROP$.
\end{proposition}

One may obviously consider {\em non-unital \PROP{s}\/} defined as
algebras over the triple
\[
\Psi_{\tt P}(E)(m,n) := \colim{{G \in \Gr mn}}{{E}}(G),\
m,n \geq 0,
\]
and develop a theory parallel to the theory of non-unital operads
reviewed in Section~\ref{sec1bis}.

\begin{example}
\label{poruchy_spanku}
{\rm
We will use
the graphical language explained in
Example~\ref{Jitka_ted_moc_nepise}.
Let $\Gamma(\gen12,\gen21)$ be the free \PROP\ generated by one
operation $\gen12$ of biarity $(1,2)$ and one operation $\gen21$ of
biarity $(2,1)$. As we noticed
already in~\cite{markl:ws93,markl:JPAA96}, the \PROP\ $\sfB$ describing
bialgebras equals
\[
\sfB = \Gamma(\gen12,\gen21)/{\sf I}_{\sf B},
\]
where ${\sf I}_{\sf B}$ is the \PROP{ic} ideal generated by
\begin{equation}
\label{Jo_ty_Jitky}
\ZbbZb - \bZbbZ,\
\ZvvZv - \vZvvZ\  \mbox { and }\
\dvojiteypsilon - \motylek \hskip .2em.
\end{equation}
In the above display we denoted
\begin{eqnarray*}
&\ZbbZb := \gen12\circ (\gen12 \ot e),\
\bZbbZ := \gen12\circ (e \ot \gen12),\
\ZvvZv := (\gen21 \ot e)\circ \gen21,\
\vZvvZ := (e \ot \gen21)\circ \gen21,&
\\
&\dvojiteypsilon  := \gen21 \circ \gen12
\mbox {\hskip .5em and \hskip .6em}
\motylek := (\gen12 \ot \gen12)\circ \kappa \circ
(\gen21 \ot \gen21),&
\end{eqnarray*}
where $\kappa \in \Sigma_4$ is the permutation
\begin{equation}
\label{novy_tabak}
\kappa :=
\left(
\begin{array}{cccc}
1 & 2 & 3 & 4
\\
1 & 3 & 2 & 4
\end{array}
\right)
 = \hskip 2mm
{
\unitlength=1.000000pt
\begin{picture}(30.00,30.00)(0.00,12.00)
\put(30.00,0.00){\makebox(0.00,0.00){\scriptsize$\bullet$}}
\put(30.00,30.00){\makebox(0.00,0.00){\scriptsize$\bullet$}}
\put(20.00,0.00){\makebox(0.00,0.00){\scriptsize$\bullet$}}
\put(20.00,30.00){\makebox(0.00,0.00){\scriptsize$\bullet$}}
\put(10.00,0.00){\makebox(0.00,0.00){\scriptsize$\bullet$}}
\put(10.00,30.00){\makebox(0.00,0.00){\scriptsize$\bullet$}}
\put(0.00,0.00){\makebox(0.00,0.00){\scriptsize$\bullet$}}
\put(0.00,30.00){\makebox(0.00,0.00){\scriptsize$\bullet$}}
\put(30.00,20.00){\line(0,-1){10.00}}
\put(20.00,20.00){\line(-1,-1){10.00}}
\put(10.00,20.00){\line(1,-1){10.00}}
\put(0.00,20.00){\line(0,-1){10.00}}
\end{picture}}\hskip .3em .
\end{equation}
\vglue .6em

The above description of $\sfB$ is `tautological,' but
B.~Enriquez and P.~Etingof
found in~\cite[Proposition~6.2]{enriquez-etingof:03} the following
basis of the $\bfk$-linear space $\sfB(m,n)$ for arbitrary $m,n \geq 1$.
 Let $\gen 12 \in \sfB(1,2)$ be the equivalence class, in
$\sfB = \Gamma(\gen12,\gen21)/{\sf I}_{\sf B}$, of the generator $\gen
12 \in \Gamma(\gen 21,\gen 12)(1,2)$ (we use the same symbol both for
a generator and its equivalence class). Define ${\gen 12}^{[1]} := e
\in \sfB(1,1)$ and, for $a \geq 2$, let
\[
{\gen 12}^{[a]} :=
   \gen 12\circ (\gen 12 \ot e)\circ (\gen 12 \ot \otexp {e}2) \circ\cdots
   \circ(\gen 12 \ot \otexp {e}{(a-2)}) \in \sfB(1,a).
\]
Let ${\gen 21}_{[b]} \in \sfB(b,1)$ has the obvious similar
meaning. The elements
\begin{equation}
\label{je_to_hruza_se_mnou}
({\gen 12}^{[a_1]} \ot \cdots \ot {\gen 12}^{[a_m]})
\circ \sigma \circ
({\gen 21}_{[b^1]} \ot \cdots \ot {\gen 21}_{[b^n]}),
\end{equation}
where $\sigma \in \Sigma_N$ for some $N \geq 1$, and $a_1 + \cdots +
a_m = b^1 + \cdots + b^m = N$, form a $\bfk$-linear basis
of~$\sfB(m,n)$.
This result can also be found
in~\cite{kontsevich-soibelman:hopf-algebras}. See
also~\cite{lack,pirashvili} for the bialgebra \PROP\ viewed from a
different perspective.
}
\end{example}

\begin{example}
{\rm Each operad $\calP$ generates a unique \PROP\ $\sfP$ such
that $\sfP(1,n) = \calP(n)$ for each $n \geq 0$.
The components of such a \PROP\ are given by
\[
\sfP(m,n) = \bigoplus_{r_1 + \cdots + r_k = n}
[\calP(1,r_1) \otimes \cdots \otimes \calP(1,r_k)]
\times_{\Sigma_{r_1} \times \cdots \times \Sigma_{r_k}} \Sigma_n,
\]
for each $m,n \geq 0$. The (topological) \PROP{s} considered in~\cite{boardman-vogt:73}
are all of this type. On the other hand,
Example~\ref{poruchy_spanku} shows that not each
\PROP\ is of this form. A
\PROP\ $\sfP$ is generated by an operad if and only if it has a
presentation $\sfP = \freePROP(E)/(R)$, where $E$ is a
$\Sigma$-bimodule such that $E(m,n) = 0$ for $m \not= 1$ and $R$ is
generated by elements in $\freePROP(E)(1,n)$, $n \geq 0$.
}
\end{example}

\section{Properads, dioperads and $\frac 12$PROPs}
\label{sec8}

As we saw in
Proposition~\ref{to_by_mne_zajimalo_jestli_jsem_to_s_ni_pokazil}, under
some mild assumptions, the components of free operads are
finite-dimensional.
In contrast, \PROP{s} are huge objects.
For example, the
component $\freePROP(\gen 12, \gen 21)(m,n)$ of the
free \PROP\ $\freePROP(\gen 12, \gen 21)$ used in the definition of the
bialgebra \PROP\ $\sfB$ in Example~\ref{poruchy_spanku} is infinite-dimensional
for each $m,n \geq 1$, and also the components of the bialgebra \PROP\
$\sfB$ itself are infinite-dimensional, as follows from
the fact that the
Enriquez-Etingof basis~(\ref{je_to_hruza_se_mnou}) of $\sfB(m,n)$
has, for $m,n \geq 1$, infinitely many elements.

To handle this combinatorial explosion of \PROP{s} combined with lack of
suitable filtrations, smaller versions of \PROP{s} were invented.
Let us begin with the simplest modification which
we use as an example which explains the general scheme of
modifying \PROP{s}.
Denote $\UGrc mn$ the full subcategory of $\UGr
mn$ consisting of {\em connected\/} graphs and consider the triple
defined by
\begin{equation}
\freeproperad(E)(m,n) := \colim{{G \in \UGrc mn}}{{E}}(G),\
m,n \geq 0,
\end{equation}
for $E \in \sigmabimod$. The following notion was introduced by
B.~Vallette~\cite{vallette:thesis,vallette:extract,vallette:CMR04}.

\begin{definition}
{\em Properads\/} are algebras over the triple $\freeproperad :
\sigmabimod \to \sigmabimod$.
\end{definition}

A properad is therefore a $\Sigma$-bimodule with operations that
determine coherent contractions along connected graphs.
A biased definition of properads is given
in~\cite{vallette:thesis,vallette:extract,vallette:CMR04}. Since $\Gamma_c$
is a subtriple of $\freePROP$, each \PROP\ is automatically also a
properad. Therefore one may speak about the {\em endomorphism
properad\/} $\End_V$ and define
{\em algebras\/} over a properad $P$ as
properad homomorphisms $\rho : P \to \End_V$. Algebras over other versions
of \PROP{s} recalled below can be defined in a similar way.

\begin{example}
{\rm\
Associative bialgebras reviewed in Example~\ref{poruchy_spanku} are algebras
over the properad~$B$ defined
(tautologically) as the quotient of the free properad
$\freeproperad(\gen 12,\gen 21)$ by the properadic ideal generated by
the elements listed
in~(\ref{Jo_ty_Jitky}). We leave as an exercise to describe the
sub-basis of~(\ref{je_to_hruza_se_mnou}) that span $B(m,n)$, $m,n \geq
1$.

The following slightly artifical structure exists
over \PROP{s} but not over properads.
It consists of a `multiplication' $\mu =\gen 12 :
V \ot V \to V$, a `comultiplication' $\Delta = \gen 21  :V \to V \ot
V$ and a linear map $f = \ef : V \to V$ satisfying
$\Delta \circ \mu = f \ot f$ or, diagrammatically
\[
\dvojiteypsilon = \ef\ \ef.
\]
This structure cannot be a properad algebra because the graph on
the right hand side of the above display is not connected.
}
\end{example}

Properads are still huge objects.
The first really small version of \PROP{s} were dioperads
introduced in~2003 by W.L.~Gan~\cite{gan}.
As a motivation for his definition, consider the following:

\begin{example}
\label{80}
{\rm
A {\em Lie bialgebra\/} is a vector space $V$ with a Lie
algebra structure $[-,-] = \gen 12 : V \ot V \to V$
and a Lie diagonal $\delta = \gen 21 : V
\to V \ot V$.
We assume that $[-,-]$ and $\delta$ are related by
\begin{equation}
\label{Jitka_Spetkova_napsala}
\delta[a,b] = \sum \left( [a_{(1)},b] \ot a_{(2)} +  [a,b_{(1)}] \ot b_{(2)}
               + a_{(1)} \ot [a_{(2)},b] + b_{(1)} \ot [a,b_{(2)}]  \right)
\end{equation}
for any $a,b \in V$, with the Sweedler notation
$\delta a = \sum a_{(1)} \ot a_{(2)}$ and $\delta b = \sum b_{(1)}
\ot b_{(2)}$.

Lie bialgebras are governed by the \PROP\
$\sfLieB = \Gamma(\gen12,\gen21)/{\sf I}_{\sfLieB}$,
where $\gen 12$ and $\gen 21$ are now {\em antisymmetric\/}
and  ${\sfI}_{\sfLieB}$ denotes the ideal generated by
\begin{equation}
\label{moralni_hangover}
\rule{0pt}{2em}
\Jac 123 + \Jac 231 + \Jac 312 \hskip 1mm , \
\coJac 123 +  \coJac 231 + \coJac 312 \hskip 1mm   \mbox { and }
\hskip 2mm
\dvojiteypsilonvetsi1212 -
\levaplastevvetsi1212  -
\pravaplastevvetsi1212 +
\levaplastevvetsi1221 + \pravaplastevvetsi1221 \hskip 1mm,
\end{equation}

\vskip 2mm
\noindent
with labels indicating the corresponding
permutations of the inputs and outputs.
}
\end{example}

We observe that all graphs in~(\ref{moralni_hangover}) are not only
connected as demanded for properads, but also
simply-connected.
This suggests considering
the full subcategory $\UGrd
mn$ of $\UGr mn$ consisting of {\em
connected simply-connected\/} graphs and the related triple
\begin{equation}
\label{Jitka_uz_zase_premysli}
\freedioperad(E)(m,n) := \colim{{G \in {\UGrd mn}}}{E}(G),\
m,n \geq 0.
\end{equation}

\begin{definition}
{\em Dioperads\/} are algebras over the triple $\freedioperad :
\sigmabimod \to \sigmabimod$.
\end{definition}

A biased definition of dioperads can be found in~\cite{gan}. As
observed by T.~Leinster, dioperads are more or less equivalent to
polycategories, in the sense of~\cite{szabo}, with one object.
Lie bialgebras reviewed in Example~\ref{80} are algebras over a
dioperad. Another important class of dioperad algebras is recalled in:

\begin{example}
{\rm
An {\em infinitesimal bialgebra\/}~\cite{joni-rota}
(called in~\cite[Example~11.7]{fox-markl:ContM97} {\em a mock
bialgebra\/}) is a vector space $V$ with an associative
multiplication $\cdot : V \ot V \to V$ and a coassociative
comultiplication $\Delta : V \to V \ot V$
such that
\[
\Delta(a \cdot b) = \sum\left(
a_{(1)} \ot a_{(2)}\cdot b + a \cdot b_{(1)} \ot b_{(2)}
\right)
\]
for any $a,b \in V$. It is easy to see that the axioms of infinitesimal
bialgebras are encoded by the following simply connected graphs:
\[
\ZbbZb - \bZbbZ,\
\ZvvZv - \vZvvZ\  \mbox { and }\
\dvojiteypsilon - \levaplastev  - \pravaplastev\hskip 1mm.
\]
}
\end{example}

Observe that
associative bialgebras recalled in
Example~\ref{za_chvili_Kolin} cannot be defined over dioperads,
because the rightmost graph in~(\ref{Jo_ty_Jitky})
is not simply connected.
The following
proposition, which should be compared to
Proposition~\ref{to_by_mne_zajimalo_jestli_jsem_to_s_ni_pokazil}, shows that
dioperads are of the same size as operads.

\begin{proposition}
Let $E = \bicol E0$ be a $\Sigma$-bimodule such that
\begin{equation}
\label{mam_ji_napsat?}
E(m,n) =  0 \mbox { for } m+n \leq 2
\end{equation}
and that $E(m,n)$ is  finite-dimensional for all
remaining $m,n$. Then the components
$\freedioperad(E)(m,n)$ of the free dioperad $\freedioperad(E)$ are
finite-dimensional, for all $m,n \geq 0$.
\end{proposition}

The proof, similar to the proof of
Proposition~\ref{to_by_mne_zajimalo_jestli_jsem_to_s_ni_pokazil},
is based on the observation that the
assumption~(\ref{mam_ji_napsat?}) reduces the
colimit~(\ref{Jitka_uz_zase_premysli}) to a summation over
reduced trees (trees whose all vertices have at least three
adjacent edges).

An important problem arising in connection with deformation
quantization is to find a reasonably small, explicit cofibrant
resolution of the bialgebra \PROP~$\sfB$.
Here by a resolution we mean a differential graded \PROP\ $\sfR$
together with a homomorphism $\beta : \sfR \to \sfB$ inducing
a homology isomorphism.
Cofibrant in this context means that $\sfR$ is of the form
$(\freePROP(E),\partial)$, where the generating $\Sigma$-bimodule
$E$ decomposes as $E = \bigoplus_{n \geq 0}E_n$ and the differential
decreases the filtration, that is
\[
\partial(E_n) \subset \freePROP(E)_{< n},
\mbox { for each $n \geq 0$,}
\]
where $\freePROP(E)_{< n}$ denotes the sub-\PROP\ of $\Gamma(E)$ generated by
$\bigoplus_{j < n} E_j$.
This notion  is an \PROP{ic} analog of the Koszul-Sullivan
algebra in rational homotopy theory~\cite{halperin:lect}.
Several papers devoted to
finding~$\sfR$ appeared
recently~\cite{kontsevich-soibelman:hopf-algebras,%
umble-saneblidze:KK,saneblidze-umble:bialgebras,%
saneblidze-umble:matrons,shoikhet:2/3,shoiket:explicit,shoikhet}.
The approach of~\cite{markl:ba} is based on the observation
that $\sfB$ is a deformation, in the sense explained below,
of the \PROP\  describing structures recalled in the following:

\begin{definition}
A {\em half-bialgebra\/} or simply a {\em $\frac12$bialgebra\/} is a
vector space $V$ with an associative multiplication $\mu: V \ot V \to
V$ and a  coassociative comultiplication
$\Delta: V \to V \ot V$
that satisfy
\begin{equation}
\label{pf}
\Delta(u \cdot v) = 0,\ \mbox { for each } u,v \in V.
\end{equation}
\end{definition}

We chose this strange name because~(\ref{pf}) is indeed
one half of the compatibility relation~(\ref{och_ty_Jitky}) of
associative bialgebras.
$\frac12$bialgebras are algebras over the \PROP\
\[
\textstyle\frac12 \sfB: = \Gamma(\gen 12)/(\ZbbZb = \bZbbZ, \ZvvZv =
\vZvvZ,\ \dvojiteypsilon = 0).
\]
Now define, for a formal variable $t$, $\sfB_t$ to be the quotient of the free
\PROP\ $\Gamma(\gen 12,\gen 21)$ by the ideal generated by
\[
\ZbbZb = \bZbbZ,\ \ZvvZv =
\vZvvZ,\ \dvojiteypsilon = t\cdot \motylek.
\]

Thus $\sfB_t$ is a one-parametric family of \PROP{s} with the property that
$\sfB_0 = \frac12\sfB$. At a generic~$t$, $\sfB_t$~is
isomorphic to the bialgebra \PROP\ $\sfB$.
In other words, the \PROP\
for bialgebras is a deformation of the \PROP\ for $\frac12$bialgebras.
According to general principles of homological perturbation
theory~\cite{gugenheim-stasheff:BSMB86},
one may try to construct the resolution $\sfR$ as a
perturbation of a cofibrant resolution $\frac12\sfR$  of
the \PROP\ $\frac12\sfB$. Since $\frac12\sfB$ is simpler that $\sfB$, one may
expect that resolving $\frac12\sfB$ would be a simpler task than
resolving $\sfB$.

For instance, one may realize that $\frac12$bialgebras are algebras
over a dioperad $\frac12{\tt B}$,
use~\cite{gan} to construct a
resolution $\frac12{\tt R}$ of the dioperad $\frac12{\tt B}$, and then take
$\frac12\sfR$ to be the \PROP\ generated by
$\frac12{\tt R}$. More precisely,
one denotes
\begin{equation}
\label{jsem_uplne_hotovej}
F_1 : \catDiop \to \catPROP
\end{equation}
the left adjoint to the forgetful functor
$\catPROP \stackrel{\Box_1}\longrightarrow \catDiop$ and defines
$\frac12\sfR := F_1(\frac12{\tt R})$.

The problem is that we do not know whether the functor $F_1$ is
exact, so it is not clear if $\frac12\sfR$ constructed in
this way is really a resolution of $\frac12\sfB$.
To get around this subtlety, M.~Kontsevich observed that
$\frac12$bialgebras live over a version of \PROP{s}
which is smaller than dioperads. It can be defined as follows.

\begin{figure}
\unitlength 3mm
\begin{picture}(28,12)(7,15)
\thicklines
\put(14,19){\vector(0,1){4}}
\put(14,23){\vector(1,1){3}}
\put(14,23){\vector(1,2){2}}
\put(14,23){\vector(0,1){4.5}}
\put(14,23){\vector(-1,2){2}}
\put(14,23){\vector(-1,1){3}}
\put(11,21){\vector(3,2){3}}
\put(17,21){\vector(-3,2){3}}
\put(11,16){\vector(1,1){3}}
\put(17,16){\vector(-1,1){3}}
\put(14,16){\vector(0,1){3}}
\put(25,19){\vector(0,1){4}}
\put(25,23){\vector(1,1){3}}
\put(25,23){\vector(-1,2){2}}
\put(25,23){\vector(0,1){4.5}}
\put(25,23){\vector(1,2){2}}
\put(25,23){\vector(-1,1){3}}
\put(25,19){\vector(1,2){2}}
\put(25,19){\vector(-1,2){2}}
\put(22,16){\vector(1,1){3}}
\put(28,16){\vector(-1,1){3}}
\put(25,16){\vector(0,1){3}}
\put(25,19){\vector(1,1){3}}
\put(25,19){\vector(-1,1){3}}
\put(23,15){\vector(1,2){2}}
\put(27,15){\vector(-1,2){2}}
\put(-11,0){
\put(23,15){\vector(1,2){2}}
\put(27,15){\vector(-1,2){2}}
}
\put(-11,4){
\put(23,15){\vector(1,2){2}}
\put(27,15){\vector(-1,2){2}}
}
\put(25,14.5){\line(0,1){2}}
\put(14,14.5){\line(0,1){2}}
\put(14,23){\makebox(0,0)[cc]{$\bullet$}}
\put(14,19){\makebox(0,0)[cc]{$\bullet$}}
\put(25,23){\makebox(0,0)[cc]{$\bullet$}}
\put(25,19){\makebox(0,0)[cc]{$\bullet$}}
\put(14.5,20){\makebox(0,0)[lc]{$e$}}
\put(25.5,22){\makebox(0,0)[lc]{$e$}}
\end{picture}
\caption{\label{Jitulka}%
Edges allowed in a $\frac12$graph.}
\end{figure}
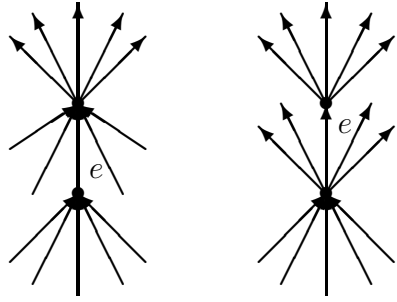

Let an {\em $(m,n)$-$\frac12$graph\/} be a connected simply-connected
directed $(m,n)$-graph whose each edge $e$ has the following
property: either $e$ is the unique outgoing edge of its initial vertex
or $e$ is the unique incoming edge of its terminal vertex, see
Figure~\ref{Jitulka}.
An example of an $(m,n)$-$\frac12$graph is given in Figure~\ref{Jit}.
Let $\Grh mn$ be the category of $(m,n)$-$\frac12$graphs and
their isomorphisms. Define a triple $\freehPROP : \sigmabimod
\to \sigmabimod$ by
\begin{figure}
\unitlength 3.5mm
\begin{picture}(23,17)(7,7)
\thicklines
\put(19,16){\vector(3,4){3}}
\put(19,16){\vector(0,1){5}}
\put(18,16){\vector(-1,2){2}}
\put(16,20){\vector(1,1){2}}
\put(16,20){\vector(0,1){2}}
\put(16,20){\vector(-1,1){2}}
\put(19,13){\vector(0,1){3}}
\put(19,13){\vector(-1,3){1}}
\put(22,20){\vector(-1,2){1}}
\put(22,20){\vector(1,2){1}}
\put(16,10){\vector(1,1){3}}
\put(22,10){\vector(-1,1){3}}
\put(19,10){\vector(0,1){3}}
\put(15,8){\vector(1,2){1}}
\put(17,8){\vector(-1,2){1}}
\put(19,8){\vector(0,1){2}}
\put(19,13){\vector(1,1){3}}
\put(21,8){\vector(1,2){1}}
\put(23,8){\vector(-1,2){1}}
\put(18,16){\vector(0,1){2}}
\put(14,22){\makebox(0,0)[cc]{$\bullet$}}
\put(16,20){\makebox(0,0)[cc]{$\bullet$}}
\put(19,21){\makebox(0,0)[cc]{$\bullet$}}
\put(22,20){\makebox(0,0)[cc]{$\bullet$}}
\put(18,16){\makebox(0,0)[cc]{$\bullet$}}
\put(19,16){\makebox(0,0)[cc]{$\bullet$}}
\put(22,16){\makebox(0,0)[cc]{$\bullet$}}
\put(18,18){\makebox(0,0)[cc]{$\bullet$}}
\put(19,13){\makebox(0,0)[cc]{$\bullet$}}
\put(16,10){\makebox(0,0)[cc]{$\bullet$}}
\put(19,10){\makebox(0,0)[cc]{$\bullet$}}
\put(22,10){\makebox(0,0)[cc]{$\bullet$}}
\put(21,8){\makebox(0,0)[cc]{$\bullet$}}
\put(16,23){\makebox(0,0)[cc]{$2$}}
\put(21,23){\makebox(0,0)[cc]{$1$}}
\put(18,23){\makebox(0,0)[cc]{$3$}}
\put(23,23){\makebox(0,0)[cc]{$4$}}
\put(19,7){\makebox(0,0)[cc]{$1$}}
\put(15,7){\makebox(0,0)[cc]{$2$}}
\put(17,7){\makebox(0,0)[cc]{$3$}}
\put(23,7){\makebox(0,0)[cc]{$4$}}
\end{picture}
\caption{\label{Jit}A graph from $\Grh44$.}
\end{figure}
\begin{equation}
\label{Jitka_uz_zase_premysli!}
\freehPROP(E)(m,n) := \colim{{G \in {\Grh mn}}}{E}(G),\
m,n \geq 0.
\end{equation}

\begin{definition}
A \hPROP\ (called a {\em meager \PROP\/} in~\cite{kontsevich:message})
is an algebra over the triple $\freehPROP : \sigmabimod \to
\sigmabimod$.
\end{definition}

A biased
definition of \hPROP{s} can be found in~\cite{kontsevich:message,markl:ba,mv}.
We followed the original
convention of~\cite{kontsevich:message} that \hPROP{s} do not have units;
the unital version of \hPROP{s} can be defined in an
obvious way, compare also the remarks in~\cite{markl:ba}.

\begin{example}
{\rm
$\frac12$bialgebras are algebras over a \hPROP\ which we denote
$\frac12\sfb$. Another example of structures that can be defined over
\hPROP{s} are {\em Lie $\frac12$bialgebras\/} consisting of a
Lie algebra bracket $[-,-] : V \ot V \to V$ and
a Lie diagonal $\delta : V
\to V \ot V$ satisfying one-half of~(\ref{Jitka_Spetkova_napsala}):
\[
\delta[a,b] = 0.
\]
}
\end{example}

Let us denote by
\[
F : \cathPROP \to \catPROP
\]
the left adjoint to the forgetful functor
$\catPROP \stackrel{\Box}\longrightarrow \cathPROP$ from the category
of \PROP{s} to the category of \hPROP{s}.
M.~Kontsevich observed that, in contrast to $F_1 : \catDiop \to
\catPROP$ in~(\ref{jsem_uplne_hotovej}), $F$~is a
{\em polynomial\/} functor, which immediately implies the following
important theorem~\cite{kontsevich:message,mv}.

\begin{theorem}
\label{polynomial}
The functor $F :  \cathPROP \to \catPROP$
is exact.
\end{theorem}

Now one may take a resolution $\frac12\sfr$ of the \hPROP\
$\frac12\sfb$
and put $\frac12\sfR := F(\frac12\sfr)$.
Theorem~\ref{polynomial} guarantees that $\frac12\sfR$ defined in
this way is indeed a resolution of the
\PROP\ $\frac12\sfB$. Let us mention that there are also
other structures invented to study resolutions of the
\PROP\ $\sfB$, as $\frac23$PROP{s} of Shoikhet~\cite{shoikhet:2/3},
matrons of Saneblidze and Umble~\cite{umble-saneblidze:KK}, or
special \PROP{s} considered in~\cite{markl:ba}.

\begin{center}
-- -- -- -- --
\end{center}

\begin{figure}
\[
\def\arraystretch{1.2}
\begin{array}{|c|c|}
\hline
\mbox {\bf Pasting schemes}  & \mbox {\rule{20pt}{0pt}\bf
  corresponding structures}\rule{20pt}{0pt}
\\
\hline\hline
\mbox {rooted trees} & \mbox {non-unital operads}
\\
\mbox {May's trees} & \mbox {non-unital May's operads}
\\
\mbox {extended rooted trees} & \mbox {operads}
\\
\mbox {cyclic trees} & \mbox {non-unital cyclic operads}
\\
\mbox {extended cyclic trees} & \mbox {cyclic operads}
\\
\mbox {stable labeled graphs} & \mbox {modular operads}
\\
\mbox {extended directed graphs} & \mbox {PROPs}
\\
\mbox {extended connected directed graphs} & \mbox {properads}
\\
\mbox {extended connected $1$-connected dir.~graphs} & \mbox {dioperads}
\\
\mbox {$\frac12$graphs} & \mbox {$\frac 12$PROPs}
\raisebox{-10pt}{\rule{0pt}{20pt}}
\\
\hline
\end{array}
\]
\caption{\label{tables}Pasting schemes and the
  structures they define.}
\end{figure}

The constructions reviewed in this section can be organized into the
following chain of inclusions of full subcategories:
\[
\catOper \subset \cathPROP \subset \catDiop \subset \catProper \subset
\catPROP.
\]
The general scheme behind all these constructions is the following.
We start by choosing a sub-groupoid $\subc = \bigsqcup_{m,n \geq 0}
\subc(m,n)$ of $\GR:=   \bigsqcup_{m,n \geq 0} \Gr
mn$ (or a subgroupoid of
$\UGR :=  \bigsqcup_{m,n \geq 0}\UGr mn$ if we want units).
Then we define a functor $\Gamma_{\tt S}
: \sigmabimod \to  \sigmabimod$ by
\[
\Gamma_{\tt S}(E)(m,n) := \colim{{G \in {\subc( m,n)}}}{E}(G),\
m,n \geq 0.
\]
It is easy to see that $\Gamma_{\tt S}$ is a
subtriple of the \PROP\ triple $\freePROP$ if and only if
the following two conditions are satisfied:
\begin{itemize}
\item[(i)]
the groupoid $\subc$ is {\em hereditary\/} in the sense that,
given a graph from
$\subc$ with vertices decorated by graphs from $\subc$, then the graph
obtained by `forgetting the braces' again belongs to $\subc$, and
\item[(ii)]
$\subc$ contains all directed corollas.
\end{itemize}

Hereditarity~(i) is necessary for $\Gamma_{\tt S}$ to be closed under the
triple multiplication of $\Gamma_{\tt P}$ while~(ii) guarantees
that $\Gamma_{\tt S}$ has an unit.
Plainly, all the three choices used above -- $\UGRC$,
$\UGRD$ and $\GRH$ -- satisfy the above assumptions.
Let us mention that one may modify the definition of \PROP{s}  also by {\em
enlarging\/} the category $\Gr mn$, as was done for  {\em wheeled
\PROP{s}\/} in~\cite{merkulov:defquant}.
Pasting schemes and the corresponding structures reviewed in this
article are listed in Figure~\ref{tables}.

\label{refs}
\def\cprime{$'$}

\end{document}